\documentclass[11pt,paper=a4,oneside]{scrartcl}
\usepackage{lmodern}
\usepackage{mathtools}
\usepackage{amssymb}
\usepackage{amsmath}
\usepackage{amsthm}
\usepackage{hyperref}
\usepackage[capitalise]{cleveref}
\usepackage{xcolor}
\usepackage{bbm}
\usepackage{nicematrix}
\usepackage{float}
\usepackage{eqparbox}
\usepackage{xparse}
\usepackage{tikz}
\usepackage{tikz-cd}
\usepackage{subcaption}
\usepackage{pgfplots}
\usepackage[shortlabels]{enumitem}
\usepackage[english]{babel}
\usepackage[style=ddmmyyyy,datesep={.}]{datetime2}

\usepgfplotslibrary{groupplots}
\pgfplotsset{compat=1.18}
\usetikzlibrary{positioning}
\usetikzlibrary{arrows.meta}

\crefname{equation}{}{}
\crefname{enumi}{}{}
\Crefname{enumi}{}{}
\crefname{assumption}{Assumption}{Assumptions}
\Crefname{assumption}{Assumption}{Assumptions}

\newtheorem{theorem}{Theorem}[section]
\newtheorem{lemma}[theorem]{Lemma}

\newtheorem{corollary}[theorem]{Corollary}
\theoremstyle{definition}
\newtheorem{remark}[theorem]{Remark}
\newtheorem{definition}[theorem]{Definition}

\newtheorem{conjecture}[theorem]{Conjecture}

\newtheorem{example}[theorem]{Example}

\newtheorem{problem}[theorem]{Open Problem}
\newtheorem{assumption}[theorem]{Assumption}

\allowdisplaybreaks

\newcommand{\C}{\mathbb{C}}
\newcommand{\F}{\mathbb{F}}
\newcommand{\N}{\mathbb{N}}
\newcommand{\R}{\mathbb{R}}
\newcommand{\Z}{\mathbb{Z}}
\newcommand{\mA}{\mathcal{A}}
\newcommand{\mB}{\mathcal{B}}
\newcommand{\mC}{\mathcal{C}}
\newcommand{\mD}{\mathcal{D}}

\newcommand{\mG}{\mathcal{G}}
\newcommand{\mL}{\mathcal{L}}
\newcommand{\mM}{\mathcal{M}}
\newcommand{\mN}{\mathcal{N}}
\newcommand{\mP}{\mathcal{P}}
\newcommand{\mR}{\mathcal{R}}
\newcommand{\mZ}{\mathcal{Z}}
\renewcommand{\d}{\mathop{}\!\mathrm{d}}
\newcommand{\abs}[1]{\left|#1\right|}
\newcommand{\floor}[1]{\left\lfloor#1\right\rfloor}
\newcommand{\ceil}[1]{\left\lceil#1\right\rceil}

\newcommand{\one}{\mathbbm 1}
\newcommand{\zero}{\mathbf{0}}
\newcommand{\mt}{\underline{\hspace{\widthof{1}}}} 

\newcommand{\concat}{\|}
\newcommand{\origin}[1]{\underline{#1}}

\DeclareMathOperator{\id}{\operatorname{id}}
\DeclareMathOperator{\ord}{\operatorname{ord}}

\DeclareMathOperator{\spn}{\operatorname{span}}
\DeclareMathOperator{\Tr}{\operatorname{Tr}}
\DeclareMathOperator{\supp}{\operatorname{supp}}
\DeclareMathOperator{\inv}{\operatorname{inv}}
\DeclareMathOperator{\ev}{\operatorname{ev}}

\setlist[enumerate,1]{label=(\alph*)}

\newcounter{eqmakeboxcounter}

\newcommand{\w}[2][\theeqmakeboxcounter]{\eqmakebox[#1][l]{\mathstrut $#2$\,}}
\ExplSyntaxOn
\NewDocumentCommand{\ws}{m}
 {
  \tl_map_inline:nn { #1 } { \w { ##1 }}
 }
\ExplSyntaxOff

\newcommand{\term}[2][1pt]{\arraycolsep=#1 \,\,\begin{matrix} #2 \end{matrix}\,\,}
\newcommand{\niceterm}[3][1pt]{\arraycolsep=#1 \,\,\begin{NiceMatrix} #2 \CodeAfter #3\end{NiceMatrix}\,\,}

\newcommand{\verticalline}{\rule[.5ex]{1em}{0.5pt}}
\newcommand{\landscapevector}[1]{\;\arraycolsep=1pt\begin{matrix} \verticalline & #1 & \verticalline \end{matrix}\;}
\newcommand{\cvector}{\landscapevector{c}}
\newcommand{\cq}{c^{(q)}}
\newcommand{\cm}{c^{(m)}}
\newcommand{\cstarvector}{\landscapevector{\cm}}

\newlength{\stretchlen}
\def\splitterm{\^}

\def\realstretch#1{%
    \def\temp{#1}%
    \ifx\temp\splitterm
    \else
    \mathmbox[\stretchlen]{\hss#1\hss}\expandafter\realstretch
\fi}

\newcommand{\wordcircle}[4]{%

\begin{tikzpicture}
\def\r{#1}
\def\n{#2}

\foreach \letter [count=\i] in {#3} {
    \pgfmathsetmacro{\angle}{90-360*(\i-1)/\n}
    \fill (\angle:\r) circle (2pt);
    \node at (\angle:{\r+0.3}) {\letter};
}

\foreach \i in {#4} {
    \pgfmathsetmacro{\sepangle}{90-360*(\i-0.5)/\n}
    \node[rotate=\sepangle-90,scale=1.4] at (\sepangle:\r) {$\|$};
}

\draw (0,0) circle (\r);
\end{tikzpicture}
}

\title{Shift-invariant functions with the polynomial composition property}
\author{Björn Kriepke}

\begin{document}

\maketitle



\begin{abstract}
    
    Shift-invariant functions on $\F_2^n$ play an important role in the design of cryptographic primitives, with the nonlinear map $\chi$ used in SHA-3 being a prominent example. The iterates of $\chi$ have an interesting algebraic behavior: they can be described using a family of functions whose composition behaves similar to polynomial multiplication. This structure leads to a connection between shift-invariant permutations and polynomial quotient rings.
    
    In this thesis, we investigate the general phenomenon underlying these constructions. We introduce the notion of a quasi-conserved landscape and a corresponding polynomial composition property. Our main result characterizes quasi-conserved landscapes as precisely those whose associated families of functions satisfy the polynomial composition property. We further classify the polynomial rings arising from this construction. This provides a general framework for the algebraic structures previously observed for $\chi$ and several of its generalizations.
    
    Finally, we determine the complete cycle structure of the functions $\chi^{(m)}$. We derive closed forms for generating functions encoding the number of cycles of each possible length by translating the problem into combinatorics on words. This description also yields asymptotic results for the number of cycles.
\end{abstract}

\newpage
\tableofcontents

\section{Introduction}

Boolean functions are an important topic in discrete mathematics, with connections to finite geometry, coding theory and cryptography. A vectorial Boolean function is a map of the form
\begin{equation*}
    f:\F_2^n\to\F_2^m
\end{equation*}
for some positive integers $n,m$. Of particular interest is the case $m=n$. In cryptography, vectorial Boolean functions are commonly used as so-called S-boxes, which play a central role in the security of many cryptographic primitives. Consequently, such functions are required to satisfy various desirable properties, including efficient implementation. One property that can contribute to efficient hardware implementations is \textit{shift-invariance}.

The most prominent shift-invariant function is $\chi:\F_2^n\to\F_2^n$, which was introduced by Daemen \cite{daemen1995cipher}. Its component functions are defined by
\begin{equation}\label{eq:intro_chi_definition}
    \chi(x)_i = x_i + (1+x_{i+1}) x_{i+2}
\end{equation}
for all $i\in\{0,\ldots,n-1\}$. The function $\chi$ is used in many cryptographic primitives, for example SHA-$3$ \cite{nistSHA3StandardPermutationBased2015}, \textsc{Ascon} \cite{dobraunig2021ascon}, Rasta \cite{dobraunig2018rasta} and Subterranean \cite{daemen2020subterranean}. Daemen \cite{daemen1995cipher} showed that $\chi$ is invertible on $\F_2^n$ if and only if $n$ is odd, and an explicit formula for its inverse was found in \cite{liu2022inverse}. More recently, \cite{schoone2024algebraic} studied various algebraic properties of $\chi$, including its univariate polynomial representations over $\F_{2^n}$. It conjectured that the map $\tilde{\chi}:\F_q^n\to\F_q^n$, defined by the same equation as \cref{eq:intro_chi_definition}, is never a permutation for $q>2$, except for trivial cases. This conjecture was further investigated in \cite{otal2023conjecture} before being solved completely in \cite{graner2025bijectivity}. Several generalizations of $\chi$ have also been considered in \cite{andreoli2025generalizations,belkheyar2025chilow,feng2026generalization}.

The topic of this thesis generalizes and unifies our studies of the function $\chi$ from \cite{kriepke2024algebraic,kriepke2026siblings}. In \cite{kriepke2024algebraic} we showed that, for odd $n$, the iterates of $\chi$ can be written as linear combinations of a family of functions $\gamma_{2k}:\F_2^n\to\F_2^n$ for $k\geq 0$. This is an unexpected result, since usually the composition of nonlinear functions does not have any additive structure. We proved that the functions $\gamma_{2k}$ satisfy the identity 
\begin{equation}\label{eq:intro_pcp_gammas}
    \gamma_{2k}\circ \left(\gamma_0 + \sum_{i=1}^s \alpha_i \gamma_{2i}\right) = \sum_{i=0}^s \alpha_i \gamma_{2i+2k}
\end{equation}
for all $\alpha_0,\ldots,\alpha_s\in\F_2$ with $\alpha_0=1$. Motivated by this, we considered the set 
\begin{equation*}
    \mG_n=\gamma_0+\spn\{\gamma_{2j}: j\geq 1\}. 
\end{equation*}
By using \cref{eq:intro_pcp_gammas}, we established that $(\mG_n,\circ)$ is a monoid isomorphic to 
\begin{equation*}
    \mM_n = 1+\spn\{X^j: j\geq 1\} \subseteq \F_2[X]/(X^{(n+1)/2}). 
\end{equation*}
Since $(\mM_n,\cdot)$ is in fact an abelian group, this shows that $\mG_n$ is an abelian group, so that every element $f\in\mG_n$ is a permutation on $\F_2^n$. The work \cite{kriepke2026siblings} then extended this analysis to even $n$. The identity \ref{eq:intro_pcp_gammas} remains valid here, while
\begin{equation*}
    \mG_n=\gamma_0+\spn\{\gamma_{2j}: j\geq 1\}
\end{equation*}
is a monoid isomorphic to 
\begin{equation*}
    \mM_n = 1+\spn\{X^j: j\geq 1\} \subseteq \F_2[X]/(X^n+X^{n/2}). 
\end{equation*}
Although $\mM_n$ is not a group in this case, this still enables a complete classification of all permutations in $\mG_n$. The results of these papers in particular yield constructions for large families of shift-invariant permutations on $\F_2^n$. 

The tools of \cite{kriepke2024algebraic} were subsequently generalized in \cite{lyu2026generalized} to the maps $\chi^{(m)}:\F_2^n\to\F_2^n$ defined for $m\geq 2$ by 
\begin{equation*}
    \chi^{(m)}(x)_i = x_i + (1+x_{i+1}) (1+x_{i+2})\cdots (1+x_{i+m-1})x_{i+m}.
\end{equation*}
It generalized the functions $\gamma_{2k}$ to a family $\theta_{m,i}$ and showed that a similar identity to \cref{eq:intro_pcp_gammas} holds. This also allowed to prove a connection to a polynomial ring $\F_2[X]/(p(X))$ for some polynomial $p(X)\in\{X^{\ceil{n/m}}, X^{n/m}(X^{n/m}+1)\}$. 

The central observation underlying all of these results is the identity \cref{eq:intro_pcp_gammas}, which is fundamental for establishing a connection between composition of functions and multiplication of polynomials. This raises the natural question of which functions satisfy such an identity, which we call polynomial composition property. One of the main achievements of this thesis is a complete answer to this question.

To establish this result, we introduce the notion of quasi-conserved landscapes and show that a landscape is quasi-conserved if and only if its associated family of functions $\gamma_j$, $j\geq0$, satisfies the polynomial composition property. Moreover, we completely classify the polynomial rings arising from this construction. They are always of the form
\begin{equation*}
    \F_2[X]/(X^j), \quad \F_2[X]/(X^j(X^r+1)) \quad \text{or}\quad \F_2[X]
\end{equation*}
for some integers $j,r$ with $j\leq r$, see \cref{thm:pcp_minimal_polynomial_conclusion}. 

The key result, \cref{thm:qc_quasi_conserved_iff_polynomial_composition_property}, was presented at WCC \cite{kriepke2026wcc}. The complete result, including \cref{thm:pcp_minimal_polynomial_conclusion}, was subsequently presented in \cite{kriepke2026pott65}. Shortly after the WCC presentation, the preprint \cite{zhao2026bifix} appeared, independently proving a weaker version of \cref{thm:qc_quasi_conserved_iff_polynomial_composition_property}. More recently, after the technical part of this thesis had been completed, the preprint \cite{lyu2026classshiftinvariant} appeared on arXiv. Its main results seem to be similar to our work. A detailed comparison of the present results with those of \cite{lyu2026classshiftinvariant} is left for future work. 

The thesis is organized as follows: In \cref{sec:prelim} we give a brief overview of shift-invariant functions. We first introduce shift-invariant functions on $\F_2^\Z$ and discuss several equivalent characterizations. We then explain how to recover the finite-dimensional setting. 

In \cref{sec:quasi_conserved_landscapes}, we introduce the concept of quasi-conserved landscapes, which is a key step in proving \cref{thm:qc_quasi_conserved_iff_polynomial_composition_property}. Then we discuss in detail the results of \cite{kriepke2024algebraic,kriepke2026siblings,lyu2026generalized}. We generalize the functions $\gamma_{2j}$ from \cite{kriepke2024algebraic,kriepke2026siblings} and $\theta_{m,i}$ from \cite{lyu2026generalized}, and define the polynomial composition property, which generalizes \cref{eq:intro_pcp_gammas}. We then state the main theorem of this thesis, \cref{thm:qc_quasi_conserved_iff_polynomial_composition_property}, which establishes the equivalence between quasi-conservedness and the polynomial composition property.

We prove \cref{thm:qc_quasi_conserved_iff_polynomial_composition_property} in \cref{sec:right_sided_landscapes} and \cref{sec:two_sided_landscapes}. In \cref{sec:right_sided_landscapes}, we first consider the special case of right-sided landscapes, which is particularly relevant since, for example, both $\chi$ and $\chi^{(m)}$ belong to this class. For this case, we present a more intuitive proof that provides useful insight into the general argument. We then turn to the general case in \cref{sec:two_sided_landscapes}, where the proof is considerably more technical.

In \cref{sec:pcp_studying_pcp}, we study functions satisfying the polynomial composition property. We introduce the notion of the minimal polynomial $m_{c,V}(X)$ of a quasi-conserved landscape $c$ over a vector space $V\in\{\F_2^n,\F_2^\Z\}$ and determine all possible minimal polynomials. In particular, for a quasi-conserved landscape on $\F_2^n$, we show that
\begin{equation*}
    m_{c,\F_2^n}(X) \in \{X^j, X^j(X^r+1)\},
\end{equation*}
where $r=n/\gcd(q,n)$ for some parameter $q$ depending on $c$ and $j\leq r$. This leads to a classification of all shift-invariant permutations in
\begin{equation*}
    \mG = \gamma_0 + \spn\{\gamma_j: j\geq 0\}.
\end{equation*} 
This classification gives constructions of large families of shift-invariant permutations. We also obtain results for the case $V=\F_2^\Z$, where
\begin{equation*}
    m_{c,\F_2^\Z}(X)\in \{0, X^j\}
\end{equation*}
and show that all functions $f\in\mG$ are surjective as functions $\F_2^\Z\to\F_2^\Z$. 

Finally, we study the cycle structure of $\chi^{(m)}$. Initial results concerning this cycle structure were obtained in \cite{kriepke2024algebraic,lyu2026generalized,schoone2024state}. In \cref{sec:cycle_structure}, we determine the cycle structure completely by deriving closed forms for generating functions that encode the number of cycles of each possible length. We obtain these generating functions by translating the problem into one in combinatorics on words. This description also allows us to derive asymptotic results for the cycle counts.

\section{Preliminaries on shift-invariant functions}\label{sec:prelim}

In this section we cover some preliminaries on shift-invariant functions. We first introduce some notation. Let $\Z$ be the sets of integers and $\N_0, \N$ be the sets of non-negative and positive integers, respectively. By $\F_2=\{0,1\}$ we denote the field with two elements. While some parts of this section also apply to other fields, the identity 
\begin{equation}\label{eq:prelim_x_times_x_plus_1_equals_zero}
    x\cdot(x+1)=0
\end{equation}
for $x\in\F_2$ is crucial for our work, so that we focus on $\F_2$.

We aim to study shift-invariant functions over $\F_2^n$. It will be more convenient to work on the space $\F_2^\Z$, i.e. the set of functions $\Z\to\F_2$. This is standard practice, see for example \cite{daemen1995cipher}. We view an element $x\in\F_2^\Z$ as an infinite row vector in both directions, i.e.
\begin{equation*}
    x = (\ldots, x_{-2},x_{-1},\origin{x_0},x_1,x_2,\ldots),
\end{equation*}
where we underline the origin. Depending on the context we may also write elements of $\F_2^\Z$ as column vectors. We later identify $\F_2^n$ with subspaces $\mP_n\leq \F_2^\Z$ and in that way recover the finite dimensional setting. 

The main object of study in this thesis are shift-invariant functions on $\F_2^\Z$, which we now define. 

\begin{definition}
    The linear left-shift operator $S:\F_2^\Z\to\F_2^\Z$ is defined by $x\mapsto y$ with $y_i=x_{i+1}$ for all $i\in\Z$. We write $Sx$ instead of $S(x)$ for $x\in\F_2^\Z$ and $Sf$ instead of $S\circ f$ for any function $f:\F_2^\Z\to\F_2^\Z$. 
\end{definition}

\begin{definition}
    A function $f:\F_2^\Z\to\F_2^\Z$ is called shift-invariant if it commutes with the left-shift $S$, i.e. 
    \begin{equation*}
        f(Sx)=Sf(x)
    \end{equation*}
    for all $x\in\F_2^\Z$. 
\end{definition}
Observe that since $S$ is invertible it follows that $f$ is shift-invariant if and only if $f(S^ix)=S^if(x)$ for all $x\in\F_2^\Z$ and $i\in\Z$. There are several equivalent characterizations of shift-invariant functions. We shortly discuss them.

\subsection{Coordinate functions}\label{subsec:prelim_coordinate_functions}

Let $f:\F_2^\Z\to\F_2^\Z$ be any function. For any coordinate $i\in\Z$ we denote by $\pi_i$ the projection $\pi_i:\F_2^\Z\to\F_2$ given by $x\mapsto x_i$. Then 
\begin{equation*}
    f_i=\pi_i\circ f: \F_2^\Z\to\F_2, \qquad x\mapsto f(x)_i
\end{equation*}
is called the $i$-th coordinate function of $f$. The following lemma classifies shift-invariant functions using their coordinate functions:

\begin{lemma}\label{lem:prelim_shift_invariant_equiv_coordinate_functions}
    Let $f:\F_2^\Z\to\F_2^\Z$ be any function with coordinate functions $f_i:\F_2^\Z\to\F_2$ for $i\in\Z$. Then the following are equivalent:
    \begin{enumerate}
        \item $f$ is shift-invariant.
        \item $f_i = f_0\circ S^i$ for all $i\in\Z$.
    \end{enumerate}
\end{lemma}

\begin{proof}
    We use the immediate identity $\pi_i\circ S^j = \pi_{i+j}$. 

    $(a)\implies(b)$. Let $f$ be shift-invariant. Then 
    \begin{equation*}
        f_i = \pi_i\circ f = \pi_0\circ S^i \circ f =\pi_0\circ f\circ S^i = f_0\circ S^i
    \end{equation*}
    and the claim follows.

    $(b)\implies(a)$. Suppose now that $f_i=f_0\circ S^i$ for all $i\in\Z$. To show $f\circ S=S\circ f$ it suffices to show that $\pi_i\circ f\circ S = \pi_i\circ S\circ f$ for all $i\in\Z$. Indeed,
    \begin{equation*}
        \pi_i\circ f\circ S = f_i \circ S = f_0\circ S^i\circ S = f_0\circ S^{i+1} = f_{i+1} = \pi_{i+1}\circ f = \pi_i\circ S\circ f.
    \end{equation*}
\end{proof}

As a consequence of \cref{lem:prelim_shift_invariant_equiv_coordinate_functions}, a shift-invariant function is completely determined by the coordinate function $f_0:\F_2^\Z\to \F_2$. In this thesis we only study shift-invariant functions for which there exist non-negative integers $k,m\geq 0$ and a function $g:\F_2^{k+m+1}\to\F_2$ such that 
\begin{equation}\label{eq:prelim_f_0_is_local_map}
    f_0(x) = g(x_{-k},\ldots,x_m).
\end{equation}
Equivalently, $f_0(x)$ only depends on a finite number of entries of $x$. In an abuse of notation we write $f_0$ instead of $g$. The function $f_0$ is then also called a local map and $f$ is called a cellular automaton.

\begin{example}
    Let $f_0:\F_2^\Z\to\F_2$ be given by $f_0(x)=1$ if and only if the sequence $(x_2,x_4,x_6,x_8,\ldots)$ contains finitely many entries equal to $1$. This coordinate function induces a shift-invariant function $f:\F_2^\Z\to\F_2^\Z$, however $f_0$ is not a local map of the form \cref{eq:prelim_f_0_is_local_map}. We are not interested in shift-invariant functions like this.
\end{example}

Any function $f_0:\F_2^{k+m+1}\to\F_2$ can be written as a multivariate polynomial $f_0\in\F_2[x_{-k},\ldots,x_m]$. As we are working over $\F_2$ and for any $y\in\F_2$ it holds that $y^2=y$, we may further assume that every exponent of each variable $x_i$ in $f_0$ is at most $1$. In other words, we consider $f_0$ as an element of the ring $\F_2[x_{-k},\ldots,x_m]/(x_{-k}^2+x_{-k},\ldots,x_m^2+x_m)$. 

\begin{example}
    The function $\chi$ is given by $f_0(x)=x_0+(x_1+1)x_2$. We then have $f_i(x)=f_0(S^i x) = x_i+(x_{i+1}+1)x_{i+2}$ for $i\in\Z$.
\end{example}

\subsection{Vectorial notation}\label{subsec:prelim_vector_notation}

In some situations it is more convenient to work with other representations of shift-invariant functions, instead of coordinate functions. The following notation was introduced in \cite{graner2025bijectivity}. For vectors $x,y\in\F_2^\Z$ we write $z=x\odot y$ for the element $z\in\F_2^\Z$ with $z_i=x_iy_i$ for all $i\in\Z$. In other words, $z$ is the componentwise product of $x$ and $y$, which is also called the Hadamard product of $x$ and $y$. It is clear that the Hadamard product is commutative, associative and distributive over addition. We extend this definition to functions $f,g:\F_2^\Z\to\F_2^\Z$ by $(f\odot g)(x)=f(x)\odot g(x)$ for $x\in\F_2^\Z$. Observe that $S(x\odot y)=Sx\odot Sy$ for all $x,y\in\F_2^\Z$. We denote the element $(\ldots,1,\origin1,1,\ldots)\in\F_2^\Z$ by $\one$. Then $\one$ is the identity with respect to $\odot$, i.e. $x\odot \one=\one\odot x=x$ for all $x\in\F_2^\Z$.

Let $f$ be a shift-invariant function with coordinate function $f_0=f_0(x_{-k},\ldots,x_m)\in\F_2[x_{-k},\ldots,x_m]/(x_{-k}^2+x_{-k},\ldots,x_m^2+x_m)$. We can then write
\begin{equation*}
    f_0(x_{-k},\ldots,x_m) = \sum_{I\subseteq \{-k,\ldots,m\}} a_I \prod_{j\in I} x_j
\end{equation*}
with some constants $a_I\in\F_2$ for every subset $I\subseteq \{-k,\ldots,m\}$. This is called the algebraic normal form (ANF) of $f_0$. Then $f=f_0(S^{-k},\ldots,S^m)$, by which we mean that
\begin{equation}\label{eq:prelim_f_is_f_0_of_shifts}
    f = \sum_{I\subseteq\{-k,\ldots,m\}} a_I \bigodot_{j\in I} S^j.
\end{equation}
It suffices to show that the $i$-th coordinate functions $f_i$ and $f_0(S^{-k},\ldots,S^m)_i$ are the same for all $i\in\Z$. Indeed, for every $x\in\F_2^\Z$ it holds that
\begin{align*}
    f_i(x) &= f_0(x_{i-k},\ldots,x_{i+m}) \\
        &= \sum_{I\subseteq \{-k,\ldots,m\}} a_I \prod_{j\in I} x_{i+j} \\
        &= \sum_{I\subseteq \{-k,\ldots,m\}} a_I \prod_{j\in I} (S^jx)_i \\
        &= \sum_{I\subseteq \{-k,\ldots,m\}} a_I \left(\bigodot_{j\in I} S^jx\right)_i \\
        &= f_0(S^{-k}x,\ldots,S^mx)_i.
\end{align*}

\begin{example}
    Let $f=\chi$. Then $f_0(x) = x_0 + (1+x_1)x_2$ and
    \begin{align*}
        \chi(x) &= 
            \begin{pmatrix}
                \vdots \\
                x_{-1} + (1+x_0)x_1 \\
                \origin{x_0    + (1+x_1)x_2} \\
                x_1    + (1+x_2)x_3 \\
                \vdots
            \end{pmatrix} \\
            &= 
            \begin{pmatrix}
                \vdots \\ x_{-1} \\ \origin{x_0} \\ x_1 \\ \vdots
            \end{pmatrix}
            +
            \left[
                \begin{pmatrix}
                    \vdots \\ 1 \\ \origin{1} \\ 1 \\ \vdots 
                \end{pmatrix}
                + 
                \begin{pmatrix}
                    \vdots \\ x_0 \\ \origin{x_1} \\ x_2 \\ \vdots 
                \end{pmatrix}
            \right] 
            \odot 
            \begin{pmatrix}
                \vdots \\ x_1 \\ \origin{x_2} \\ x_3 \\ \vdots 
            \end{pmatrix} \\
            &= S^0(x) + (\one+S^1(x))\odot S^2(x) \\
            &= (S^0+(\one+S^1)\odot S^2)(x),
    \end{align*}
    so that
    \begin{equation*}
        \chi = S^0 + (\one+S^1)\odot S^2.
    \end{equation*}
\end{example}

\subsection{Complementing landscapes} \label{subsec:prelim_landscapes}

The third view that we discuss are so-called (complementing) landscapes. Let $f$ be a shift-invariant function, let $x\in\F_2^\Z$ and set $y=f(x)$. As $x$ and $y$ are elements of $\F_2^\Z$, each entry is a bit in $\{0,1\}$. Hence, to find $y$ it is enough to know when a coordinate $x_i$ gets flipped under application of $f$, i.e. we need to know whether $y_i=x_i$ or $y_i=x_i+1$ for each $i\in\Z$. We focus on shift-invariant functions for which $y_i=x_i+1$ occurs if and only if a certain pattern occurs in the neighborhood of $x_i$. Such a pattern is then called a complementing landscape. 

\begin{example}
    Consider $\chi$, fix $x\in\F_2^\Z$ and set $y=\chi(x)$. Then $y_i = x_i+(x_{i+1}+1)x_{i+2}$ for all $i\in\Z$. In particular, $y_i = x_i+1$ if and only if $(x_{i+1},x_{i+2})=(0,1)$. In other words, the bit $x_i$ gets flipped by $\chi$ if and only if it is followed by the pattern $01$. We then say that $\chi$ has the complementing landscape $*01$. 
\end{example}

\begin{definition}\label{def:prelim_landscapes}
    Let $k,m\geq 0$ be non-negative integers. A landscape $c$ is an element $c\in\{0,1,\mt\}^{k+m+1}$ where $\mt$ is a placeholder symbol. We write $c$ as a string
    \begin{equation*}
        c = a_{-k}a_{-k+1}\ldots a_{-1} \hat{a_0} a_1 a_2\ldots a_m.
    \end{equation*}
    The position $0$ is highlighted, if $a_0=\mt$ then we write $*$, otherwise we write $\hat 0$ or $\hat 1$, depending on $a_0=0$ or $a_0=1$, respectively. 
    
    Landscapes of the form $*a_1a_2\ldots a_m$ and $a_{-k} \ldots a_{-1} *$ are called right-sided and left-sided, respectively. Other landscapes are called two-sided.
    
    To a landscape $c=a_{-k}\ldots a_m$ we associate the indicator function of $c$, which in an abuse of notation we also call $c:\F_2^\Z\to\F_2^\Z$. It is defined by $c:x\mapsto y$ where $y_i=1$ if and only if 
    \begin{equation}\label{eq:prelim_x_matches_c_at_i}
        (x_{i-k},\ldots,x_{i-1},x_i,x_{i+1},\ldots,x_{i+m}) = (a_{-k},\ldots, a_{-1}, \hat{a_0}, a_1, \ldots, a_m),
    \end{equation}
    where a condition of the form $x_j=\mt$ means that $x_j\in \F_2$ can be arbitrary. 
    If \cref{eq:prelim_x_matches_c_at_i} holds we say that $x$ matches the landscape $c$ at position $i$.

    For a landscape $c=a_{-k}\ldots a_m$ we define its support $\supp(c)$ by 
    \begin{equation*}
        \supp(c) = \{j\in\{-k,\ldots,m\} : a_j \neq \mt\}.
    \end{equation*}
\end{definition}

\begin{remark}
    In the literature it is often written that $\chi=*01$, see for example \cite{daemen1995cipher}. With our notation $\chi=\id+*01$, which is more convenient for our studies.
\end{remark}

\begin{example}
    Let $c=1*10$. Then $\supp(c)=\{-1,1,2\}$. For the following $x\in\F_2^\Z$ we have:
    \begin{align*}
        \term[2pt]{
            x    & = & (\ldots & \origin0 & 1 & 0 & 1 & 0 & 0 & 1 & 1 & 1 & 0 & 0 & 1 & 1 & 0 & 1 & 1 & 1 & 0 & 1 & 0 & 1 & 0 & 0 & \ldots), \\
            c(x) & = & (\ldots & \origin. & . & 1 & . & . & . & . & 1 & . & . & . & . & . & . & . & 1 & . & 1 & . & 1 & . & . & . & \ldots).
        }
    \end{align*}
    For readability we have replaced every $0$ in $c(x)$ by a dot. We see that $x$ matches $c$ at positions $2,7,15,17$ and $19$.
\end{example}

From a landscape $c=a_{-k}\ldots a_m$ it is easy to obtain the coordinate function $c_0:\F_2^{k+m+1}\to\F_2$ and the vectorial notation of $c:\F_2^\Z\to\F_2^\Z$. Then
\begin{equation}\label{eq:prelim_coordinate_function_from_landscape}
    c_0(x) = \varepsilon_{-k}\cdot \varepsilon_{-k+1}\cdots \varepsilon_m
\end{equation}
where
\begin{equation*}
    \varepsilon_\ell = \begin{cases}
        x_\ell         & \text{if } a_\ell = 1, \\
        1+x_\ell       & \text{if } a_\ell = 0, \\
        1           & \text{if } a_\ell = \mt.
    \end{cases}
\end{equation*}
Indeed, observe that $c_0(x)=1$ if and only if \cref{eq:prelim_x_matches_c_at_i} holds. Furthermore, by \cref{eq:prelim_f_is_f_0_of_shifts},
\begin{equation}\label{eq:prelim_vector_notation_from_landscape}
    c = A_{-k}\odot A_{-k+1} \odot \cdots \odot A_m
\end{equation}
where $A_\ell:\F_2^\Z\to\F_2^\Z$ are given by
\begin{equation*}
    A_\ell = \begin{cases}
        S^\ell         & \text{if } a_\ell = 1, \\
        \one + S^\ell  & \text{if } a_\ell = 0, \\
        \one           & \text{if } a_\ell = \mt.
    \end{cases}
\end{equation*}
For example, for $c=*01$,
\begin{equation*}
    c_0(x) = (1+x_1)x_2
\end{equation*}
and
\begin{equation*}
    c = \one\odot(\one + S^1)\odot S^2=(\one+S^1)\odot S^2.
\end{equation*}

An advantage of the landscape notation is that it is easy to compute the Hadamard product of two shift-invariant functions given in landscape notation. As an example, let 
\begin{align*}
    c &= 0*\mt10 = (\one+S^{-1})\odot S^2\odot(\one+S^3), \\
    d &= 10*01 = S^{-2}\odot (\one+S^{-1})\odot (\one + S^1)\odot S^2. 
\end{align*}
Then
\begin{align*}
    c\odot d &= \left[(\one+S^{-1})\odot S^2\odot(\one+S^3)\right] \odot \left[S^{-2}\odot (\one+S^{-1})\odot (\one + S^1)\odot S^2\right] \\
    &= S^{-2} \odot (\one+S^{-1})\odot (\one + S^1)\odot S^2 \odot (\one+S^3) \\
    &= 10*010.
\end{align*}
We write this product in shorter form as 
\begin{equation}\label{eq:prelim_diagram_computation_of_product_of_landscapes}
    c\odot d = 
    \term[2pt]{
          & 0 & * & \mt & 1 & 0 \\
        1 & 0 & * &   0 & 1
    } 
    =
    \term[2pt]{
        \mt & 0 & * & \mt & 1 & 0 \\
        1   & 0 & * &   0 & 1 & \mt 
    } 
    = 10*010,
\end{equation}
where in the first row is the landscape $c$ and in the second row is the landscape $d$, which are aligned at the highlighted positions. The result of the product can then easily read off by combining both rows. More explicitly, in each column we apply the following rule to the pair of entries:
\begin{align*}
    (0,0) &\mapsto 0, \\
    (1,1) &\mapsto 1, \\
    (z,\mt), (\mt,z) &\mapsto z,
\end{align*}
for $z\in\{0,1,\mt\}$ and any empty position is interpreted as $\mt$. Note that in this example there was no position where one row had a $0$ and the other row had a $1$. Such a situation we call a conflict. In that case the product would be the zero function. For example, let 
\begin{align*}
    c &= 0*10=(\one+S^{-1})\odot S^1\odot(\one+S^2), \\
    d &= 10*00=S^{-2}\odot (\one+S^{-1})\odot (\one + S^1)\odot (\one+S^2). 
\end{align*}
Then
\begin{align*}
    c\odot d &= \left[(\one+S^{-1})\odot S^1\odot(\one+S^2)\right] \odot \left[S^{-2}\odot (\one+S^{-1})\odot (\one + S^1)\odot (\one+S^2)\right] \\
    &= \left[S^1\odot(\one+S^1)\right]\odot(\ldots) \\
    &= \zero
\end{align*}
because
\begin{equation}\label{eq:prelim_S_odot_S_plus_one_equals_zero}
    S^i \odot (\one + S^i) = \zero
\end{equation}
for all $i\in\Z$ due to \cref{eq:prelim_x_times_x_plus_1_equals_zero}. 
Visually, 
\begin{equation*}
    c\odot d = 
    \term[3pt]{
          & 0 & * & 1 & 0 \\
        1 & 0 & * & 0 & 0
    }
    = \zero,
\end{equation*}
and we see that in the fourth column there is a $1$ in the first row and a $0$ in the second row, so the product is the zero function. 

Furthermore, we can also conveniently represent the composition of the shift $S$ with a function given by a landscape. Indeed, if $c=a_{-k}\ldots a_m$, then the landscape of $Sc=S\circ c$ is simply $c$ shifted one position to the right, or equivalently, shifting the highlighted position of $c$ to the left. For example, for $c=0*10=(\one+S^{-1})\odot S^1\odot(\one+S^2)$ we have
\begin{equation*}
    Sc = S\left((\one+S^{-1})\odot S^1\odot(\one+S^2)\right) = (\one+S^0)\odot S^2\odot (\one+S^3) = \hat{0}\mt 10.
\end{equation*}

\subsection{Shift-invariant functions in finite dimension}\label{subsec:prelim_finite_dim}

While for a large part of this thesis we consider shift-invariant functions as functions $\F_2^\Z\to\F_2^\Z$, for practical purposes like cryptography we need functions defined on finitely many bits, i.e. functions $\F_2^n\to\F_2^n$. 

To this end we define a subset $\mP_n\subseteq \F_2^\Z$ with $\mP_n\cong \F_2^n$ and then identify $f:\mP_n\to\mP_n$ with a function $\F_2^n\to\F_2^n$. More precisely, for $n\in\N$ we define the set
\begin{equation*}
    \mP_n = \{x\in\F_2^\Z: x_i = x_{i+n} \text{ for all } i\in\Z\} = \{x\in\F_2^\Z: x = S^nx\}
\end{equation*}
of the elements $x\in\F_2^\Z$ with period $n$. It is clear that $\mP_n$ is a subspace of $\F_2^\Z$ isomorphic to $\F_2^n$. Indeed, the map $\theta:\F_2^n\to\mP_n$ given by 
\begin{equation*}
    (x_0,x_1,\ldots,x_{n-1})\mapsto (\ldots, x_{n-1},\origin{x_0},x_1,\ldots,x_{n-1},x_0,x_1,\ldots,x_{n-1},x_0,\ldots)
\end{equation*}
is an isomorphism. Let $f:\mP_n\to\mP_n$ be any function. Then $f$ induces a function $f^\theta :\F_2^n\to\F_2^n$ by 
\begin{equation*}
    f^\theta = \theta^{-1} \circ f \circ \theta. 
\end{equation*}
The functions $f$ and $f^\theta$ are called conjugate and if the context is clear we write $f$ instead of $f^\theta$. Observe that for functions $f,g:\mP_n\to\mP_n$ it holds that $(f\circ g)^\theta = f^\theta \circ g^\theta$. 

The following lemma shows that $f:\F_2^\Z\to\F_2^\Z$ restricts to a map $f:\mP_n\to\mP_n$ for any $n\in\N$ if $f$ is shift-invariant.

\begin{lemma}\label{lem:prelim_shift_invariant_implies_periodic_stays_periodic}
    Let $f:\F_2^\Z\to\F_2^\Z$ be shift-invariant. Then $f(\mP_n)\subseteq \mP_n$ for any $n\in\N$.
\end{lemma}

\begin{proof}
    Let $x\in \mP_n$ arbitrary, i.e. $S^nx=x$. Set $y=f(x)$. Then 
    \begin{equation*}
        S^ny = S^nf(x) = f(S^nx) = f(x) = y
    \end{equation*}
    and $y\in \mP_n$.
\end{proof}

The left-shift $S:\F_2^\Z\to\F_2^\Z$ is clearly shift-invariant and restricts to a map $S:\mP_n\to\mP_n$, so that we obtain the induced map $S^\theta:\F_2^n\to\F_2^n$. Let $x=(x_0,\ldots,x_{n-1})\in\F_2^n$. Then
\begin{align*}
    S^\theta (x_0,x_1\ldots,x_{n-1}) 
        &= \theta^{-1}(S(\ldots, x_{n-1},\origin{x_0},x_1,\ldots,x_{n-1},x_0,x_1,\ldots,x_{n-1},x_0,\ldots)) \\
        &= \theta^{-1}(\ldots, x_0, \origin{x_1}, x_2, \ldots, x_0, x_1, x_2, \ldots, x_0, x_1, \ldots) \\
        &= (x_1,x_2,\ldots,x_{n-1},x_0),
\end{align*}
and $S^\theta$ is the cyclic left-shift on $\F_2^n$. Now let $f:\F_2^\Z\to\F_2^\Z$ be any shift-invariant function. By \cref{lem:prelim_shift_invariant_implies_periodic_stays_periodic} the map $f^\theta:\F_2^n\to\F_2^n$ is well-defined. Then $f^\theta$ commutes with $S^\theta$, i.e. $f^\theta \circ S^\theta = S^\theta \circ f^\theta$. Indeed, 
\begin{equation*}
    f^\theta \circ S^\theta = (f\circ S)^\theta = (S\circ f)^\theta = S^\theta \circ f^\theta.
\end{equation*}
This is the standard definition of a shift-invariant function $f:\F_2^n\to\F_2^n$ in the literature, see for example \cite{haugland2025shift,haugland2026new,kriepke2024algebraic,kriepke2026siblings,liu2026finding,lyu2026generalized,omland2022permutation,schoone2024algebraic}.

\begin{definition}
    A function $f:\F_2^n\to\F_2^n$ is called shift-invariant if it commutes with the cyclic left-shift $S$, i.e. 
    \begin{equation*}
        f(Sx)=Sf(x)
    \end{equation*}
    for all $x\in\F_2^n$. 
\end{definition}

Let now $f_0\in\F_2[x_{-k},\ldots,x_m]/(x_{-k}^2+x_{-k},\ldots,x_m^2+x_m)$ be a coordinate function and set $f:\F_2^\Z\to\F_2^\Z$ as the corresponding shift-invariant function. It then holds that
\begin{equation*}
    f^\theta(x)_i = f_0(x_{i-k},\ldots,x_{i+m})
\end{equation*}
where all indices of $x=(x_0,\ldots,x_{n-1})\in\F_2^n$ are computed modulo $n$. Given a landscape $c=a_{-k}\ldots a_m$ we obtain a function $c^\theta:\F_2^n\to\F_2^n$ defined by $c^\theta:x\mapsto y$ where $y_i=1$ if and only if 
\begin{equation*}
    (x_{i-k},\ldots,x_{i-1},x_i,x_{i+1},\ldots,x_{i+m}) = (a_{-k},\ldots, a_{-1}, \hat{a_0}, a_1, \ldots, a_m),
\end{equation*}
where an equality of the form $x_j=\mt$ means that $x_j\in \F_2$ can be arbitrary. Again, the indices here are computed modulo $n$. Similarly, the support of $c^\theta$ also is computed modulo $n$. Note that it is possible that $c^\theta=\zero$ on $\F_2^n$, in which case $\supp(c^\theta)=\emptyset$ is the empty set.

\begin{example}
    \begin{enumerate}
        \item Let $c=*01$ and $x=(x_0,\ldots,x_5)=(0,1,0,1,1,1)\in\F_2^6$. Then 
        \begin{equation*}
            c^\theta(x) = (0, 1, 0, 0, 0, 1),
        \end{equation*}
        because $x$ matches $c$ at positions $1$ and $5$. 

        \item Let $c=*0\mt\,\mt1 = (\one+S)\odot S^4$ and $n=3$. Then $c^\theta = (\one+S)\odot S = \zero$ on $\F_2^3$.
    \end{enumerate}
\end{example}

Given a shift-invariant function $f:\F_2^\Z\to\F_2^\Z$, it is a natural question to ask for which $n\in\N$ the function $f:\F_2^n\to\F_2^n$ is a permutation. If $f$ is a permutation on $\F_2^\Z$ then the answer is clear.

\begin{lemma}\label{lem:prelim_f_permutation_F_2_Z_implies_permutation_on_F_2_n}
    Let $f:\F_2^\Z\to\F_2^\Z$ be a shift-invariant permutation. Then $f:\F_2^n\to\F_2^n$ is a permutation for all $n\in\N$.
\end{lemma}

\begin{proof}
    If $f:\F_2^\Z\to\F_2^\Z$ is a permutation, then it is injective, so $f:\mP_n\to\mP_n$ is also injective. But $\mP_n$ is finite. The claim follows.
\end{proof}

The converse of \cref{lem:prelim_f_permutation_F_2_Z_implies_permutation_on_F_2_n} is also true under the additional assumption that $f:\F_2^\Z\to\F_2^\Z$ is a cellular automaton, i.e. $f_0$ is a local map of the form \cref{eq:prelim_f_0_is_local_map}. We do not prove this fact here and instead refer to the survey \cite{kari2005CA} and the references therein. 

However, many interesting functions for cryptography are not invertible for all $n\in\N$. In the paper \cite{omland2022permutation} the set
\begin{equation*}
    \inv(f) = \{n\in\N: f\text{ is a permutation on }\F_2^n\}
\end{equation*}
was introduced. 
The following lemma shows that $\inv(f)$ is closed under taking divisors.

\begin{lemma}[{c.f. \cite[Prop 6.1]{daemen1995cipher}}]\label{lem:prelim_inv_f_is_lower_set}
    Let $f:\F_2^n\to\F_2^n$ be a permutation and let $d\mid n$ be a divisor of $n$. Then $f:\F_2^d\to\F_2^d$ is also a permutation.
\end{lemma}

\begin{proof}
    View $f$ as a function $f:\mP_n\to\mP_n$. As $d\mid n$ it follows that $\mP_d\subseteq \mP_n$. As $f$ is a permutation on $\mP_n$, it is injective on $\mP_n$ and hence also injective on $\mP_d$. Therefore, $f:\mP_d\to\mP_d$ is a permutation.
\end{proof}

\begin{remark}\label{rem:prelim_upper_and_lower_sets}
    Let $(X,\preceq)$ be a partially ordered set. A non-empty subset $L\subseteq X$ is called a lower set of $X$ if whenever $x,y\in X$ with $x\in L$ and $y\preceq x$, then $y\in L$. In other words, $L$ is downward-closed. Dually, a non-empty subset $U\subseteq X$ is called an upper set of $X$ if whenever $x,y\in X$ with $x\in U$ and $x\preceq y$, then $y\in U$, i.e. $U$ is upward-closed. Given a lower set $L$, the complement $X\setminus L$ is an upper set and vice versa.

    Given some subset $A\subseteq X$, we define
    \begin{align*}
        \underline{A} &= \{x\in X: x\preceq a \text{ for some } a\in A\} \\
        \overline{A}  &= \{x\in X: a\preceq x \text{ for some } a\in A\}.
    \end{align*}
    The sets $\underline A, \overline A$ are called the lower and upper sets generated by $A$ and $A$ is called a generating set of $\underline A, \overline A$. A generating set of a lower/upper set always exists.
    However, it is generally not unique.  Of particular interest are inclusion-minimal generating sets, if they exist.
\end{remark}

With the terminology from \cref{rem:prelim_upper_and_lower_sets}, \cref{lem:prelim_inv_f_is_lower_set} states that the set $\inv(f)$ is a lower set in the partially ordered set $(\N, \mid)$, where $\mid$ is the divisibility relation. In particular, we can describe the invertibility properties of $f$ using $4$ different sets: the lower set $\inv(f)$, a generating set $\zeta(f)$ of $\inv(f)$, the upper set $\N\setminus\inv(f)$ or a generating set $\xi(f)$ of $\N\setminus\inv(f)$. Determining these sets given some shift-invariant $f$ is hard in general.

\begin{example}\label{ex:prelim_inv_set_of_chi}
    The function $\chi:\F_2^n\to\F_2^n$ is invertible if and only if $n$ is odd \cite{daemen1995cipher}. Hence, $\inv(\chi)$ is the set of all odd natural numbers. A generating set $\zeta(\chi)$ of $\inv(\chi)$ necessarily contains infinitely many numbers and there is no minimal generating set. The upper set $\N\setminus\inv(\chi)$ is the set of all even natural numbers and $\xi(\chi)=\{2\}$. In particular, $\xi(\chi)$ is finite. 
\end{example}

\begin{remark}
    For fixed $n$, choosing a basis of $\F_{2^n}$ over $\F_2$ gives an isomorphism $L:\F_{2^n}\to\F_2^n$, 
    and hence a map $f^L:\F_{2^n}\to\F_{2^n}$ corresponding to $f$. Every such map is uniquely represented by a polynomial in $\F_{2^n}[X]$ of degree at most $2^n-1$. By choosing a normal basis, the shift $S$ corresponds to the Frobenius automorphism $\sigma:\F_{2^n}\to\F_{2^n}$ given by $x\mapsto x^2$. Consequently, a shift-invariant map $f:\F_2^n\to\F_2^n$ corresponds to a map $f^L:\F_{2^n}\to\F_{2^n}$ satisfying $f^L(x^2)=(f^L(x))^2$ for all $x\in\F_{2^n}$. This condition is equivalent to $f^L$ being represented by a polynomial in $\F_2[X]$.

    The difficulty with this description is that the isomorphism $L$, and consequently the polynomial representing $f^L$, depends on $n$. As $n$ varies, there is therefore no canonical choice of polynomial representation that is useful for our purposes. For further details, as well as a study of the univariate representations of $\chi$, we refer to \cite{schoone2024algebraic}.
\end{remark}

\subsection{Conjugate equivalence}\label{subsec:prelims_notions_of_equivalences}

Let $R:\F_2^\Z\to\F_2^\Z$ be the reflection at the origin, i.e. 
\begin{equation*}
    R(\ldots, x_{-2},x_{-1},\origin{x_0},x_1,x_2,\ldots) = (\ldots, x_2, x_1, \origin{x_0}, x_{-1}, x_{-2}, \ldots)
\end{equation*}
and let $T:\F_2^\Z\to\F_2^\Z$ be given by $T(x)=x+\one$. Observe that $R^2=T^2=\id$ and $TR=RT$. In other words, $\{\id, R, T, RT\}$ is a group isomorphic to the Klein four-group.

Let $c=a_{-k}\ldots a_m$ be a landscape. We consider the set
\begin{equation*}
    \{(\id+c)^g: g\in\{\id, R, T, RT\}\} = \{\id+c, (\id+c)^R, (\id+c)^T, (\id+c)^{RT}\}
\end{equation*}
of conjugate maps of $\id+c$. Intuitively, conjugating by $R$ is the same as reading an element $x\in\F_2^\Z$ from right to left instead of left to right and conjugating by $T$ is the same as switching all $0$s and $1$s of $x$. We see that
\begin{align*}
    (\id+c)^R &= R\circ(\id+c)\circ R = \id+R\circ c\circ R = \id+c^R, \\
    (\id+c)^T &= T\circ(\id+c)\circ T = T+c\circ T + \one = \id+c\circ T, \\
    (\id+c)^{RT} &= \left((\id+c)^R\right)^T = (\id+c^R)^T = \id+c^R\circ T,
\end{align*}
and observe that these maps are given by the landscapes $c^R, c\circ T$ and $c^R\circ T$. These landscapes can be easily described:
\begin{itemize}[labelwidth=2cm,leftmargin=1.5cm]
    \item[$c^R$:] The landscape is $a_m a_{m-1}\ldots a_{-k}$, i.e. the landscapes of $c$ reflected at the highlighted position.
    \item[$c\circ T$:] The landscape is $b_{-k}\ldots b_m$ with $b_\ell=a_\ell+1$ if $a_\ell\neq \mt$ and $b_\ell=a_\ell=\mt$ otherwise, i.e. all $0$s are replaced by $1$s and vice versa.
    \item[$c^R\circ T$:] The landscape is $b_m b_{m-1}\ldots b_{-k}$ with $b_\ell$ as above.   
\end{itemize}

\begin{definition}
    The landscapes $c, c^R, c\circ T$ and $c^R\circ T$ are called conjugate equivalent. 
\end{definition}

For example, if $c=*01$ then the conjugate equivalent landscapes are $*10, 10*,01*$. 

\begin{remark}
    A related notion which is called elementary equivalence in \cite{haugland2025shift} and essential equivalence in \cite{omland2022permutation} also includes the four further functions 
    \begin{align*}
        &T\circ(\id+c) = \id+c+\one = \id + T\circ c, \\
        &T\circ(\id+c)^R = \id+c^R+\one = \id + T\circ c^R, \\
        &T\circ(\id+c)^T = \id+c\circ T+\one = \id+T\circ c\circ T, \\
        &T\circ(\id+c)^{TR} = \id+c^R\circ T + \one = \id+T\circ c^R\circ T.
    \end{align*}
    However, the functions $T\circ c, T\circ c^R, c^T, c^{TR}$ are not landscapes. Furthermore, as these functions are not conjugate to $\id+c$, they have different behavior under composition and we therefore do not include them here.
\end{remark}

\section{Conserved and quasi-conserved landscapes}\label{sec:quasi_conserved_landscapes}

For some landscapes $c$ it is easy to see that the associated function $f=\id+c$ is a permutation on $\F_2^\Z$. We consider the Patt function given by $c=0*10$ 
from \cite{amoroso1972decision}. Let $x\in\F_2^\Z$ such that $x$ matches $c$ at position $i$, for some $i\in\Z$. Then
\begin{equation*}
    x_{i-1} x_i x_{i+1} x_{i+2} = 0x_i10
\end{equation*}
and $x_i$ gets flipped under $f$, i.e. $y_i = x_i+1$, where $y=f(x)$. Our next goal is to determine $f(y)_i$. For this purpose we need to compute the values $y_{i-1},y_{i+1},y_{i+2}$. To determine $y_j$ we need to check whether $x$ matches $c$ at position $j$.
\begin{itemize}
    \item The vector $x$ matches $c$ at position $i-1$ if and only if $x_{i-2} x_{i-1} x_i x_{i+1} = 0x_{i-1}10$, in particular $x_{i+1}=0$. We know that $x_{i+1}=1$, so this is not the case. This can also be seen from the diagram obtained by aligning the pattern $c=0*10$ so that the $*$ is at the position $i-1$:
    \begin{equation*}
        \term[3pt]{
              & 0 & x_i & 1 & 0 \\
            0 & * &   1 & 0
        }
    \end{equation*}
    We have a conflict in the fourth column, because in the first row is a $1$ indicating that $x_{i+1}=1$, while in the second row there is a $0$. This means that $x$ does not match $c$ at position $i-1$ and therefore $y_{i-1}=x_{i-1}=0$.
    \item Similarly, for position $i+1$ we need that $x_i x_{i+1} x_{i+2} x_{i+3} = 0x_{i+1}10$. This is not satisfied, as $x_{i+2}=0\neq 1$. Visually, we see a conflict in the fourth column of 
    \begin{equation*}
        \term[3pt]{
            0 & x_i & 1 & 0 &   \\
              & 0   & * & 1 & 0
        },
    \end{equation*}
    since there is a $0$ in the first row and a $1$ in the second row. Hence, $x$ also does not match $c$ at position $i+1$ and $y_{i+1}=x_{i+1}=1$.
    \item Finally, for position $i+2$ we see that the third column of 
    \begin{equation*}
        \term[3pt]{
            0 & x_i & 1 & 0 &   \\
              &     & 0 & * & 1 & 0
        }
    \end{equation*}
    has a conflict, i.e. $x_{i+1}=1\neq 0$, so that $y_{i+2}=x_{i+2}=0$.
\end{itemize}
In all cases we see that $y_j=x_j$ for $j\in\{i-1,i+1,i+2\}$. We obtain
\begin{equation*}
    y_{i-1}y_iy_{i+1}y_{i+2} = 0 y_i 1 0.
\end{equation*}
In conclusion, if $x$ matches $c$ at position $i$, then $y$ also matches $c$ at position $i$.

Furthermore, we can also check that if $y$ matches $c$ at position $i$, then $x$ matches $c$ at position $i$. Indeed, assume otherwise, then $y_{i-1}y_iy_{i+1}y_{i+2}=0y_i10$ and $x_{i-1}x_ix_{i+1}x_{i+2}\neq 0x_i10$. In particular, at least one of $x_{i-1},x_{i+1}, x_{i+2}$ needs to get flipped under $f$, which implies that $x$ needs to match $c$ at one position $j\in\{i-1,i+1,i+2\}$. But then by the same argument as before, $y$ has to match $c$ at position $j$ as well. But we showed that it is impossible for a vector $y$ to match $c$ at position $i$ and another position in $\{i-1,i+1,i+2\}$, hence we obtain a contradiction.

To summarize, $x$ matches $c$ at position $i$ if and only if $y$ matches $c$ at position $i$. In particular, the set of positions $i$, such that $x$ matches $c$ at position $i$ is conserved under application of $f$. The landscape $c$ is called conserved in \cite{toffoli1990invertible}. Furthermore, we see that $f(y)=x$, implying that $f$ is an involution on $\F_2^\Z$ and in particular a permutation.

Note that $c=0*10$ is a two-sided landscape. A natural question is whether a one-sided landscape can also lead to an involution. The answer is negative. Indeed, let $c=*a_1\ldots a_m$ be a right-sided landscape and assume that $x$ matches $c$ at position $i$. We now need that any such $x$ does not match $c$ at all positions $i+\ell$ for $\ell\in\supp(c)$. We can again check this visually by constructing the diagram
\begin{equation}\label{eq:qc_x_matches_c_at_i_to_i_plus_m}
    \term[3pt]{
        x_i & x_{i+1} & x_{i+2} & x_{i+3} & \ldots & x_{i+m-1} & x_{i+m}     \\
        *   & a_1     & a_2     & a_3     & \ldots & a_{m-1}   & a_m          \\
            & *       & a_1     & a_2     & \ldots & a_{m-2}   & a_{m-1} & a_m \\
            &         & *       & a_1     & \ldots & a_{m-3}   & a_{m-2} & a_{m-1} & a_m \\
            &         &         &         &        & \ddots    & \ddots  & \ddots  \\
            &         &         &         &        & *         & a_1     & a_2     & a_3 & a_4 & \ldots \\
            &         &         &         &        &           & *       & a_1     & a_2 & a_3 & a_4    & \ldots 
    }
\end{equation}
We know that $x$ matches $c$ at position $i$, which means that in the first two rows of \cref{eq:qc_x_matches_c_at_i_to_i_plus_m} there is no conflict and $x_{i+\ell}=a_\ell$ for all $\ell\in\supp(c)$. Next, to determine whether $x$ matches $c$ at position $i+1$, we need to compare the first and third row, or equivalently, the second and third row. It is possible that there is a conflict in those rows, for example if $a_1=1$ but $a_2=x_{i+2}=0$. This then implies that $x$ does not match $c$ at position $i+1$, no matter the other entries of $x$. By a similar argument it might be impossible that $x$ matches $c$ at all positions $i+\ell$ for $\ell\in\supp(c)\setminus\{m\}$. 

For position $i+m$ we have a different situation. The diagram \cref{eq:qc_x_matches_c_at_i_to_i_plus_m} tells us that whether $x$ matches $c$ at position $i+m$ is completely independent of the entries $x_{i+1},\ldots,x_{i+m}$ in the first row. In particular, it is possible to choose $x_{i+m+\ell}=a_\ell$ for every $\ell\in\supp(c)$, which then means that $x$ matches $c$ at position $i+m$. 
We take this as motivation for the following definition.

\begin{definition}[Conserved and quasi-conserved landscapes]\label{def:qc_conserved_and_quasi_conserved_landscapes}
    Let $c=a_{-k}\ldots a_m$ be a landscape. 
    \begin{itemize}
        \item We say that $c$ is a conserved landscape, if whenever $x\in\F_2^\Z$ matches $c$ at position $i\in\Z$, then $x$ does not match $c$ at all positions $i+\ell$ for $\ell\in\supp(c)$.
        \item We say that $c$ is a quasi-conserved landscape, if $c$ is not conserved and there exists one index $q\in\supp(c)$ such that whenever $x\in\F_2^\Z$ matches $c$ at position $i\in\Z$, then $x$ does not match $c$ at all positions $i+\ell$ for $\ell\in\supp(c)\setminus\{q\}$. Furthermore, we also say that $(c,q)$ is quasi-conserved or $c$ is quasi-conserved with special index $q$.
    \end{itemize}
\end{definition}

\begin{assumption}\label{ass:qc_assume_not_affine}
    From now on we assume that $\abs{\supp(c)}>1$, as otherwise $\supp(c)\setminus\{q\}$ is empty and many of the subsequent claims are either trivial or vacuous. Furthermore, if $\abs{\supp(c)}=1$, then $c$ is affine and affine functions are already well-understood. 
\end{assumption}

\begin{example}\label{ex:qc_conserved_and_quasi_conserved_landscapes}
    Let $x\in\F_2^\Z$. 
    \begin{enumerate}
        \item The previous discussion shows that $0*10$ is a conserved landscape.
        \item The landscape $c=*01$ is a quasi-conserved landscape with $q=2$. Indeed, suppose that $x_i x_{i+1} x_{i+2} = x_i01$. Then $x$ does not match $c$ at position $i+1$, as $x_{i+2}=1$ while it would need to be $0$. Visually we read this off from the conflict in the third column of the diagram
        \begin{equation*}
            \term[3pt]{
                x_i & 0 & 1 \\
                    & * & 0 & 1.
            }
        \end{equation*}
        \item More generally, $*0\ldots01$ is quasi-conserved for any $m\geq 2$. Indeed, suppose $x_i x_{i+1}\ldots x_{i+m-1} x_{i+m} = x_i 0\ldots 01$. From
        \begin{equation*}
            \term[3pt]{
                x_i & 0 & \ldots & 0 & \ldots & 0 & 1 \\
                    &   &        & * & \ldots & 0 & 0 & 0 & \ldots & 0 & 1 
            }
        \end{equation*}
        it follows that $x$ does not match $*0\ldots01$ at any position $i+\ell$ for $\ell\in\{1,\ldots,m-1\}$.
        \item Let $d=0*0\mt\,\mt01$. Then $(d,5)$ is quasi-conserved. Indeed, suppose that $x$ matches $d$ at position $i$. Then $x$ does not match $d$ at positions $i-1, i+1, i+4$ as we see from the seventh and eighth column of
        \begin{equation*}
            \term[3pt]{
                  & 0 & x_i &   0 & x_{i+2} & x_{i+3} & 0   & 1 \\
                0 & * &   0 & \mt &     \mt &       0 & 1   \\
                  &   &   0 & *   &     0   &     \mt & \mt & 0 & 1 \\
                  &   &     &     &         &       0 & *   & 0 & \mt & \mt & 0 & 1.
            }
        \end{equation*}
        However, $x$ may match $d$ at positions $i+2,i+3,i+5$. Note that this does not go against the definition of a quasi-conserved landscape, as $2,3\notin\supp(d)=\{-1,1,4,5\}$. 
        \item The landscape $0*\mt10$ is a quasi-conserved landscape with $q=2$. Indeed, suppose that $x$ matches $0*\mt10$ at position $i$. Then $x$ does not match $0*\mt10$ at positions $i-1$ and $i+3$ as we see from the fifth column of
        \begin{equation*}
            \term[3pt]{
                  & 0 & x_i & \mt & 1 & 0 \\
                0 & * & \mt &   1 & 0   \\
                  &   &     &     & 0 & * & \mt & 1 & 0.
            }
        \end{equation*}
        \item The landscape $*01\mt1$ is neither conserved nor quasi-conserved. For example, if $x_i x_{i+1} x_{i+2} x_{i+3} x_{i+4} x_{i+5} x_{i+6}= 0010111$ then $x$ matches $*01\mt1$ at positions $i$ and $i+2$. 
    \end{enumerate}
\end{example}

Many of the shift-invariant functions considered in the literature are given by quasi-conserved landscapes. 

\begin{example}
    \begin{table}[!ht]
        \centering
        \begin{tabular}{|c|c|c|}
            \hline
            Landscape & conserved or quasi-conserved & $q$ if quasi-conserved \\
            \hline
            $*01$ & quasi-conserved & $2$ \\
            $0*01$ & quasi-conserved & $2$ \\
            $*001$ & quasi-conserved & $3$ \\
            $*100$ & quasi-conserved & $3$ \\
            $0*\mt10$ & quasi-conserved & $2$ \\
            $*01\mt1$ & neither & - \\
            $*0\mt01$ & quasi-conserved & $4$ \\
            $*01\mt0$ & quasi-conserved & $4$ \\
            $0\mt_r*1\mt_r0$ & quasi-conserved & $1$ \\
            $*0100$ & neither & - \\
            $0*010$ & quasi-conserved & $2$ \\
            $0*001$ & quasi-conserved & $3$ \\
            $0*100$ & quasi-conserved & $3$ \\
            $0*101$ & quasi-conserved & $3$ \\
            $*0111$ & quasi-conserved & $4$ \\
            $*0001$ & quasi-conserved & $4$ \\
            $*0011$ & quasi-conserved & $4$ \\
            $0*01\mt0$ & quasi-conserved & $4$ \\
            $0*011$ & conserved & - \\
            $0\mt*\mt11\mt0$ & conserved & - \\
            $0*1\mt_r10$ & conserved & - \\
            $01*01$ & conserved & - \\
            $0\mt_{r+1}*1\mt_r10$ & conserved & - \\
            $0\mt*100$ & conserved & - \\
            \hline
        \end{tabular}
        \caption{\label{tb:landscapes_from_daemen} (Table A.1 in \cite{daemen1995cipher}) Here, $\mt_r$ denotes the symbol $\mt$ repeated $r$ times.}
    \end{table}
    The function $\chi$ is given by $\chi=\id+*01$ \cite{daemen1995cipher}. By \cref{ex:qc_conserved_and_quasi_conserved_landscapes} $(b)$ the landscape $*01$ is quasi-conserved. A generalization of $\chi$ given by $\id+*0\ldots01$ is introduced in \cite{lyu2026generalized}. Note that $*0\ldots 01$ is also quasi-conserved by \cref{ex:qc_conserved_and_quasi_conserved_landscapes} $(c)$. We call the corresponding function $\chi^{(m)}$.

    The primitive \textsc{Monolith} \cite{grassi2024monolith} uses the function $\id+*100$. The landscape $*100$ is also quasi-conserved.

    The landscapes $*01, *001$ and $*100$ are given in Table A.1 of Daemen's thesis \cite{daemen1995cipher} which in total contains $24$ landscapes. Of these, $16$ are quasi-conserved, $6$ are conserved and $2$ are neither, see \cref{tb:landscapes_from_daemen}. In \cite[Table 1]{liu2026finding} another list of landscapes is given for which the tools in their paper can compute the inverse function. In fact, all of these landscapes are already contained in \cite[Table A.1]{daemen1995cipher} and are either conserved or quasi-conserved.
\end{example}

The following example shows that quasi-conservedness is not preserved under shifting the landscape. This is important, since otherwise we could for example shift any quasi-conserved two-sided landscape to obtain a quasi-conserved one-sided landscape, which are easier to work with.

\begin{example}\label{ex:two_sided_quasi_conserved_not_preserved_under_shift}
    Let $c=0*\mt10$ with $\supp(c)=\{-1,2,3\}$. Then $(c,2)$ is quasi-conserved (see \cref{ex:qc_conserved_and_quasi_conserved_landscapes}). 
    Now consider $\tilde{c} = S^2c = *0\mt\,\mt 10$ with $\supp(\tilde c) = \{1,4,5\}$. Then $\tilde c$ is not quasi-conserved. To see this, take $x\in\F_2^\Z$ which matches $\tilde c$ in position $i\in\Z$. From the diagram
    \begin{equation*}
        \term[2pt]{
            x_i & 0 & \mt & \mt & 1 & 0 \\
                &   &     &     & * & 0 & \mt & \mt & 1   & 0 \\
                &   &     &     &   & * & 0   & \mt & \mt & 1 & 0
        }
    \end{equation*}
    we can read off that $x$ could match $c$ in positions $i+4, i+5$. Hence, $\tilde c$ is not quasi-conserved.

    In fact, the only shift of $c$ which is quasi-conserved is $S^{-2}c = 0\mt\,\mt\hat{1}0$, which is also two-sided.

\end{example}

\begin{example}
    Let $c=*01$. Then by \cref{ex:qc_conserved_and_quasi_conserved_landscapes} $(b)$ we know that $c$ is quasi-conserved. If $x=(x_0,\ldots,x_8)\in\F_2^9$ is specified as follows, then
    \begin{align*}
             x  &= (0, 1, 1, 0, 1, 0, 1, 0, 0), \\
        \chi(x) &= (0, 1, 0, 0, 0, 0, 1, 0, 1).
    \end{align*}
    We see that $x$ matches $c$ at positions $i=2,4,8$ while $y=\chi(x)$ matches $c$ at positions $i=4,6,8$. For conserved landscapes we saw that $x$ matches $c$ at position $i$ if and only if $y$ matches $c$ at position $i$. Neither implication is true for quasi-conserved landscapes.
\end{example}

From all the previous examples it is clear that, given a landscape $c=a_{-k}\ldots a_m$ and some vector $x\in\F_2^\Z$ which matches $c$ at some position $i\in\Z$, we can decide whether $x$ matches $c$ at position $i+\ell$ by looking at a diagram such as 
\begin{equation}\label{eq:qc_checking_x_matches_c_position_i_plus_ell}
    \term[2pt]{
        a_{-k} & \ldots & a_{-1}      & \hat{a_0} & a_1         & \ldots & a_{\ell-1} & a_\ell    & a_{\ell+1} & \ldots & a_m \\
               & \ldots & a_{-\ell-1} & a_{-\ell} & a_{-\ell+1} & \ldots & a_{-1}     & \hat{a_0} & a_1        & \ldots & a_{m-\ell} & \ldots & a_m
    }
\end{equation}
and checking whether there is one column with a $0$ and a $1$ (if there is such a column, then $x$ does not match $c$ at position $i+\ell$). Note that this diagram is essentially equivalent to the diagram that is obtained when we compute $S^\ell c\odot c$, which is given by 
\begin{equation}\label{eq:qc_landscape_depiction_S_ell_c_odot_c}
    S^\ell c\odot c = 
    \term[2pt]{
               &        &      a_{-k} & \ldots & a_{-\ell-1} & \hat{a_{-\ell}} & a_{-\ell+1} & \ldots & a_{-1}     & a_0    & a_1        & \ldots & a_{m-\ell} & \ldots & a_m \\
        a_{-k} & \ldots & a_{-k+\ell} & \ldots & a_{-1}      & \hat{a_0}       & a_1         & \ldots & a_{\ell-1} & a_\ell & a_{\ell+1} & \ldots & a_m
    }
\end{equation}
Furthermore, $S^\ell c\odot c=\zero$ if and only if there exists a column with a $0$ and a $1$ in \cref{eq:qc_landscape_depiction_S_ell_c_odot_c}, which is equivalent to there being such a column in \cref{eq:qc_checking_x_matches_c_position_i_plus_ell}. We obtain the subsequent lemma:

\begin{lemma}\label{lem:qc_x_matches_at_pos_i_plus_ell_iff_S_ell_c_odot_c_equals_zero}
    Let $c=a_{-k}\ldots a_m$ be a landscape and $\ell\in\supp(c)$. Then the following are equivalent:
    \begin{enumerate}
        \item For all $x\in\F_2^\Z$ and all $i\in\Z$ it holds that if $x$ matches $c$ at position $i$, then $x$ does not match $c$ at position $i+\ell$.
        \item It holds that $S^\ell c\odot c=\zero$. 
        \item There exists a conflict in \cref{eq:qc_checking_x_matches_c_position_i_plus_ell} and \cref{eq:qc_landscape_depiction_S_ell_c_odot_c}.
        \item There exists $\rho\in\supp(c)$ such that $\rho+\ell\in\supp(c)$ and $(a_\rho,a_{\rho+\ell})\in\{(1,0),(0,1)\}$.
        \item There exists $\rho\in\supp(c)$ such that $\rho+\ell\in\supp(c)$ and $S^\ell A_\rho \odot A_{\rho+\ell}=\zero$.
    \end{enumerate}
\end{lemma}

\begin{proof}
    By the previous discussion, $(a),(b)$ and $(c)$ are equivalent. Statement $(d)$ is a rephrasing of $(c)$. Note that
    \begin{equation*}
        S^\ell A_{\rho} \odot A_{\rho+\ell}
            = \term[2pt]{
                * & \mt & \ldots & \mt & a_\rho \\
                * & \mt & \ldots & \mt & a_{\rho+\ell} 
            }
    \end{equation*}
    so that $(e)$ is another rephrasing of $(d)$.
\end{proof}

We then obtain the following characterization of conserved and quasi-conserved landscapes.

\begin{lemma}\label{lem:qc_equivalent_characterization_conserved_quasi_conserved_with_S_ell_c_odot_c}
    Let $c=a_{-k}\ldots a_m$ be a landscape. Then:
    \begin{enumerate}
        \item $c$ is conserved if and only if $S^\ell c\odot c=\zero$ for all $\ell\in\supp(c)$.
        \item $(c,q)$ is quasi-conserved if and only if $S^\ell c\odot c=\zero$ for all $\ell\in\supp(c)\setminus\{q\}$.
    \end{enumerate}
\end{lemma}

We note that \cref{lem:qc_equivalent_characterization_conserved_quasi_conserved_with_S_ell_c_odot_c} combined with \cref{lem:qc_x_matches_at_pos_i_plus_ell_iff_S_ell_c_odot_c_equals_zero} gives an easy way to automatically check whether a landscape is quasi-conserved with a computer. 

More generally, for some landscape $c$ we define the set
\begin{equation*}
    \mL_c = \{\ell\in\Z: S^\ell c\odot c=\zero\}.
\end{equation*}
It is clear that $\mL_c$ is finite, as the supports of $S^\ell c$ and $c$ are disjoint for $\ell$ large enough.
\cref{lem:qc_equivalent_characterization_conserved_quasi_conserved_with_S_ell_c_odot_c} then states that $c$ is conserved if and only if $\supp(c)\subseteq \mL_c$ and quasi-conserved if and only if $\supp(c)\setminus\{q\}\subseteq \mL_c$.

\begin{example}
    Let $c=*00\mt101$. By direct calculation we see that $\mL_c=\{-5,\ldots,5\}\setminus\{0\}$. Observe that $\supp(c)=\{1,2,4,5,6\}$. With the choice $q=6$ we see $\supp(c)\setminus\{6\}\subseteq \mL_c$, hence $(c,6)$ is quasi-conserved.
\end{example}

The subsequent lemma shows that quasi-conservedness is preserved under conjugate equivalence.

\begin{lemma}\label{lem:qc_conjugate_equivalence_preserves_quasi_conserved}
    Let $c$ be a landscape. Then the following are equivalent:
    \begin{enumerate}
        \item $c=a_{-k}\ldots a_m$ is quasi-conserved with special index $q$.
        \item $c\circ T=b_{-k}\ldots b_m$ is quasi-conserved with special index $q$. 
        \item $c^R=a_m\ldots a_{-k}$ is quasi-conserved with special index $-q$.
        \item $c^R\circ T=b_m\ldots b_{-k}$ is quasi-conserved with special index $-q$. 
    \end{enumerate}
\end{lemma}

\begin{proof}
    By choosing $\tilde c\in\{c,c\circ T, c^R, c^R\circ T\}$ and using $T^2=R^2=\id$ and $TR=RT$ it suffices to show that $(a)$ implies $(b)$ and $(c)$.

    $(a)\implies (b)$: By \cref{lem:qc_equivalent_characterization_conserved_quasi_conserved_with_S_ell_c_odot_c} $(c\circ T, q)$ is quasi-conserved if $(S^\ell (c\circ T))\odot (c\circ T)=\zero$ for all $\ell\in\supp(c\circ T)\setminus\{q\}=\supp(c)\setminus\{q\}$. Indeed, we have
    \begin{equation*}
        (S^\ell (c\circ T))\odot (c\circ T) = (S^\ell c\odot c)\circ T = \zero
    \end{equation*}
    as $S^\ell c\odot c=\zero$ because $c$ is quasi-conserved.

    $(a)\implies (c)$: Note that $S^R = S^{-1}$ and $f^R\odot g^R = (f\odot g)^R$ for functions $f,g:\F_2^\Z\to\F_2^\Z$. We show that $S^\ell c^R\odot c^R=\zero$ for all $\ell\in\supp(c^R)\setminus\{-q\}=-\supp(c)\setminus\{q\}$. Indeed,
    \begin{equation*}
        S^\ell c^R\odot c^R = (S^{-\ell}c)^R\odot c^R = (S^{-\ell}c\odot c)^R = \zero
    \end{equation*}
    because $c$ is quasi-conserved and $-\ell\in\supp(c)\setminus\{q\}$. 
\end{proof}

\subsection{The polynomial composition property}

In \cite{kriepke2024algebraic} we studied the iterates of the function $\chi=\id+*01$ over $\F_2^n$ with $n$ odd. This resulted in the definition of a family of functions $\gamma_{2j}:\F_2^n\to\F_2^n$ for $j\geq 0$, where $\gamma_0=\id$ is the identity function and 
\begin{equation*}
    \gamma_{2j}(x)_i = x_{i+2j}(1+x_{i+2j-1})(1+x_{i+2j-3})\cdots (1+x_{i+i})
\end{equation*}
for $j\geq 1$. In landscape notation, 
\begin{equation}\label{eq:qc_landscape_notation_gamma_2j}
    \gamma_{2j} = *0\mt 0\mt \ldots 0\mt 01.
\end{equation}
We showed that
\begin{equation}\label{eq:qc_polynomial_composition_property_gamma_2j}
    \gamma_{2j}\circ \left(\sum_{i=0}^s \alpha_i \gamma_{2i}\right) = \sum_{i=0}^s \alpha_i \gamma_{2i+2j}
\end{equation}
if $\alpha_i\in\F_2$ with $\alpha_0=1$. Observe that \cref{eq:qc_polynomial_composition_property_gamma_2j} looks similar to polynomial multiplication. Indeed, by replacing $\gamma_{2j}$ with $X^j$ and composition with multiplication, \cref{eq:qc_polynomial_composition_property_gamma_2j} turns into
\begin{equation*}
    X^j\cdot \left(\sum_{i=0}^s \alpha_i X^i\right) = \sum_{i=0}^s \alpha_i X^{i+j}.
\end{equation*}
This leads to a connection of the iterates of $\chi$ and a polynomial quotient ring, from which several properties about $\chi$ are proven. More precisely, set
\begin{equation*}
    \Gamma = \spn\{\gamma_{2j}: j\geq 0\}
\end{equation*}
and 
\begin{equation*}
    \mG = \gamma_0+\spn\{\gamma_{2j}: j\geq 1\}.
\end{equation*}
Observe that $\chi=\gamma_0+\gamma_2\in\mG$. By using \cref{eq:qc_polynomial_composition_property_gamma_2j} we proved that $(\mG,\circ)$ is a monoid. Define the quotient ring
\begin{equation*}
    \mR = \F_2[X]/(X^{(n+1)/2})
\end{equation*}
and the subset
\begin{equation*}
    \mM = \left\{1+\sum_{i=1}^s \alpha_i X^i \right\}\subseteq \mR.
\end{equation*}
Then $\varphi: \Gamma\to \mR$ given by 
\begin{equation*}
    \varphi\left(\sum_{i=0}^s \alpha_i \gamma_{2i}\right) = \sum_{i=0}^s \alpha_i X^i 
\end{equation*}
is a linear isomorphism and its restriction $\varphi:\mG\to\mM$ is a monoid isomorphism. Since $\mM$ is a group, it follows that $\mG$ is also a group. This shows that every map $f\in\mG$ is a permutation on $\F_2^n\to\F_2^n$. While \cite{kriepke2024algebraic} focused on the $n$ odd case, the case $n$ even was then covered in \cite{kriepke2026siblings}. Here, the quotient ring has to be chosen as 
\begin{equation*}
    \mR = \F_2[X]/(X^n+X^{n/2}),
\end{equation*}
but many of the same properties still hold. This time $\mM$ is not a group, so that $\mG$ is also only a monoid. However, it is well-known that $F(X)\in\mM$ is a unit if and only if $\gcd(F(X), X^n+X^{n/2})=1$. This enables a complete classification of all permutations in $\mG$. The paper \cite{lyu2026generalized} generalized the tools of \cite{kriepke2024algebraic} to the case of the functions $\chi^{(m)}$ based on the landscapes $*0\ldots 01$ with some functions $\theta_{m,j}$ replacing the functions $\gamma_{2j}$ and for $n$ not a multiple of $m$.

Our goal from here is to generalize the identity \cref{eq:qc_polynomial_composition_property_gamma_2j}. We begin by defining functions that include the functions $\gamma_{2j}$ from \cite{kriepke2024algebraic,kriepke2026siblings} and $\theta_{m,j}$ from \cite{lyu2026generalized} as special cases.

\begin{definition}\label{def:qc_gamma_functions}
    Let $c=a_{-k}\ldots a_m$ be a landscape and take $q\in\supp(c)$ arbitrary. Define the landscape $\cq$ by replacing $a_q$ in $c$ with $\mt$, i.e. 
    \begin{equation*}
        \cq = a_{-k}\ldots a_{q-1}\,\mt\, a_{q+1}\ldots a_m. 
    \end{equation*}
    We define functions $\gamma_j^{(c,q)}:\F_2^\Z\to\F_2^\Z$ for $j\in\N_0$ by $\gamma_0^{(c,q)}=\id$, $\gamma_1^{(c,q)}=c$ and 
    \begin{equation*}
        \gamma_{j+1}^{(c,q)}=\cq \odot S^q \gamma_j^{(c,q)}
    \end{equation*}
    for $j\geq 1$. If $c$ and $q$ are clear from context then we simply write $\gamma_j$. 
\end{definition}

\begin{example}\label{ex:qc_gammas_for_*01}
    Let $c=*01$ which is quasi-conserved with $q=2$ by \cref{ex:qc_conserved_and_quasi_conserved_landscapes}. Then $c^{(2)}=*0\mt$ and therefore
    \begin{align*}
        \gamma_0^{(*01,2)} &=\id = \hat 1, \\
        \gamma_1^{(*01,2)} &=c = *01,   \\
        \gamma_2^{(*01,2)} &=c^{(2)}\odot S^2(*01) = *0\mt\odot *\mt\,\mt 01 =*0\mt 01, \\
        \gamma_3^{(*01,2)} &=c^{(2)}\odot S^2(*0\mt01) = *0\mt\odot *\mt\,\mt0\mt01 = *0\mt0\mt01.
    \end{align*}
    These are exactly the functions $\gamma_0,\gamma_2,\gamma_4,\gamma_6$ of \cite{kriepke2024algebraic,kriepke2026siblings}, see \cref{eq:qc_landscape_notation_gamma_2j}. 
\end{example}

We want to convey some intuition about why the functions $\gamma_j^{(c,q)}$ appear when considering the iterates of $\id+c$. Take the landscape $*001$, then $\gamma_2=*00\mt001$ and $\gamma_3=*00\mt00\mt001$. Observe that $\id+*001=\chi^{(3)}$ and $\gamma_2=\theta_{3,2}, \gamma_3=\theta_{3,3}$ from \cite{lyu2026generalized}. Take $x\in\F_2^\Z$ and set $x^{(j)} = (\chi^{(3)})^j(x)$. Our first goal is to study $x_0^{(2)}$.

By the definition of the function $\chi^{(3)}$ the value of $x_0^{(2)}$ depends on whether $x_1^{(1)}x_2^{(1)}x_3^{(1)}=001$. Assume first that $x_1x_2\neq 00$. Then at least one of them needs to be flipped by $\chi^{(3)}$ in order for the equality $x_1^{(1)}x_2^{(1)}=00$ to hold. This requires at the very least that $x_3=0$, meaning it also needs to get flipped. As $*001$ is quasi-conserved it is impossible that $x_3$ and at least one of $x_1$ or $x_2$ can get flipped. This leads to a contradiction, so we need $x_1x_2=00$. Next, the value of $x_3^{(1)}$ depends on whether $x_4x_5x_6=001$. Combining this with the previous argument, the value of $x_0^{(2)}$ is related to whether $x_1x_2x_3x_4x_5x_6=00x_3001$, i.e. whether $x$ matches $\gamma_2=*00\mt001$ at position $0$. In fact, \cref{ex:qc_wpcp_for_right_sided_landscape} shows that
\begin{equation*}
    (\chi^{(3)})^2 = \id+*00\mt001.
\end{equation*}
Now we want to determine $x_0^{(3)}$ which depends on whether $x_1^{(2)}x_2^{(2)}x_3^{(2)}=001$. By a similar argument as before, using that $*001$ is quasi-conserved, the equality $x_1^{(2)}x_2^{(2)}=00$ is influenced by the condition $x_1x_2=00$. Furthermore, by the result of the previous paragraph, the value of $x_2^{(3)}$ depends on the equality $x_4x_5x_6x_7x_8x_9=00x_6001$. Combining these we need to check the equality $x_1x_2x_3x_4x_5x_6x_7x_8x_9=00x_300x_6001$, i.e. whether $x$ matches $\gamma_3$ at position $0$. With a computation similar to \cref{ex:qc_wpcp_for_right_sided_landscape} or via the results of \cite{lyu2026generalized} one can show that
\begin{equation*}
    (\chi^{(3)})^3 = \id+*001+*00\mt001+*00\mt00\mt001.
\end{equation*}

\begin{lemma}\label{lem:qc_gamma_functions_of_conjugates}
    Let $c$ be a landscape. Then:
    \begin{align*}
        \gamma_j^{(c\circ T,q)} &= \gamma_j^{(c,q)}\circ T, \\
        \gamma_j^{(c^R,-q)} &= (\gamma_j^{(c,q)})^R, \\
        \gamma_j^{(c^R\circ T,-q)} &= (\gamma_j^{(c,q)})^R\circ T
    \end{align*}
    for $j\geq 1$.
\end{lemma}

\begin{proof}
    By the same argument as in the proof of \cref{lem:qc_conjugate_equivalence_preserves_quasi_conserved} it is enough to show the claim for $\gamma_j^{(c\circ T,q)}$ and $\gamma_j^{(c^R,-q)}$. We proceed by induction on $j$. The base case $j=1$ is clear. Then
    \begin{align*}
        \gamma_{j+1}^{(c\circ T,q)} &= (c\circ T)^{(q)}\odot (S^q (\gamma_j^{(c\circ T,q)}\circ T)) \\
            &= (\cq\circ T) \odot ((S^q \gamma_j^{(c,q)})\circ T )\\
            &= (\cq \odot S^q \gamma_j^{(c,q)})\circ T \\
            &= \gamma_{j+1}^{(c,q)}\circ T 
    \end{align*}
    and 
    \begin{align*}
        \gamma_{j+1}^{(c^R,-q)} &= (c^R)^{(-q)}\odot S^{-q} \gamma_j^{(c^R, -q)} \\
            &= (\cq )^R\odot S^{-q} (\gamma_j^{(c,q)})^R \\
            &= (\cq)^R \odot (S^q \gamma_j^{(c,q)})^R \\
            &= (\cq \odot S^q \gamma_j^{(c,q)})^R \\
            &= (\gamma_{j+1}^{(c,q)})^R.
    \end{align*}
\end{proof}

Note that so far we have allowed that $0\in\supp(c)$. However, if $c$ is quasi-conserved and $0\in\supp(c)$, then $q=0$ and the functions $\gamma_j^{(c,q)}$ are degenerate. Indeed, as $S^0c \odot c=c\neq \zero$ but $S^\ell c\odot c=\zero$ for all $\ell\in\supp(c)\setminus\{q\}$ by \cref{lem:qc_x_matches_at_pos_i_plus_ell_iff_S_ell_c_odot_c_equals_zero}, we need $q=0$. But then $\gamma_1=c$ and $\gamma_2=\cq\odot S^q c=c^{(0)}\odot c=c$ and similarly $\gamma_j=c$ for every $j\geq 2$. 

There is also a second reason to exclude the case $0\in\supp(c)$. Let $c$ be a landscape with $0\in\supp(c)$. Then $f=\id+c$ is never a permutation on $\F_2^n$ for $n$ large enough and also not on $\F_2^\Z$. Indeed, take any $x\in \F_2^n$ such that $x$ matches $c$ at position $0$. Such a vector $x$ exists as long as $n$ is large enough. Suppose without loss of generality that $a_0=0$, i.e. $x_0=0$. Consider the sequence $x^{(j)}=f^j(x)$ we obtain from repeatedly applying $f$ to $x$. If $f$ is a permutation on $\F_2^n$, then eventually $x^{(j)}=x$ for some $j\geq 1$. We know that $x_0=0$ and $x$ matches $c$ at position $0$, implying $x^{(1)}_0=1$. In particular, $x^{(1)}$ does not match $c$ at position $0$ because $x^{(1)}_0=1\neq 0=a_0$. Then $x^{(2)}_0=1$. It follows that $x^{(j)}_0=1$ for all $j\geq 1$, so that $x^{(j)}\neq x$ for all $j\geq 1$. Hence, $f$ is not a permutation on $\F_2^n$. 

Observe that the recurrence $\gamma_{j+1}^{(c,q)}=\cq\odot S^q \gamma_j^{(c,q)}$ also holds for $j=0$ if and only if $a_q=1$. Next, $q>0$ will be convenient sometimes. Observe that by replacing $c$ by one of the conjugate equivalent landscapes $c\circ T, c^R, c^R\circ T$ from \cref{subsec:prelims_notions_of_equivalences} we may assume that both $q>0$ and $a_q=1$ hold. 

\begin{definition}\label{def:qc_degenerate_and_normed_landscape}
    Let $c$ be a landscape. We say that
    \begin{enumerate}
        \item[(i)] $c$ is degenerate if $0\in\supp(c)$ and nondegenerate otherwise. 
        \item[(ii)] $(c,q)$ is normed if $q>0$ and $a_q=1$.
    \end{enumerate}
\end{definition}

We now generalize the identity \cref{eq:qc_polynomial_composition_property_gamma_2j}. For this we simply replace the functions $\gamma_{2j}$ by $\gamma_j$. We still only limit ourselves to the case $\alpha_0=1$, as this is the only case that is relevant to us. 

\begin{definition}\label{def:qc_polynomial_composition_property}
    Let $\{\gamma_j: j\geq 0\}$ be a family of functions $\F_2^\Z\to\F_2^\Z$. We say that this family satisfies the polynomial composition property, if for every $j\geq 0$ and $\alpha_0,\ldots,\alpha_s\in\F_2$ with $\alpha_0=1$ it holds that
    \begin{equation*}
        \gamma_j \circ \left(\gamma_0+\sum_{i=1}^s \alpha_i \gamma_i\right) = \sum_{i=0}^s \alpha_i \gamma_{i+j}.
    \end{equation*}
\end{definition}

Obviously, if $L$ is linear then $\{L^j: j\geq 0\}$ satisfies the polynomial composition property also for $\alpha_0=0$. But, as for example \cite[Remark 2]{kriepke2024algebraic} shows, if $\alpha_0=0$ then \cref{eq:qc_polynomial_composition_property_gamma_2j} is false in general. Whether there are nonlinear functions which satisfy the polynomial composition property even for $\alpha_0=0$ is not clear. 

\begin{problem}
    Are there any nonlinear functions which satisfy the polynomial composition property also for $\alpha_0=0$ ?
\end{problem}

\begin{definition}\label{def:qc_polynomial_composition_property_for_landscape_c}
    Let $c$ be a landscape and $q\in\supp(c)$. We say that $(c,q)$ satisfies the polynomial composition property if the family $\{\gamma_j^{(c,q)}: j\geq 0\}$ satisfies the polynomial composition property.
\end{definition}

\begin{lemma}\label{lem:qc_conjugate_equivalence_preserves_pcp}
    Let $c$ be a landscape. Then the following are equivalent:
    \begin{enumerate}
        \item $(c,q)$ satisfies the polynomial composition property.
        \item $(c\circ T,q)$ satisfies the polynomial composition property.
        \item $(c^R,-q)$ satisfies the polynomial composition property.
        \item $(c^R\circ T,-q)$ satisfies the polynomial composition property.
    \end{enumerate}
\end{lemma}

\begin{proof}
    It suffices to show that $(a)$ implies $(b)$ and $(c)$. 

    $(a)\implies (b)$: Using \cref{lem:qc_gamma_functions_of_conjugates},
    \begin{align*}
        \gamma_j^{(c\circ T,q)}\circ \left(\gamma_0^{(c\circ T,q)}+\sum_{i=1}^s \alpha_i \gamma_i^{(c\circ T,q)}\right)
            &= \gamma_j^{(c,q)}\circ T\circ \left(\id+\sum_{i=1}^s \alpha_i \gamma_i^{(c,q)}\circ T\right) \\
            &= \gamma_j^{(c,q)}\circ \left(\id+\sum_{i=1}^s \alpha_i \gamma_i^{(c,q)}\circ T+\one\right) \\
            &= \gamma_j^{(c,q)}\circ \left(T+\sum_{i=1}^s \alpha_i \gamma_i^{(c,q)}\circ T\right) \\
            &= \gamma_j^{(c,q)}\circ \left(\id+\sum_{i=1}^s \alpha_i \gamma_i^{(c,q)}\right)\circ T \\
            &= \left(\sum_{i=0}^s \alpha_i \gamma_{i+j}^{(c,q)}\right)\circ T \\
            &= \sum_{i=0}^s \alpha_i \gamma_{i+j}^{(c,q)}\circ T \\
            &= \sum_{i=0}^s \alpha_i \gamma_{i+j}^{(c\circ T,q)}.
    \end{align*}

    $(a)\implies (c)$: Note that $(f+g)^R=f^R+g^R$ for functions $f,g:\F_2^\Z\to\F_2^\Z$ because $R$ is linear. It follows that
    \begin{align*}
        \gamma_j^{(c^R,-q)}\circ \left(\id+\sum_{i=1}^s \alpha_i \gamma_i^{(c^R,-q)}\right)
            &= (\gamma_j^{(c,q)})^R\circ \left(\id+\sum_{i=1}^s \alpha_i (\gamma_i^{(c,q)})^R\right) \\
            &= (\gamma_j^{(c,q)})^R\circ \left(\gamma_0^{(c,q)}+\sum_{i=1}^s \alpha_i \gamma_i^{(c,q)}\right)^R \\
            &= \left(\gamma_j^{(c,q)}\circ \left(\gamma_0^{(c,q)}+\sum_{i=1}^s \alpha_i \gamma_i^{(c,q)}\right)\right)^R \\
            &= \left(\sum_{i=0}^s \alpha_i \gamma_{i+j}^{(c,q)}\right)^R \\
            &= \sum_{i=0}^s \alpha_i (\gamma_{i+j}^{(c,q)})^R \\
            &= \sum_{i=0}^s \alpha_i \gamma_{i+j}^{(c^R,-q)}.
    \end{align*}
    The statement for $(c^R\circ T,-q)$ follows.
\end{proof}

\subsection{The main result}

We are now ready to state the main result of this thesis.

\begin{theorem}\label{thm:qc_quasi_conserved_iff_polynomial_composition_property}
    Let $c$ be a nondegenerate landscape. Then the following are equivalent:
    \begin{enumerate}
        \item $(c,q)$ is quasi-conserved.
        \item $(c,q)$ satisfies the polynomial composition property.
    \end{enumerate}
\end{theorem}

\begin{example}\label{ex:qc_wpcp_for_right_sided_landscape}
    Let $c=*001$ which is quasi-conserved with special index $3$. The simplest case of the polynomial composition property states that 
    \begin{equation*}
        \gamma_1 \circ (\gamma_0+\gamma_1) = \gamma_1+\gamma_2.
    \end{equation*}
    Here $\gamma_0=\id, \gamma_1=c$ and $\gamma_2=c^{(3)}\odot S^3c=*00\mt001$. We are going to compute $c\circ(\id+c)$ and show that this indeed equals $c+\gamma_2$. We get
    \begin{align*}
        c\circ (\id+c)
            &= \left((\one+S)\odot(\one+S^2)\odot S^3\right)\circ (\id+c) \\
            &= (\one+S(\id+c))\odot(\one+S^2 (\id+c))\odot S^3(\id+c) \\
            &= (\one+S + Sc)\odot (\one+S^2 + S^2c)\odot(S^3 + S^3c) \\
            &= (*0 + *\mt 001) \odot (*\mt 0 + *\mt\,\mt001) \odot (*\mt\,\mt1 + *\mt\,\mt\,\mt001) \\
            &=\term{
                * & 0 \\
                * & \mt & 0 \\
                * & \mt & \mt & 1 
            }
            + 
            \term{
                * & 0 \\
                * & \mt & 0 \\
                * & \mt & \mt & \mt & 0 & 0 & 1 
            }
            +
            \term{
                * & 0 \\
                * & \mt & \mt & 0 & 0 & 1 \\
                * & \mt & \mt & 1 
            }
            +
            \term{
                * & 0 \\
                * & \mt & \mt & 0 & 0 & 1 \\
                * & \mt & \mt & \mt & 0 & 0 & 1 
            }
            \\[0.5cm]
            &\quad +
            \term{
                * & \mt & 0 & 0 & 1 \\
                * & \mt & 0 \\
                * & \mt & \mt & 1 
            }
            +
            \term{
                * & \mt & 0 & 0 & 1 \\
                * & \mt & 0 \\
                * & \mt & \mt & \mt & 0 & 0 & 1
            }
            +
            \term{
                * & \mt & 0 & 0 & 1 \\
                * & \mt & \mt & 0 & 0 & 1 \\
                * & \mt & \mt & 1 
            }
            +
            \term{
                * & \mt & 0 & 0 & 1 \\
                * & \mt & \mt & 0 & 0 & 1 \\
                * & \mt & \mt & \mt & 0 & 0 & 1
            } \\
            &= *001 + *00\mt001 = c+\gamma_2,
    \end{align*}
    as claimed.
\end{example}

The preprint \cite{zhao2026bifix} independently proved a weaker version of \cref{thm:qc_quasi_conserved_iff_polynomial_composition_property}. It shows that if $c$ is a right-sided landscape with full support $\supp(c)=\{1,\ldots,m\}$, then $(c,m)$ satisfies the polynomial composition property if and only if $c$ corresponds to a so-called bifix-free string, which is equivalent to being quasi-conserved under these special assumptions. For more details we refer to \cref{rem:right_sided_bifix_free}. 

The results of the recent preprint \cite{lyu2026classshiftinvariant} are also similar to \cref{thm:qc_quasi_conserved_iff_polynomial_composition_property}. The main distinction is that their results are formulated over $\F_2^n$, rather than $\F_2^\Z$. As this work appeared only after the technical part of this thesis had been completed, during the final stages of preparing this thesis, we have not had sufficient time to carry out a detailed comparison of the two results.

We split the proof of \cref{thm:qc_quasi_conserved_iff_polynomial_composition_property} into two sections. \cref{sec:right_sided_landscapes} focuses on right-sided landscapes while \cref{sec:two_sided_landscapes} proves the general case of two-sided landscapes. Although the argument in \cref{sec:two_sided_landscapes} also applies to right-sided landscapes, we treat them separately because they allow for a more intuitive proof, providing useful insight into the general case. Additionally, we are able to prove some stronger claims for right-sided landscapes. 

We finish this section by stating the following lemma, implying that we may assume without loss of generality that $c$ is normed when proving \cref{thm:qc_quasi_conserved_iff_polynomial_composition_property}.

\begin{lemma}\label{lem:qc_may_assume_normed_in_main_theorem}
    Let $c$ be a nondegenerate landscape. Then \cref{thm:qc_quasi_conserved_iff_polynomial_composition_property} is true for $c$ if and only if it is true for any conjugate equivalent landscape. 
\end{lemma}

\begin{proof}
    Combine \cref{lem:qc_conjugate_equivalence_preserves_quasi_conserved} and \cref{lem:qc_conjugate_equivalence_preserves_pcp}.
\end{proof}

\section{Right-sided quasi-conserved landscapes}\label{sec:right_sided_landscapes}

In this section we prove \cref{thm:qc_quasi_conserved_iff_polynomial_composition_property} for the case of right-sided landscapes $c$. While this also follows from the results in \cref{sec:two_sided_landscapes}, we split up this proof because for right-sided landscapes we are able to give more intuitive proofs, using diagrams such as \cref{eq:qc_landscape_depiction_S_ell_c_odot_c}.

By \cref{lem:qc_may_assume_normed_in_main_theorem} we make the following assumption for this section.

\begin{assumption}
    Throughout this section we assume that $c=*a_1\ldots a_m$ is a normed nondegenerate quasi-conserved right-sided landscape with special index $q$, unless mentioned otherwise.
\end{assumption}

Our first observation is that we can determine $q$.

\begin{lemma}\label{lem:right_sided_special_index_is_m}
    For right-sided landscapes it holds that $q=m$.
\end{lemma}

\begin{proof}
    By \cref{lem:qc_equivalent_characterization_conserved_quasi_conserved_with_S_ell_c_odot_c} we only need to check that $S^mc\odot c\neq \zero$. But this is obvious, as $S^m c$ and $c$ have disjoint supports. Alternatively, consider the discussion preceding \cref{def:qc_conserved_and_quasi_conserved_landscapes}.
\end{proof}

The following remark highlights a connection between right-sided landscapes and so-called bifix-free strings \cite{nielsen1973note}. 

\begin{remark}\label{rem:right_sided_bifix_free}
    Set $m=11$ and let $c=*a_1\ldots a_{11}$ be a right-sided landscape with full support $\supp(c)=\{1,\ldots,11\}$, i.e. $a_\ell\in\{0,1\}$ for all $\ell\in\{1,\ldots,11\}$. Consider for example 
    \begin{equation}\label{eq:qc_S_3_c_odot_c_example}
        S^3 c\odot c = 
            \term[2pt]{
                * & \mt & \mt & \mt & a_1 & a_2 & a_3 & a_4 & a_5 & a_6 & a_7    & a_8 & a_9 & a_{10} & a_{11}\\
                * & a_1 & a_2 & a_3 & a_4 & a_5 & a_6 & a_7 & a_8 & a_9 & a_{10} & a_{11}
            }
    \end{equation}
    and suppose that $S^3 c\odot c\neq \zero$, i.e. there is no conflict in \cref{eq:qc_S_3_c_odot_c_example}. It then follows that $a_\ell=a_{\ell+3}$ for all $\ell\in\{1,\ldots,8\}$, i.e. $c = *a_1 a_2 a_3 a_1 a_2 a_3 a_1 a_2 a_3 a_1 a_2$. The number $3$ is then called a period of the string $a_1\ldots a_{11}$. More generally, there is no conflict in $S^\ell c\odot c$ if and only if $\ell$ is a period of the string $a_1\ldots a_m$. Note that the existence of a period of the string implies that there exists a prefix and a suffix of the string which are equal. Indeed, here $a_1a_2 = a_{10}a_{11}$. The substring $a_1a_2$ is then called a bifix of $a_1\ldots a_{11}$. A string without any proper bifix is called a bifix-free string \cite{nielsen1973note}. By \cref{lem:qc_equivalent_characterization_conserved_quasi_conserved_with_S_ell_c_odot_c} combined with \cref{lem:qc_x_matches_at_pos_i_plus_ell_iff_S_ell_c_odot_c_equals_zero} the following statements are equivalent:
    \begin{enumerate}
        \item The landscape $c=*a_1\ldots a_m$ with full support is quasi-conserved with special index $m$.
        \item The binary string $a_1\ldots a_m$ is bifix-free.
    \end{enumerate}
    The paper \cite{zhao2026bifix} studies exactly the functions given by these landscapes and independently of our work proves \cref{thm:qc_quasi_conserved_iff_polynomial_composition_property} for this specific class of landscapes. 
    
    We want to emphasize that our approach is more general, even if we restrict ourselves solely to right-sided landscapes. Indeed, we additionally allow the symbol $\mt$ as part of our string, which gives exponentially many more landscapes to consider. Furthermore, the previous argument does not apply. Even if we know that there is no conflict in \cref{eq:qc_S_3_c_odot_c_example}, then we cannot for example conclude that $a_1=a_4$, because it is possible that $a_1=\mt$ or $a_4=\mt$. In fact, now neither implication is true:
    \begin{itemize}
        \item $c=*10\mt1$ is quasi-conserved but the ternary string $10\mt1$ is not bifix-free ($1$ is a bifix).
        \item $c=*00\mt\,\mt\,\mt1$ is not quasi-conserved but the ternary string $00\mt\,\mt\,\mt1$ is bifix-free. 
    \end{itemize}
    The reason for this is that on one hand, we only need $S^\ell c\odot c=\zero$ for $\ell\in\supp(c)$ instead of all $\ell\in\{1,\ldots,m-1\}$ so that certain periods/bifixes are allowed and on the other hand the symbol $\mt$ cannot lead to a conflict but may prevent a period/bifix. 

    In \cite[Theorem 1]{nielsen1973note} the number of bifix-free strings was studied. It was shown that asymptotically $27\%$ of all binary strings and $56\%$ of all ternary strings are bifix-free. The precise count of all quasi-conserved landscapes (right-sided and two-sided) is an interesting problem which we leave open. We only mention that the sequence of the number of right-sided quasi-conserved landscapes of length $m\geq 1$ starts with the values 
    \begin{equation*}
        2, 4, 6, 14, 26, 74, 178, 534, 1468, 4412, 12508, 37652, 110418, 331994, \ldots
    \end{equation*}
    and does not seem to satisfy any simple type of recurrence. 
\end{remark}

\begin{problem}
    Let $k,m\geq 0$ be non-negative integers. Determine the number $q(k,m)$ of quasi-conserved landscapes $c=a_{-k}\ldots a_m$. 
\end{problem}

Next we want to study the functions $\gamma_j$. Let us first look at some examples.

\begin{example}\label{ex:right_sided_gammas_for_right_sided_landscape}
    \begin{enumerate}
        \item Let $c=*0111$. Then $c^{(4)}=*011\mt$ and
        \begin{align*}
            \gamma_0 &= \id = \hat 1 \\
            \gamma_1 &= c = *0111 \\
            \gamma_2 &= c^{(4)}\odot S^4(*0111) = *011\mt\odot *\mt\,\mt\,\mt\,\mt 0111 = *011\mt0111 \\
            \gamma_3 &= c^{(4)}\odot S^4(*011\mt0111) = *011\mt\odot *\mt\,\mt\,\mt\,\mt 011\mt0111 = *011\mt011\mt0111.
        \end{align*}
        \item Let $c=*10\mt1$, then $c^{(4)} = *10\mt\,\mt$ and hence
        \begin{align*}
            \gamma_0 &= \id = \hat 1 \\
            \gamma_1 &= c = *10\mt1 \\
            \gamma_2 &= c^{(4)}\odot S^4(*10\mt1) = *10\mt\,\mt\odot *\mt\,\mt\,\mt\,\mt 10\mt1 = *10\mt\,\mt10\mt1 \\
            \gamma_3 &= c^{(4)}\odot S^4(*10\mt\,\mt10\mt1) = *10\mt\,\mt\odot *\mt\,\mt\,\mt\,\mt 10\mt\,\mt10\mt1 = *10\mt\,\mt10\mt\,\mt10\mt1.
        \end{align*}
    \end{enumerate}
\end{example}

From \cref{def:qc_gamma_functions} we quickly get the following formulas for the functions $\gamma_j$ via induction. This also holds for landscapes that are not quasi-conserved.

\begin{lemma}\label{lem:right_sided_landscape_formulas_for_gamma_j}
    It holds that
    \begin{equation*}
        \gamma_j = \cm \odot S^m \cm\odot S^{2m} \cm \odot\cdots\odot S^{(j-2)m}\cm\odot S^{(j-1)m}c.
    \end{equation*}
    In particular, the landscape of $\gamma_j$ has the form
    \begin{equation*}
        \cstarvector \cstarvector \ldots \cstarvector \cvector,
    \end{equation*}
    i.e. the landscape of $\cm$ is repeated $j-1$ times and then followed by the landscape $c$. More explicitly, 
    \begin{equation*}
        \gamma_j = *a_1\ldots a_{m-1}\mt a_1\ldots a_{m-1}\mt \ldots a_1\ldots a_{m-1}\mt a_1\ldots a_{m-1} a_m.
    \end{equation*}
    Furthermore, if $i\geq j$, then
    \begin{equation*}
        \gamma_i = \cm\odot S^m \cm\odot\cdots\odot S^{(j-2)m}\cm \odot S^{(j-1)m} \gamma_{i-j+1}.
    \end{equation*}
\end{lemma}

The following lemma quickly follows. 

\begin{lemma}\label{lem:right_sided_gammas_never_zero}
    Let $c$ be a right-sided landscape. Then $\gamma_j:\F_2^\Z\to\F_2^\Z$ is nonzero for all $j\geq 0$.
\end{lemma}

\begin{proof}
    The claim follows immediately from \cref{lem:right_sided_landscape_formulas_for_gamma_j} as
    \begin{equation*}
        \gamma_j = \cm \odot S^m \cm\odot S^{2m} \cm \odot\cdots\odot S^{(j-2)m}\cm\odot S^{(j-1)m}c
    \end{equation*}
    and the supports of the factors are disjoint.
\end{proof}

We note that \cref{lem:right_sided_gammas_never_zero} does not necessarily hold for two-sided landscapes or when we consider $\gamma_j$ as a function $\F_2^n\to\F_2^n$.

From here on we want to show that if $(c,q)$ is quasi-conserved, then it satisfies the polynomial composition property. The proof strategy here is similar in spirit to that of \cite{kriepke2024algebraic}. When computing the expression
\begin{equation*}
    \gamma_j \circ \left(\gamma_0+\sum_{i=1}^s \alpha_i\gamma_i\right)
\end{equation*}
we obtain many terms containing products of the form $S^\rho \gamma_s\odot \gamma_t$ for certain parameters $\rho,s,t$. Our goal is to show that most of these products vanish. We start by showing the following result. Although we only need a few cases of this lemma, we state and prove the general version for completeness and intuition. Since the corresponding general statement for two-sided landscapes is overly technical and will not be needed for our purposes, we do not state it later.

\begin{lemma}\label{lem:right_sided_equivalent_conditions_product_gammas_zero}
    Let $\rho\in\{1,\ldots,m-1\}$ and $i,j\geq 1$. Then $S^\rho \gamma_i \odot \gamma_j=\zero$ if and only if at least one of the following conditions is satisfied:
    \begin{enumerate}
        \item $S^\rho c\odot \cm =\zero$ and $j\geq 2$.
        \item $S^{m-\rho}c\odot \cm=\zero$ and $i,j\geq 2$.
        \item $S^{m-\rho}c\odot c=\zero$ and $i<j$.
        \item $S^\rho c\odot c=\zero$ and $i\geq j$.
    \end{enumerate}
\end{lemma}

\begin{proof}
    \hbadness=10000
    There are the two cases $i<j$ and $i\geq j$, which we illustrate with examples. First, let $c=*a_1a_2a_3a_4a_5, \rho=2, i=2, j=3$. Then
    \begin{equation*}
        S^\rho \gamma_i \odot \gamma_j = \niceterm[2pt]{
            * & \mt & \mt & a_1 & a_2 & a_3 & a_4 & \mt & a_1 & a_2 & a_3 & a_4 & a_5 \\
            * & a_1 & a_2 & a_3 & a_4 & \mt & a_1 & a_2 & a_3 & a_4 & \mt & a_1 & a_2 & a_3 & a_4 & a_5
        }{
            \UnderBrace[shorten,yshift=3pt]{2-4}{2-6}{(a)}
            \UnderBrace[shorten,yshift=3pt]{2-7}{2-8}{(b)}
            \UnderBrace[shorten,yshift=3pt]{2-9}{2-11}{(a)}
            \UnderBrace[shorten,yshift=3pt]{2-12}{2-13}{(c)}
        }.
    \end{equation*}\\
    
    Similarly, for $i=3$ and $j=2$,
    \begin{equation*}
        S^\rho \gamma_i \odot \gamma_j = \niceterm[2pt]{
            * & \mt & \mt & a_1 & a_2 & a_3 & a_4 & \mt & a_1 & a_2 & a_3 & a_4 & \mt & a_1 & a_2 & a_3 & a_4 & a_5\\
            * & a_1 & a_2 & a_3 & a_4 & \mt & a_1 & a_2 & a_3 & a_4 & a_5
        }{
            \UnderBrace[shorten,yshift=3pt]{2-4}{2-6}{(a)}
            \UnderBrace[shorten,yshift=3pt]{2-7}{2-8}{(b)}
            \UnderBrace[shorten,yshift=3pt]{2-9}{2-11}{(d)}
        }.
    \end{equation*}\\

    It is clear that $S^\rho \gamma_i \odot \gamma_j = \zero$ if and only if one of the indicated blocks $(a)-(d)$ has a conflict, i.e. a column with a $0$ and a $1$. In general the possible conflicts that could occur are the following:
    \begin{enumerate}
        \item The beginning of one of the landscapes $\cm$ or $c$ in $\gamma_i$ is aligned with the end of one of the landscapes $\cm$ in $\gamma_j$. This happens if and only if $j\geq 2$ and
        \begin{equation*}
            S^\rho c\odot \cm = 
                \term[2pt]{
                    * & \mt & \mt & a_1 & a_2 & a_3 & a_4 & a_5    \\
                    * & a_1 & a_2 & a_3 & a_4 & \mt
                }
                = \zero.
        \end{equation*}
        \item The end of one of the landscapes $\cm$ in $\gamma_i$ is aligned with the beginning of a landscape $\cm$ or $c$ in $\gamma_j$. This happens if and only if $i,j\geq 2$ and
        \begin{equation*}
            S^{m-\rho} c\odot \cm = 
                \term[2pt]{
                    * & \mt & \mt & \mt & a_1 & a_2 & a_3 & a_4 & a_5    \\
                    * & a_1 & a_2 & a_3 & a_4 & \mt
                }
                =\zero.
        \end{equation*}
        \item The end of the landscape $c$ in $\gamma_i$ is aligned with the beginning of a landscape $\cm$ or $c$ in $\gamma_j$. This happens if and only if $i<j$ and 
        \begin{equation*}
            S^{m-\rho}c\odot c = 
                \term[2pt]{
                    * & \mt & \mt & \mt & a_1 & a_2 & a_3 & a_4 & a_5    \\
                    * & a_1 & a_2 & a_3 & a_4 & a_5
                }
                =\zero.
        \end{equation*}
        \item The beginning of one of the landscapes $\cm$ or $c$ in $\gamma_i$ is aligned with the end of the landscape $c$ in $\gamma_j$. This happens if and only if $i\geq j$ and
        \begin{equation*}
            S^\rho c\odot c = 
                \term[2pt]{
                    * & \mt & \mt & a_1 & a_2 & a_3 & a_4 & a_5    \\
                    * & a_1 & a_2 & a_3 & a_4 & a_5
                }
                = \zero.
        \end{equation*}
    \qedhere
    \end{enumerate}
\end{proof}

\begin{remark}\label{rem:right_sided_no_diff_between_S_ell_cm_odot_c_and_S_ell_c_odot_c}
    There are two important observations which follow from the proof. The first is that because $q=m$ is the length of the landscapes since $c$ is right-sided, the only relevant products are of the form $S^\rho \hat c\odot \tilde c$ or $S^{m-\rho}\hat c\odot \tilde c$ with $\hat c, \tilde c\in\{c,\cm\}$, i.e. with the offsets $\rho$ and $m-\rho$. 

    The second observation is that $S^\rho \cm \odot \tilde c=\zero$ if and only if $S^\rho c\odot \tilde c=\zero$. 
    Indeed, 
    \begin{equation*}
        S^\rho \cm\odot \tilde c=
            \term[2pt]{
                * & \mt & \ldots & \mt    & a_1        & \ldots & a_{m-\rho-1} & a_{m-\rho} & \ldots & a_{m-1} & \mt \\
                * & a_1 & \ldots & a_\rho & a_{\rho+1} & \ldots & a_{m-1}      & b_m
            }
    \end{equation*}
    has a conflict if and only if 
    \begin{equation*}
        S^\rho c\odot \tilde c=
            \term[2pt]{
                * & \mt & \ldots & \mt    & a_1        & \ldots & a_{m-\rho-1} & a_{m-\rho} & \ldots & a_{m-1} & a_m \\
                * & a_1 & \ldots & a_\rho & a_{\rho+1} & \ldots & a_{m-1}      & b_m
            }
    \end{equation*}
    has a conflict, where $b_m\in\{a_m,\mt\}$.

    Both of these facts are not true for two-sided landscapes.
\end{remark}

The identities $S^{m-\rho}c\odot c=\zero$ and $S^\rho c\odot c=\zero$ from conditions $(c)$ and $(d)$ in \cref{lem:right_sided_equivalent_conditions_product_gammas_zero} are satisfied if $m-\rho\in\supp(c)\setminus\{m\}$ or $\rho\in\supp(c)\setminus\{m\}$, respectively, due to \cref{lem:qc_equivalent_characterization_conserved_quasi_conserved_with_S_ell_c_odot_c}. For conditions $(a)$ and $(b)$ we need the following lemma:

\begin{lemma}\label{lem:right_sided_S_ell_c_odot_cm_equals_zero}
    Let $\ell\in\supp(c)\setminus\{m\}$. Then $S^\ell c\odot \cm=\zero$ or $m-\ell\in\supp(c)$. 
\end{lemma}

\begin{proof}
    By \cref{lem:qc_equivalent_characterization_conserved_quasi_conserved_with_S_ell_c_odot_c} we know that $S^\ell c\odot c=\zero$, so there is at least one column in
    \begin{equation*}
        S^\ell c\odot c=
            \term[2pt]{
                * & \mt & \ldots & \mt    & a_1        & \ldots & a_{m-\ell} & \ldots & a_m\\
                * & a_1 & \ldots & a_\ell & a_{\ell+1} & \ldots & a_m
            }
    \end{equation*}
    with a conflict, i.e. a column containing both a $0$ and a $1$. If this column is not $(a_{m-\ell},a_m)$, then
    \begin{equation*}
        S^\ell c\odot \cm=
            \term[2pt]{
                * & \mt & \ldots & \mt    & a_1        & \ldots & a_{m-\ell-1} & a_{m-\ell} & \ldots & a_m\\
                * & a_1 & \ldots & a_\ell & a_{\ell+1} & \ldots & a_{m-1}      & \mt 
            }
    \end{equation*}
    also has a conflict, implying $S^\ell c\odot \cm=\zero$. Otherwise, we need $a_{m-\ell}\in\{0,1\}$ and therefore $m-\ell\in\supp(c)$. 
\end{proof}

We obtain the following lemma. We exclude the cases $i=0, j\geq 0$ and $i=j=1$, since the equality $S^{m-\ell}\gamma_i \odot \gamma_j=\zero$ does not necessarily hold in these cases. 

\begin{lemma}\label{lem:right_sided_S_m_minus_ell_gamma_i_odot_gamma_j_equals_zero}
    Let $\ell\in\supp(c)\setminus\{m\}$ and $i\geq 1, j\geq 2$. Then $S^{m-\ell}\gamma_i \odot \gamma_j=\zero$.
\end{lemma}

\begin{proof}
    Set $\rho=m-\ell$. By \cref{lem:right_sided_equivalent_conditions_product_gammas_zero} we need to show that at least one of the following conditions is true:
    \begin{enumerate}
        \item $S^{m-\ell} c\odot \cm =\zero$ and $j\geq 2$.
        \item $S^\ell c\odot \cm=\zero$ and $i,j\geq 2$.
        \item $S^\ell c\odot c=\zero$ and $i<j$.
        \item $S^{m-\ell} c\odot c=\zero$ and $i\geq j$.
    \end{enumerate}
    By \cref{lem:qc_equivalent_characterization_conserved_quasi_conserved_with_S_ell_c_odot_c} we know that $S^\ell c\odot c=\zero$, so if $i<j$ we have shown condition $(c)$.

    So now let $i\geq j$. Note that this implies $i\geq 2$. By \cref{lem:right_sided_S_ell_c_odot_cm_equals_zero} we know that condition $(b)$ holds or $m-\ell\in\supp(c)$. In the latter case $S^{m-\ell}c\odot c=\zero$ by another application of \cref{lem:qc_equivalent_characterization_conserved_quasi_conserved_with_S_ell_c_odot_c}, which gives us $(d)$. 
\end{proof}

This brings us to the main technical lemma of our proof.

\begin{lemma}\label{lem:right_sided_c^mg_odot_S_m_f}
    Let $\alpha_0,\ldots,\alpha_s, \beta_1,\ldots,\beta_t\in\F_2$. Then 
    \begin{equation*}
        \cm\left(\gamma_0+\sum_{j=1}^t \beta_j \gamma_j\right)\odot S^m \left(\sum_{i=0}^s \alpha_i \gamma_i\right) = \sum_{i=0}^s \alpha_i \gamma_{i+1}.
    \end{equation*} 
\end{lemma}

\begin{proof}
    Write $f=\sum_{i=0}^s \alpha_i \gamma_i$ and $g=\gamma_0+\sum_{j=1}^t \beta_j \gamma_j$. Consider first
    \begin{equation*}
        \cm(g) = \left(\bigodot_{\ell\in\supp(c)\setminus\{m\}} A_\ell\right)\circ g = \bigodot_{\ell\in\supp(c)\setminus\{m\}} A_\ell\circ g
    \end{equation*}
    with
    \begin{equation*}
        A_\ell =
            \begin{cases}
                S^\ell              & : a_\ell = 1   \\
                \one + S^\ell       & : a_\ell = 0.
            \end{cases}
    \end{equation*}
    Assume $a_\ell = 1$. Then 
    \begin{equation*}
        A_\ell \circ g = S^\ell\left(\gamma_0+\sum_{j=1}^t \beta_j \gamma_j\right) = S^\ell + S^\ell\left(\sum_{j=1}^t \beta_j \gamma_j\right) = A_\ell + \sum_{j=1}^t \beta_j S^\ell\gamma_j.
    \end{equation*}
    Similarly, if $a_\ell=0$ then
    \begin{equation*}
        A_\ell \circ g = \one+S^\ell\left(\gamma_0+\sum_{j=1}^t \beta_j \gamma_j\right) = \one+S^\ell + S^\ell\left(\sum_{j=1}^t \beta_j \gamma_j\right) = A_\ell + \sum_{j=1}^t \beta_j S^\ell\gamma_j.
    \end{equation*}
    We then obtain 
    \begin{align}\label{eq:right_sided_before_multiplying_out_c^mg_odot_S_m_f}
        \cm(g)\odot S^m f
            &= \left(\bigodot_{\ell\in\supp(c)\setminus\{m\}} A_\ell\circ g\right)\odot S^m f \notag\\
            &= \bigodot_{\ell\in\supp(c)\setminus\{m\}} \left(A_\ell + \sum_{j=1}^t \beta_jS^\ell\gamma_j\right)\odot \left(\sum_{i=0}^s \alpha_i S^m \gamma_i\right)
    \end{align}
    We are now going to multiply this out. By choosing the terms $A_\ell$ of the first factors we get
    \begin{align}\label{eq:right_sided_surviving_terms_c^mg_odot_S_m_f}
        \left(\bigodot_{\ell\in\supp(c)\setminus\{m\}} A_\ell \right)\odot \left(\sum_{i=0}^s \alpha_i S^m \gamma_i\right) 
            &= \cm \odot \left(\sum_{i=0}^s \alpha_i S^m \gamma_i\right) \notag \\
            &= \sum_{i=0}^s \alpha_i \cm \odot S^m \gamma_i \notag \\
            &= \sum_{i=0}^s \alpha_i \gamma_{i+1}
    \end{align}
    by the definition of the functions $\gamma_i$ and using that $c$ is normed.

    Let $P$ be any other term that we get from multiplying out. Then $P$ has the form
    \begin{equation*}
        P = \bigodot_{\ell\in N_1} A_\ell \odot \bigodot_{\ell\in N_2} S^\ell \gamma_{j_\ell} \odot S^m \gamma_i
    \end{equation*}
    where $N_1\dot\cup N_2=\supp(c)\setminus\{m\}$ is a disjoint union, $N_2$ is non-empty, $j_\ell\geq 1$ for all $\ell\in N_2$ and $i\geq 0$. Assume first that $i\geq 1$ and there exists some $\ell\in N_2$ such that $j_\ell\geq 2$. In this case we can apply \cref{lem:right_sided_S_m_minus_ell_gamma_i_odot_gamma_j_equals_zero} to obtain
    \begin{equation*}
        S^\ell \gamma_{j_\ell} \odot S^m \gamma_i = S^\ell (\gamma_{j_\ell}\odot S^{m-\ell}\gamma_i) = \zero,
    \end{equation*}
    which implies $P=\zero$. We are left with the case that $j_\ell=1$ for all $\ell\in N_2$ or $i=0$. Take $\mu\in N_2$ maximal, then
    \begin{equation*}
        P = P' \odot S^\mu \gamma_1 \odot \bigodot_{\ell\in N_3} A_\ell \odot S^m \gamma_i,
    \end{equation*}
    where $N_3 = \{\ell \in N_1: \ell > \mu\}$ and $P'$ is the product of the remaining factors, i.e.
    \begin{equation*}
        P' = \bigodot_{\ell\in N_1\setminus N_3}A_\ell \odot \bigodot_{\ell\in N_2\setminus\{\mu\}} S^\ell \gamma_1.
    \end{equation*}
    Let first $i\geq 1$. Observe that
    \begin{align*}
        S^\mu \gamma_1 \odot \bigodot_{\ell\in N_3} A_\ell \odot S^m \gamma_i 
            &= \term[2pt]{
                * & \mt    & \cdots & \mt    & a_1        & a_2    & \cdots & a_{m-\mu+1} & a_{m-\mu} & a_{m-\mu+1} & \cdots & a_m \\
                * & \mt    & \cdots & \mt    & a_{\mu+1} \\
                * & \vdots &        & \vdots &            & \ddots \\
                * & \mt    & \cdots & \cdots &            & \cdots & \mt    & a_{m-1}         \\
                * & \mt    & \cdots & \cdots &            & \cdots & \mt    & \mt          & \mt        & a_1          & \cdots & a_\mu & \cdots 
            } \\[0.5cm]
            &=\term[2pt]{
                * & \mt    & \cdots & \mt    & a_1       & a_2       & \cdots & a_{m-\mu+1} & a_{m-\mu} & a_{m-\mu+1} & \cdots & a_m \\
                * & \mt    & \cdots & \mt    & a_{\mu+1} & a_{\mu+2} & \cdots & a_{m-1}     & \mt       & a_1          & \cdots & a_\mu & \cdots 
            } 
    \end{align*}
    has a conflict if and only if 
    \begin{equation*}
        S^\mu c\odot \gamma_2 
            = \term[2pt]{
                * & \mt & \cdots & \mt    & a_1       & a_2       & \cdots & a_{m-\mu+1} & a_{m-\mu} & a_{m-\mu+1} & \cdots & a_m \\
                * & a_1 & \cdots & a_\mu  & a_{\mu+1} & a_{\mu+2} & \cdots & a_{m-1}     & \mt       & a_1          & \cdots & a_\mu & \cdots & a_m
            } 
    \end{equation*}
    has a conflict, which is equivalent to $S^\mu c\odot \gamma_2=\zero$. Thus, to show that $P=\zero$ it suffices to show that $S^\mu c\odot \gamma_2=\zero$. By \cref{lem:right_sided_S_ell_c_odot_cm_equals_zero} combined with \cref{lem:qc_equivalent_characterization_conserved_quasi_conserved_with_S_ell_c_odot_c} we know that at least one of 
    \begin{enumerate}
        \item $S^\mu c \odot \cm = \zero$
        \item[(c)] $S^{m-\mu} c\odot c=\zero$
    \end{enumerate}
    is true, so that \cref{lem:right_sided_equivalent_conditions_product_gammas_zero} implies $S^\mu c\odot \gamma_2=\zero$. Hence, $P=\zero$ if $i\geq 1$. Suppose now that $i=0$. Then similarly
    \begin{align*}
        S^\mu \gamma_1 \odot \bigodot_{\ell\in N_3} A_\ell \odot S^m \gamma_i 
            &= \term[2pt]{
                * & \mt    & \cdots & \mt    & a_1        & a_2    & \cdots & a_{m-\mu+1} & a_{m-\mu} & a_{m-\mu+1} & \cdots & a_m \\
                * & \mt    & \cdots & \mt    & a_{\mu+1} \\
                * & \vdots &        & \vdots &            & \ddots \\
                * & \mt    & \cdots & \cdots &            & \cdots & \mt    & a_{m-1}         \\
                * & \mt    & \cdots & \cdots &            & \cdots & \mt    & \mt          & a_m       
            } \\[0.5cm]
            &=\term[2pt]{
                * & \mt    & \cdots & \mt    & a_1       & a_2       & \cdots & a_{m-\mu+1} & a_{m-\mu} & a_{m-\mu+1} & \cdots & a_m \\
                * & \mt    & \cdots & \mt    & a_{\mu+1} & a_{\mu+2} & \cdots & a_{m-1}     & a_m   
            } 
    \end{align*}
    which has a conflict if and only if 
    \begin{equation*}
        S^\mu c\odot c 
            = \term[2pt]{
                * & \mt & \cdots & \mt    & a_1       & a_2       & \cdots & a_{m-\mu+1} & a_{m-\mu} & a_{m-\mu+1} & \cdots & a_m \\
                * & a_1 & \cdots & a_\mu  & a_{\mu+1} & a_{\mu+2} & \cdots & a_{m-1}     & a_m      
            } 
    \end{equation*}
    has a conflict, i.e. $S^\mu c\odot c=\zero$, which is the case due to \cref{lem:qc_equivalent_characterization_conserved_quasi_conserved_with_S_ell_c_odot_c}. Thus, $P=\zero$ for $i=0$. 
    Hence, the only terms from \cref{eq:right_sided_before_multiplying_out_c^mg_odot_S_m_f} that are nonzero are those in \cref{eq:right_sided_surviving_terms_c^mg_odot_S_m_f} and
    \begin{equation*}
        \cm(g)\odot S^mf = \sum_{i=0}^s \alpha_i \gamma_{i+1},
    \end{equation*}
    as desired.
\end{proof}

We then obtain the following result, showing the first implication of \cref{thm:qc_quasi_conserved_iff_polynomial_composition_property}.

\begin{lemma}\label{lem:right_sided_quasi_conserved_implies_pcp}
    Let $j\geq 0$ and $\alpha_0,\ldots,\alpha_s\in\F_2$ with $\alpha_0=1$. Then
    \begin{equation*}
        \gamma_j\circ \left(\gamma_0+\sum_{i=1}^s \alpha_i\gamma_i\right) = \sum_{i=0}^s \alpha_i \gamma_{i+j}.
    \end{equation*}
    In other words, $(c,q)$ satisfies the polynomial composition property.
\end{lemma}

\begin{proof}
    We proceed by induction on $j$. The base case $j=0$ is clear, as $\gamma_0=\id$. Note that since $c$ is normed we can use the recurrence $\gamma_{j+1}=\cm \odot S^m \gamma_j$ starting from $j=0$. Write $g=\gamma_0+\sum_{i=1}^s \alpha_i\gamma_i$. Then
    \begin{align*}
        \gamma_{j+1} \circ \left(\gamma_0+\sum_{i=1}^s \alpha_i\gamma_i\right)
            &= (\cm \odot S^m \gamma_j)\circ \left(\gamma_0+\sum_{i=1}^s \alpha_i\gamma_i\right) \\
            &= \cm(g) \odot S^m \gamma_j\circ \left(\gamma_0+\sum_{i=1}^s \alpha_i\gamma_i\right) \\
            &= \cm(g) \odot S^m \left(\sum_{i=0}^s \alpha_i\gamma_{i+j}\right) \\
            &= \sum_{i=0}^s \alpha_i \gamma_{i+j+1},
    \end{align*}
    where in the third line we used the induction hypothesis and the last equality follows from \cref{lem:right_sided_c^mg_odot_S_m_f}.
\end{proof}

The remainder of this section is dedicated to proving the converse implication of \cref{thm:qc_quasi_conserved_iff_polynomial_composition_property}, i.e. that if $(c,q)$ satisfies the polynomial composition property, then $c$ is quasi-conserved. We accomplish this by showing that if $c$ is not quasi-conserved, then
\begin{equation*}
    c\circ (\id+c) \neq c+\gamma_2.
\end{equation*}
To show that these two functions are different we find a vector $x\in\F_2^\Z$ with $f(x)\neq g(x)$, where $f=c\circ (\id+c)$ and $g=c+\gamma_2$. 

Let $c$ be a right-sided landscape which is not quasi-conserved. By \cref{def:qc_conserved_and_quasi_conserved_landscapes} and the discussion preceding it we know that there exist at least one index $\ell\in\supp(c)\setminus\{m\}$ and a vector $x\in\F_2^\Z$ which matches $c$ at positions $0$ and $\ell$. Consider the diagram
\begin{equation*}
    \term[2pt]{
        x_0 & \ldots & x_\ell & \ldots & x_m        & \ldots & x_{m+\ell} & \ldots & x_{2m} \\
        *   & \ldots & a_\ell & \ldots & a_m  \\
            &        & *      & \ldots & a_{m-\ell} & \ldots & a_m   \\
            &        &        &        & *          & \ldots & a_\ell     & \ldots & a_m, 
    }
\end{equation*}
which describes the entries of $x$ around $0$. Now change $x_{2m}$ so that $x$ does not match $c$ at position $m$, i.e. let $x_{2m}=a_m+1$. We are now going to compute $f(x)_0$ and $g(x)_0$. 

Set $y=(\id+c)(x)$. Then $y_\ell=x_\ell+1$ because $x$ matches $c$ at position $\ell$. This implies that $y$ does not match $c$ at position $0$, so that $f(x)_0=c(y)_0=0$. Next, set $u=c(x)$ and $w=\gamma_2(x)$. Then $u_0=1$ because $x$ matches $c$ at position $0$. To determine $w_0$ we need to check whether $x$ matches $\gamma_2=\cm\odot S^m c$ at $0$. We see from 
\begin{equation*}
    \term[2pt]{
        x_0 & \ldots & x_\ell & \ldots & x_{m-1} & x_m & x_{m+1} & \ldots & x_{2m} \\
        *   & \ldots & a_\ell & \ldots & a_{m-1} & \mt & a_1     & \ldots & a_m
        }
\end{equation*}
that this is not the case, as $x_{2m}=a_m+1$. Hence, $w_0=0$. Therefore,
\begin{equation*}
    f(x)_0 = 0 \neq 1 = u_0 + w_0 = g(x)_0,
\end{equation*}
implying
\begin{equation*}
    c\circ(\id+c) \neq c+\gamma_2.
\end{equation*}
We have proven:

\begin{lemma}\label{lem:right_sided_pcp_implies_quasi_conserved}
    Let $(c,q)$ satisfy the polynomial composition property. Then $(c,q)$ is quasi-conserved.
\end{lemma}

We arrive at the following theorem.

\begin{theorem}
    Let $c$ be a one-sided landscape. Then \cref{thm:qc_quasi_conserved_iff_polynomial_composition_property} is true for $c$.
\end{theorem}

\begin{proof}
    If $c$ is right-sided then this follows from \cref{lem:right_sided_quasi_conserved_implies_pcp} and \cref{lem:right_sided_pcp_implies_quasi_conserved}. Since any left-sided landscape is conjugate equivalent to a right-sided one, this also proves \cref{thm:qc_quasi_conserved_iff_polynomial_composition_property} for left-sided landscapes.
\end{proof}

\section{Two-sided quasi-conserved landscapes}\label{sec:two_sided_landscapes}

We now prove \cref{thm:qc_quasi_conserved_iff_polynomial_composition_property} for the case of two-sided landscapes. By \cref{lem:qc_may_assume_normed_in_main_theorem} we may assume without loss of generality the following.

\begin{assumption}
    Throughout this section we assume that $c=a_{-k}\ldots a_m$ is a normed nondegenerate quasi-conserved two-sided landscape with special index $q$, unless mentioned otherwise.
\end{assumption}

The following example gives the functions $\gamma_j:\F_2^\Z\to\F_2^\Z$ defined in \cref{def:qc_gamma_functions}. It shows that \cref{lem:right_sided_landscape_formulas_for_gamma_j,lem:right_sided_gammas_never_zero} are not true for two-sided landscapes.

\begin{example}\label{ex:two_sided_gammas}
    \begin{enumerate}
        \item Let $c=0*\mt10$ as in \cref{ex:qc_conserved_and_quasi_conserved_landscapes} and \cref{ex:two_sided_quasi_conserved_not_preserved_under_shift}.
        Then $q=2$ and $c^{(2)} = 0*\mt\,\mt 0$. Then
        \begin{align*}
            \gamma_1 &= c^{(2)} \odot S^2 \gamma_0 = 0*\mt 10 = c, \\
            \gamma_2 &= c^{(2)} \odot S^2 \gamma_1 = 0*\mt\,\mt 0 \odot S^2(0*\mt 10) = 0*\mt\,\mt 0 \odot *0\mt\,\mt 10 = 0*0\mt010, \\
            \gamma_3 &= c^{(2)} \odot S^2 \gamma_2 = 0*\mt\,\mt 0 \odot S^2(0*0\mt010) = 0*\mt\,\mt 0 \odot *0\mt0\mt010 = 0*0\mt0\mt010.
        \end{align*}
        \item Let $c=0\mt*1\mt0$ which is quasi-conserved with special index $q=1$. We compute
        \begin{align*}
            \gamma_2 &= c^{(1)}\odot Sc = 
                \term[2pt]{
                    0 & \mt & * & \mt & \mt & 0 \\
                      &   0 & * & \mt & 1   & \mt & 0
                }
                = 00*\mt100, \\
            \gamma_3 &= c^{(1)}\odot S\gamma_2 = 
                \term[2pt]{
                    0 & \mt & *     & \mt & \mt & 0 \\
                      &   0 & \hat0 & \mt & \mt & 1 & 0 & 0 
                }
                = \zero
        \end{align*}
        and see that $\gamma_3=\zero$. 
    \end{enumerate}
\end{example}

We restate \cref{lem:right_sided_landscape_formulas_for_gamma_j} for two-sided landscapes, which also holds for landscapes which are not quasi-conserved.

\begin{lemma}\label{lem:two_sided_gammas_formulas}
    Let $i\geq j\geq 1$. Then it holds that
    \begin{equation}\label{eq:two_sided_formula_gamma_j}
        \gamma_j = \cq \odot S^q \cq\odot S^{2q} \cq \odot\cdots\odot S^{(j-2)q}\cq\odot S^{(j-1)q}c
    \end{equation}
    and
    \begin{equation}\label{eq:two_sided_recursive_formula_gamma_i_using_gamma_j}
        \gamma_i = \cq\odot S^q \cq\odot\cdots\odot S^{(j-2)q}\cq \odot S^{(j-1)q} \gamma_{i-j+1}.
    \end{equation}
\end{lemma}

\cref{lem:two_sided_gammas_formulas} gives an explanation for why it is possible that $\gamma_j=\zero$ for some $j> 0$. In the case of right-sided landscapes we have $q=m$, so that all the factors in \cref{eq:two_sided_formula_gamma_j} have disjoint support. Now this is not the case and we can for example obtain a product $\cq \odot S^{\zeta q}\cq =\zero$. We study this more closely in \cref{subsec:pcp_computing_min_polys}.

We also need the following lemma which follows immediately from \cref{eq:two_sided_formula_gamma_j}.

\begin{lemma}\label{lem:two_sided_support_gamma_i}
    Let $j\geq 1$ and suppose $\gamma_j\neq 0$. Then
    \begin{equation*}
        \supp(\gamma_j) = \bigcup_{\zeta=0}^{j-1} (\zeta q+\supp(c)\setminus\{q\}) \cup \{jq\}.
    \end{equation*}
\end{lemma}

Our goal from here is to prove an analogous version of \cref{lem:right_sided_c^mg_odot_S_m_f} for two-sided landscapes. For this purpose we need results similar to \cref{lem:right_sided_equivalent_conditions_product_gammas_zero} and \cref{lem:right_sided_S_m_minus_ell_gamma_i_odot_gamma_j_equals_zero}. However, the argument is more involved here. One reason is the following lemma.

\begin{lemma}\label{lem:two_sided_S^ell_c_or_c*_equals_zero}
    Take $\ell\in\supp(c)\setminus\{q\}$. Then:
    \begin{enumerate}
        \item $S^\ell c\odot \cq = \zero$ or $q-\ell\in\supp(c)$.
        \item $S^\ell \cq\odot c = \zero$ or $q+\ell\in\supp(c)$.
        \item $S^\ell c\odot \cq=\zero$ or $S^\ell \cq \odot c=\zero$.
    \end{enumerate}
\end{lemma}

\begin{proof}
    By \cref{lem:qc_equivalent_characterization_conserved_quasi_conserved_with_S_ell_c_odot_c} we know that $S^\ell c\odot c=\zero$, hence in
    \begin{equation}\label{eq:two_sided_S^ell_c_odot_c_in_proof_S^ell_c_or_cq}
        S^\ell c\odot c=
            \term[2pt]{
            \ldots & a_{q-\ell-1} & a_{q-\ell} & a_{q-\ell+1} & \ldots & a_{-1}     & a_0    & a_1        & \ldots & a_{q-1}      & a_q        & a_{q+1}      & \ldots \\
            \ldots & a_{q-1}      & a_q        & a_{q+1}      & \ldots & a_{\ell-1} & a_\ell & a_{\ell+1} & \ldots & a_{\ell+q-1} & a_{\ell+q} & a_{q+\ell+1} & \ldots 
        }
    \end{equation}
    there is at least one column with a $0$ and a $1$. 

    Suppose that $S^\ell c\odot \cq\neq \zero$. We need to show that $q-\ell\in\supp(c)$ and $S^\ell\cq\odot c=\zero$. The only possible column in
    \begin{equation*}
        S^\ell c\odot \cq = 
            \term[2pt]{
            \ldots & a_{q-\ell-1} & a_{q-\ell} & a_{q-\ell+1} & \ldots & a_{-1}     & a_0    & a_1        & \ldots & a_{q-1}      & a_q        & a_{q+1}      & \ldots \\
            \ldots & a_{q-1}      & \mt        & a_{q+1}      & \ldots & a_{\ell-1} & a_\ell & a_{\ell+1} & \ldots & a_{\ell+q-1} & a_{\ell+q} & a_{q+\ell+1} & \ldots 
        }
    \end{equation*}
    that is different compared to \cref{eq:two_sided_S^ell_c_odot_c_in_proof_S^ell_c_or_cq} is the column $(a_{q-\ell},a_q)$, which is now $(a_{q-\ell},\mt)$. It follows that $(a_{q-\ell},a_q)\in\{(0,1),(1,0)\}$, which in particular means that $a_{q-\ell}\neq \mt$, i.e. $q-\ell\in\supp(c)$. Furthermore, the product $S^\ell \cq\odot c$ also includes the column $(a_{q-\ell},a_q)$, hence $S^\ell \cq\odot c=\zero$.

    The remaining case is analogous. 
\end{proof}

\begin{remark}\label{rem:two_sided_S^ell_c_or_c*_equals_zero_also_for_ell_notin_supp}
    Note that the proof of \cref{lem:two_sided_S^ell_c_or_c*_equals_zero} only uses that $S^\ell c\odot c=\zero$, which may also be the case for some $\ell\notin\supp(c)\setminus\{q\}$. 
\end{remark}

One of the reason right-sided landscapes were easier to deal with is that $q=m$ and $\ell\geq 1$ for all $\ell\in\supp(c)$, so that $q+\ell\notin\supp(c)$ and therefore always $S^\ell \cq \odot c=\zero$, see also \cref{rem:right_sided_no_diff_between_S_ell_cm_odot_c_and_S_ell_c_odot_c}. This is not the case for two-sided landscapes, as the following example shows.

\begin{example}
    Let $c=0*0\mt\,\mt01$ (see \cref{ex:qc_conserved_and_quasi_conserved_landscapes}). Then $(c,5)$ is quasi-conserved. Take $\ell=-1\in\supp(c)$. We see that 
    \begin{equation*}
        S^{-1}\cq\odot c = 
            \term[2pt]{
                0 & \mt & \hat0 & \mt & \mt &   0 & \mt \\
                  & 0   & *     & 0   & \mt & \mt &   0 & 1
            }
            = 00\hat{0}0\mt001 \neq \zero.
    \end{equation*}
\end{example}

\begin{lemma}\label{lem:two_sided_S_ell_gamma_i_odot_c_equals_zero}
    Take $\ell\in\supp(c)\setminus\{q\}$ and $i\geq 1$. Then $S^\ell \gamma_i\odot c=\zero$.
\end{lemma}

\begin{proof}
    We let $\ell\in\supp(c)\setminus\{q\}$ be arbitrary and proceed by induction on $i$. Let $i=1$. Then $\gamma_1=c$ and $S^\ell c\odot c=\zero$ follows from \cref{lem:qc_equivalent_characterization_conserved_quasi_conserved_with_S_ell_c_odot_c}. So now suppose that $S^\ell \gamma_i\odot c=\zero$ for all $\ell\in\supp(c)\setminus\{q\}$. We want to show that $S^\ell \gamma_{i+1}\odot c=\zero$. It holds that $\gamma_{i+1}=\cq\odot S^q \gamma_i$. Hence,
    \begin{equation*}
        S^\ell\gamma_{i+1}\odot c = S^\ell(\cq\odot S^q \gamma_i) \odot c = (S^\ell \cq\odot c) \odot (S^{q+\ell}\gamma_i \odot c),
    \end{equation*}
    using that $c=c\odot c$. By \cref{lem:two_sided_S^ell_c_or_c*_equals_zero} $S^\ell \cq\odot c=\zero$ or $q+\ell\in\supp(c)$. But if $q+\ell\in\supp(c)$, then in particular also $q+\ell\in\supp(c)\setminus\{q\}$ and by the induction hypothesis $S^{q+\ell}\gamma_i\odot c=\zero$. In either case it follows that $S^\ell \gamma_{i+1}\odot c=\zero$, as desired.
\end{proof}

\begin{lemma}\label{lem:two_sided_S_ell_c_odot_gamma_j_equals_zero}
    Take $\ell\in\supp(c)\setminus\{q\}$ and $j\geq 1$. Then $S^\ell c\odot \gamma_j=\zero$.
\end{lemma}

\begin{proof}
    Similarly to the proof of \cref{lem:two_sided_S_ell_gamma_i_odot_c_equals_zero}, we let $\ell\in\supp(c)\setminus\{q\}$ be arbitrary and proceed by induction on $j$. The base case $j=1$ is clear and for the induction step we write
    \begin{equation*}
        S^\ell c\odot \gamma_{j+1} = S^\ell c\odot(\cq\odot S^q \gamma_j) = (S^\ell c\odot \cq) \odot S^\ell(c\odot S^{q-\ell}\gamma_j).
    \end{equation*}
    If $S^\ell c\odot \cq =\zero$ we are done. Otherwise, by \cref{lem:two_sided_S^ell_c_or_c*_equals_zero} we have $q-\ell\in\supp(c)\setminus\{q\}$ and the claim follows from \cref{lem:two_sided_S_ell_gamma_i_odot_c_equals_zero}.
\end{proof}

These results yield the following lemma, which corresponds to \cref{lem:right_sided_equivalent_conditions_product_gammas_zero}. We only give a sufficient condition here, as this is all that is needed for our purposes.

\begin{lemma}\label{lem:two_sided_S_ell_gamma_i_odot_gamma_j_equals_zero}
    Take $\ell\in\supp(c)\setminus\{q\}$ and $i,j\geq 1$ arbitrary. Then $S^\ell \gamma_i\odot \gamma_j=\zero$.
\end{lemma}

\begin{proof}
    First, suppose that $i\geq j$. Using \cref{eq:two_sided_recursive_formula_gamma_i_using_gamma_j} we write $\gamma_i=f\odot S^{(j-1)q}\gamma_{i-j+1}$ and $\gamma_j=g\odot S^{(j-1)q}c$ for some functions $f,g$. Then
    \begin{align*}
        S^\ell \gamma_i\odot \gamma_j &= S^\ell(f\odot S^{(j-1)q}\gamma_{i-j+1})\odot (g\odot S^{(j-1)q}c) \\
            &= S^\ell f\odot g \odot S^{(j-1)q}(S^\ell \gamma_{i-j+1}\odot c)=\zero
    \end{align*}
    by \cref{lem:two_sided_S_ell_gamma_i_odot_c_equals_zero}. If $i<j$ then we write $\gamma_i=f\odot S^{(i-1)q}c$ and $\gamma_j=g\odot S^{(i-1)q}\gamma_{j-i+1}$ and then
    \begin{align*}
        S^\ell\gamma_i\odot \gamma_j &= S^\ell (f\odot S^{(i-1)q}c)\odot (g\odot S^{(i-1)q}\gamma_{j-i+1}) \\
            &= S^\ell f\odot g \odot S^{(i-1)q}\left(S^\ell c\odot \gamma_{j-i+1}\right)=\zero
    \end{align*}
    by \cref{lem:two_sided_S_ell_c_odot_gamma_j_equals_zero}.
\end{proof}

\begin{lemma}\label{lem:two_sided_S_q_minus_ell_gamma_i_odot_gamma_j}
    Take $\ell\in\supp(c)\setminus\{q\}$ and $i\geq 1,j\geq 2$. Then $S^{q-\ell}\gamma_i\odot\gamma_j=\zero$.
\end{lemma}

\begin{proof}
    As $j\geq 2$ we may write $\gamma_j = \cq\odot S^q \gamma_{j-1}$ with $j-1\geq 1$. Then
    \begin{equation*}
        S^{q-\ell}\gamma_i\odot \gamma_j = S^{q-\ell}\gamma_i\odot (\cq\odot S^q\gamma_{j-1}) = \cq \odot S^{q-\ell}\left(\gamma_i\odot S^\ell\gamma_{j-1}\right)=\zero
    \end{equation*}
    by \cref{lem:two_sided_S_ell_gamma_i_odot_gamma_j_equals_zero}.
\end{proof}

In the case of right-sided landscapes, it sufficed to only consider $\ell\in\supp(c)\setminus\{q\}$ in \cref{lem:two_sided_S_ell_gamma_i_odot_gamma_j_equals_zero}, essentially due to the first observation of \cref{rem:right_sided_no_diff_between_S_ell_cm_odot_c_and_S_ell_c_odot_c}. Here we more generally allow $\mu\in\supp(\gamma_j)\setminus\{jq\}$. We assume that $\gamma_j\neq \zero$, otherwise the claim is trivial.

\begin{lemma}\label{lem:two_sided_S_mu_gamma_i_odot_gamma_j_equals_zero}
    Let $i,j\geq 1$ with $\gamma_j\neq \zero$ and take $\mu\in\supp(\gamma_j)$ with $\mu\neq jq$. Then $S^\mu \gamma_i\odot \gamma_j=\zero$.
\end{lemma}

\begin{proof}
    By \cref{lem:two_sided_support_gamma_i} and using $\mu\neq jq$ we know that $\mu = \ell+\zeta q$ for some $\ell\in\supp(c)\setminus\{q\}$ and some $0\leq \zeta\leq j-1$. By \cref{eq:two_sided_recursive_formula_gamma_i_using_gamma_j} applied to $j\geq \zeta+1\geq 1$ we have $\gamma_j = g\odot S^{\zeta q}\gamma_{j-\zeta}$ for some function $g$. Then
    \begin{equation*}
        S^\mu \gamma_i\odot \gamma_j = S^{\ell+\zeta q} \gamma_i \odot g\odot S^{\zeta q} \gamma_{j-\zeta} = g\odot S^{\zeta q} (S^\ell \gamma_i\odot \gamma_{j-\zeta})=\zero
    \end{equation*}
    by \cref{lem:two_sided_S_ell_gamma_i_odot_gamma_j_equals_zero}.
\end{proof}

\begin{remark}
    \cref{lem:two_sided_S_mu_gamma_i_odot_gamma_j_equals_zero} also implies that $(\gamma_j^{(c,q)},jq)$ is quasi-conserved. In particular, the families $\{\gamma_j^{(c,q)}: j\geq 0\}$ and $\{\gamma_j^{(d,r)}: j\geq 0\}$ for different quasi-conserved landscapes $c,d$ are not necessarily disjoint.
\end{remark}

We now state the generalization of \cref{lem:right_sided_c^mg_odot_S_m_f} for the case of two-sided landscapes. The general proof strategy stays the same, but the argument is more involved.

\begin{lemma}\label{lem:two_sided_c^qg_odot_S^q_f}
    Let $\alpha_0,\ldots,\alpha_s, \beta_1,\ldots,\beta_t\in\F_2$. Then 
    \begin{equation*}
        \cq\left(\gamma_0+\sum_{j=1}^t \beta_j \gamma_j\right)\odot S^q \left(\sum_{i=0}^s \alpha_i \gamma_i\right) = \sum_{i=0}^s \alpha_i \gamma_{i+1}.
    \end{equation*} 
\end{lemma}

\begin{proof}
    Write $f=\sum_{i=0}^s \alpha_i \gamma_i$ and $g=\gamma_0+\sum_{j=1}^t \beta_j \gamma_j$. Consider first
    \begin{equation*}
        \cq(g) = \left(\bigodot_{\ell\in\supp(c)\setminus\{q\}} A_\ell\right)\circ g = \bigodot_{\ell\in\supp(c)\setminus\{q\}} A_\ell\circ g
    \end{equation*}
    with
    \begin{equation*}
        A_\ell =
            \begin{cases}
                S^\ell              & : a_\ell = 1   \\
                \one + S^\ell       & : a_\ell = 0.
            \end{cases}
    \end{equation*}
    Analogously to the proof of \cref{lem:right_sided_c^mg_odot_S_m_f} we obtain
    \begin{equation*}
        A_\ell\circ g = A_\ell + \sum_{j=1}^t \beta_j S^\ell\gamma_j
    \end{equation*}
    for every $\ell\in\supp(c)\setminus\{q\}$.
    Therefore,
    \begin{equation}\label{eq:two_sided_before_multiplying_out_c^qg_odot_S_q_f}
        \cq(g)\odot S^qf = \bigodot_{\ell\in\supp(c)\setminus\{q\}} \left[A_\ell + \sum_{j=1}^t \beta_jS^\ell\gamma_j\right] \odot \left(\sum_{i=0}^s \alpha_iS^q\gamma_i\right).
    \end{equation}
    Consider now any term we obtain from multiplying out \cref{eq:two_sided_before_multiplying_out_c^qg_odot_S_q_f}. First we have all terms for which we only take the left summand $A_\ell$ from each of the factors in the big product, so 
    \begin{equation}\label{eq:two_sided_surviving_terms_c^(q)g_odot_S_q_f}
        \left(\bigodot_{\ell\in\supp(c)\setminus\{q\}} A_\ell\right) \odot \sum_{i=0}^s \alpha_i S^q\gamma_i = \cq\odot \sum_{i=0}^s \alpha_i S^q \gamma_i = \sum_{i=0}^s \alpha_i\gamma_{i+1},
    \end{equation}
    using the recurrence $\gamma_{i+1}=\cq\odot S^q \gamma_i$, which is valid for $i=0$ because $c$ is normed.

    Let $P$ be any other term that we get from multiplying out \cref{eq:two_sided_before_multiplying_out_c^qg_odot_S_q_f}. Our goal is to show that $P=\zero$. The term $P$ can be written as 
    \begin{equation*}
        P = \bigodot_{\ell\in N_1} A_\ell \odot\bigodot_{\ell\in N_2} S^\ell \gamma_{j_\ell}\odot S^q \gamma_i,
    \end{equation*}
    where $N_1\dot\cup N_2=\supp(c)\setminus\{q\}$ is a disjoint union, $N_2$ is non-empty and $j_\ell\geq 1$ for all $\ell\in N_2$. Take $\lambda\in N_2$ maximal and assume $\gamma_{j_\lambda}\neq \zero$, otherwise $P=\zero$ is obvious. We have $\gamma_{j_\lambda}=\bigodot_{r\in\supp(\gamma_{j_\lambda})} B_r$ with $B_r\in\{S^r, \one+S^r\}$ for all $r\in\supp(\gamma_{j_\lambda})$. As $S^\lambda \gamma_{j_\lambda}\odot c=\zero$ by \cref{lem:two_sided_S_ell_gamma_i_odot_gamma_j_equals_zero}, there exists $\mu\in\supp(\gamma_{j_\lambda})$ such that $\mu+\lambda\in\supp(c)$ and $S^\lambda B_\mu \odot A_{\mu+\lambda}=\zero$. There are several cases to consider:
    \begin{enumerate}
        \item[(1)] $\mu+\lambda\in N_1$. In this case
        \begin{equation*}
            \bigodot_{\ell\in N_1} A_\ell \odot \bigodot_{\ell\in N_2} S^\ell \gamma_{j_\ell} = A_{\mu+\lambda}\odot S^\lambda \gamma_{j_\lambda}\odot(\ldots) = \zero
        \end{equation*}
        using $S^\lambda B_\mu \odot A_{\mu+\lambda}=\zero$. This shows $P=\zero$.
        \item[(2)] $\mu+\lambda\in N_2$. We have two subcases:
        \begin{enumerate}
            \item[(A)] $\mu=j_\lambda q$. But then $\mu+\lambda=j_\lambda q+\lambda > \lambda$ using that $q>0$ because $c$ is normed. This contradicts the maximality assumption on $\lambda$.
            \item[(B)] $\mu\neq j_\lambda q$. In this case we are able to apply \cref{lem:two_sided_S_mu_gamma_i_odot_gamma_j_equals_zero} to obtain
            \begin{equation*}
                \bigodot_{\ell\in N_2} S^\ell \gamma_{j_\ell} = S^\lambda \gamma_{j_\lambda}\odot S^{\mu+\lambda}\gamma_{j_{\mu+\lambda}} \odot(\ldots) = S^\lambda (\gamma_{j_\lambda}\odot S^\mu \gamma_{j_{\mu+\lambda}})\odot (\ldots) = \zero
            \end{equation*}
            and therefore $P=\zero$.
        \end{enumerate}
        \item[(3)] $\mu+\lambda=q$. We again have some subcases:
        \begin{enumerate}
            \item[(A)] $i=0$. Then $S^q\gamma_i=S^q \id = A_q = A_{\mu+\lambda}$ and then similar to case (1)
            \begin{equation*}
                \bigodot_{\ell\in N_2} S^\ell \gamma_{j_\ell} \odot S^q\gamma_i= S^\lambda \gamma_{j_\lambda}\odot A_{\mu+\lambda}\odot(\ldots) = \zero,
            \end{equation*}
            implying $P=\zero$.
            \item[(B)] $i\geq 1$. We again have two subcases:
            \begin{enumerate}
                \item[(a)] $\mu\in \supp(c)$. Note that in particular $\mu\neq q$ because $\mu+\lambda=q$ and $\lambda\neq 0$, since $\lambda\in\supp(c)\setminus\{q\}$ and $c$ is nondegenerate. Therefore, $\mu\in\supp(c)\setminus\{q\}$. Then
                \begin{equation*}
                    \bigodot_{\ell\in N_2} S^\ell \gamma_{j_\ell} \odot S^q\gamma_i = S^\lambda \gamma_{j_\lambda} \odot S^q\gamma_i\odot (\ldots) = S^\lambda (\gamma_{j_\lambda}\odot S^\mu \gamma_i)\odot(\ldots)=\zero
                \end{equation*}
                by \cref{lem:two_sided_S_ell_gamma_i_odot_gamma_j_equals_zero}. Hence, $P=\zero$.
                \item[(b)] $\mu\notin\supp(c)$. Then in particular $j_\lambda\geq 2$, as otherwise $\gamma_{j_\lambda}=c$ and $\mu\in\supp(\gamma_{j_\lambda})=\supp(c)$. Then 
                \begin{equation*}
                    \bigodot_{\ell\in N_2} S^\ell \gamma_{j_\ell} \odot S^q\gamma_i = S^\lambda \gamma_{j_\lambda} \odot S^q\gamma_i\odot (\ldots) = S^\lambda (\gamma_{j_\lambda}\odot S^{q-\lambda} \gamma_i)\odot(\ldots)=\zero
                \end{equation*}
                by \cref{lem:two_sided_S_q_minus_ell_gamma_i_odot_gamma_j}. Again it follows that $P=\zero$.
            \end{enumerate}
        \end{enumerate}
    \end{enumerate}
    We see that in all cases $P=\zero$. Therefore, the only terms that survive are those in \cref{eq:two_sided_surviving_terms_c^(q)g_odot_S_q_f} and the proof is completed.
\end{proof}

We can now prove the first implication of \cref{thm:qc_quasi_conserved_iff_polynomial_composition_property} for two-sided landscapes. The proof is identical to \cref{lem:right_sided_quasi_conserved_implies_pcp}.

\begin{lemma}\label{lem:two_sided_quasi_conserved_implies_pcp}
    Let $j\geq 0$ and $\alpha_0,\ldots,\alpha_s\in\F_2$ with $\alpha_0=1$. Then
    \begin{equation*}
        \gamma_j\circ \left(\gamma_0+\sum_{i=1}^s \alpha_i\gamma_i\right) = \sum_{i=0}^s \alpha_i \gamma_{i+j}.
    \end{equation*}
\end{lemma}

\begin{proof}
    We proceed by induction on $j$. The base case $j=0$ is clear, as $\gamma_0=\id$. Note that since $c$ is normed we are able to use the recurrence $\gamma_{j+1}=\cq \odot S^q \gamma_j$ starting from $j=0$. Write $g=\gamma_0+\sum_{i=1}^s \alpha_i\gamma_i$. Then
    \begin{align*}
        \gamma_{j+1} \circ \left(\gamma_0+\sum_{i=1}^s \alpha_i\gamma_i\right)
            &= (\cq \odot S^q \gamma_j)\circ \left(\gamma_0+\sum_{i=1}^s \alpha_i\gamma_i\right) \\
            &= \cq(g) \odot S^q \gamma_j\circ \left(\gamma_0+\sum_{i=1}^s \alpha_i\gamma_i\right) \\
            &= \cq(g) \odot S^q \left(\sum_{i=0}^s \alpha_i\gamma_{i+j}\right) \\
            &= \sum_{i=0}^s \alpha_i \gamma_{i+j+1},
    \end{align*}
    where in the third line we used the induction hypothesis and the last equality follows from \cref{lem:two_sided_c^qg_odot_S^q_f}.
\end{proof}

This shows the first implication of \cref{thm:qc_quasi_conserved_iff_polynomial_composition_property}, i.e. that if $(c,q)$ is quasi-conserved, then it satisfies the polynomial composition property. To show the other direction we prove the contrapositive. 

Suppose that $(c,q)$ is not quasi-conserved. The polynomial composition property would imply that
\begin{equation*}
    \gamma_j\circ(\id+c) = \gamma_j+\gamma_{j+1}
\end{equation*}
for all $j\geq 0$, as functions on $\F_2^\Z$. Set $f_j=\gamma_j\circ(\id+c)$ and $g_j=\gamma_j+\gamma_{j+1}$. It suffices to show that there exist $x\in\F_2^\Z$ and $j\geq 0$ such that $f_j(x)\neq g_j(x)$.

Since $(c,q)$ is not quasi-conserved, there exists $\ell\in\supp(c)\setminus\{q\}$ and $x\in\F_2^\Z$ so that $x$ matches $c$ at positions $0$ and $\ell$. Suppose for now that $\ell<0$. We may assume that $c$ is normed, i.e. $q>0$. Consider the diagram
\begin{equation*}
    \term[2pt]{
        x_{-k+\ell} & \ldots & x_{-k}      & \ldots & x_{-k+q}      & \ldots & x_\ell     & \ldots & x_0       & \ldots & x_q        & \ldots & x_{m+\ell}   & \ldots & x_m     & \ldots & x_{m+q} \\
                    &        & a_{-k}      & \ldots & a_{-k+q}      & \ldots & a_\ell     & \ldots & a_0       & \ldots & a_q        & \ldots & a_{m+\ell}   & \ldots & a_m     & \\
        a_{-k}      & \ldots & a_{-k-\ell} & \ldots & a_{-k+q-\ell} & \ldots & a_0        & \ldots & a_{-\ell} & \ldots & a_{q-\ell} & \ldots & a_m \\
                    &        &             &        & a_{-k}        & \ldots & a_{\ell-q} & \ldots & a_{-q}    & \ldots & a_0        & \ldots & a_{m+\ell-q} & \ldots & a_{m-q} & \ldots & a_m,
    }
\end{equation*}
which describes the entries of $x$ around $0$. We know that it is possible to choose $x$ so that it matches $c$ at positions $0$ and $\ell$ because there is no conflict between the second and third row. Since this selection does not involve the entry $x_{m+q}$, we can then afterward choose $x_{m+q}$ in such a way that $x$ does not match $c$ at position $q$. 

We now compute 
\begin{equation*}
    f_1(x)=c\circ(\id+c)(x)
\end{equation*}
and 
\begin{equation*}
    g_1(x)=c(x)+\gamma_2(x)=c(x)+(\cq\odot S^qc)(x). 
\end{equation*}

Set $y=(\id+c)(x)$ and $z=c(y)$. Since $x$ matches $c$ at position $\ell$, we know that $y_\ell=x_\ell+1$. In particular, since $y_\ell\neq x_\ell=a_\ell$, this implies that $y$ does not match $c$ at position $0$. Therefore, $z_0=0$. 

Now, set $w=c(x)$ and $u=\gamma_2(x)=(\cq\odot S^qc)(x)$. Since $x$ matches $c$ at position $0$, $w_0=1$. Now consider $u$. From the definition it holds that $u_0=1$ if and only if $x$ matches $\cq\odot S^q c$ at position $0$. This is equivalent to $x$ matching $\cq$ at position $0$ and $c$ at position $q$. By construction $x$ does not match $c$ at position $q$. Hence, $u_0=0$. It follows that 
\begin{equation*}
    f_1(x)_0 = z_0 = 0 \neq 1 = w_0+u_0 = (c+\gamma_2)(x)_0=g_1(x)_0.
\end{equation*}
Then $f_1\neq g_1$ and hence $c$ does not satisfy the polynomial composition property.

We now consider the case $\ell>0$. There are two subcases, either $0<\ell<q$ or $0<q<\ell$. Let first $0<\ell<q$. Like before, we consider a diagram depicting of the entries of $x$ around $0$:
\begin{equation*}
    \term[2pt]{
        x_{-k} & \ldots & x_{-k+\ell} & \ldots & x_{-k+q}      & \ldots & x_0       & \ldots & x_\ell     & \ldots & x_q        & \ldots & x_m        & \ldots & x_{m+\ell}  & \ldots & x_{m+q} \\
        a_{-k} & \ldots & a_{-k+\ell} & \ldots & a_{-k+q}      & \ldots & a_0       & \ldots & a_\ell     & \ldots & a_q        & \ldots & a_m  \\
               &        & a_{-k}      & \ldots & a_{-k+q-\ell} & \ldots & a_{-\ell} & \ldots & a_0        & \ldots & a_{q-\ell} & \ldots & a_{m-\ell} & \ldots & a_m   \\
               &        &             &        & a_{-k}        & \ldots & a_{-q}    & \ldots & a_{\ell-q} & \ldots & a_0        & \ldots & a_{m-q}    & \ldots & a_{m+\ell-q}& \ldots & a_m
    }
\end{equation*}
We know that it is possible to choose $x$ so that it matches $c$ at positions $0$ and $\ell$. Furthermore, we can read off that this choice is independent of the entry $x_{m+q}$. From here we can repeat the same argument as before and obtain that $f_1(x)\neq g_1(x)$, implying that $c$ does not satisfy the polynomial composition property.

Finally, consider the case $0<q<\ell$. Here, the entries of $x$ around $0$ look as follows:
\begin{equation*}
    \term[2pt]{
        x_{-k} & \ldots & x_{-k+q} & \ldots & x_{-k+\ell}   & \ldots & x_0       & \ldots & x_\ell     & \ldots & x_q        & \ldots & x_m        & \ldots & x_{m+q}      & \ldots & x_{m+\ell} \\
        a_{-k} & \ldots & a_{-k+q} & \ldots & a_{-k+\ell}   & \ldots & a_0       & \ldots & a_\ell     & \ldots & a_q        & \ldots & a_m  \\
               &        &          &        & a_{-k}        & \ldots & a_{-\ell} & \ldots & a_0        & \ldots & a_{q-\ell} & \ldots & a_{m-\ell} & \ldots & a_{m+q-\ell} & \ldots & a_m  \\
               &        & a_{-k}   & \ldots & a_{-k+\ell-q} & \ldots & a_{-q}    & \ldots & a_{\ell-q} & \ldots & a_0        & \ldots & a_{m-q}    & \ldots & a_m 
    }
\end{equation*}
We see that the previous strategy does not work here. By choosing some $x$ which matches $c$ at positions $0$ and $\ell$ it might automatically also match $c$ at position $q$. 

Instead, we introduce the concept of $j$-nice vectors. We say that $x\in\F_2^\Z$ is $j$-nice if 
\begin{itemize}
    \item $x$ matches $c$ at position $\ell$.
    \item $x$ matches $\gamma_j$ at position $0$. 
\end{itemize}
Observe that if $x$ matches $\gamma_{j+1}=\cq\odot S^q \gamma_j$ at $0$, then it matches $\gamma_j$ at position $q$. With this terminology, the vectors we have previously chosen are $1$-nice but not $2$-nice. Observe that $1$-nice vectors always exist because $(c,q)$ is not quasi-conserved. Furthermore, for large enough $j$ there exist vectors which are not $j$-nice. Indeed, $m+\ell < m+jq$ for large enough $j$ and $m+jq\in\supp(\gamma_j)$ by \cref{lem:two_sided_gammas_formulas}, so that by choosing $x_{m+jq}$ appropriately we find some $x$ which matches $c$ at $\ell$ but not $\gamma_j$ at $0$.

Choose now $j>0$ for which there exists $x\in\F_2^\Z$ which is $j$-nice but not $(j+1)$-nice. By the previous argument such an index $j$ exists. We want to compute $\gamma_j\circ(\id+c)(x)$ and $\gamma_j(x)+\gamma_{j+1}(x)$. 

Set first $y=(\id+c)(x)$ and $z=\gamma_j(y)$. Since $x$ is $j$-nice, it matches $c$ at position $\ell$, so $y_\ell=x_\ell+1$. This implies that $y$ does not match $\cq$ at position $0$ since $y_\ell=x_\ell+1=a_\ell+1$.  In particular, $y$ does not match $\gamma_j$ at $0$. Therefore, $z_0=0$.

Now set $w=\gamma_j(x)$ and $u=\gamma_{j+1}(x)$. Since $x$ is $j$-nice but not $(j+1)$-nice, it matches $\gamma_j$ at position $0$ but not $\gamma_{j+1}$. Hence, $w_0=1$ but $u_0=0$. It follows that
\begin{equation*}
    f_j(x)_0 = z_0 = 0 \neq 1 = w_0+u_0 = g_j(x)_0.
\end{equation*}
This shows that $c$ does not satisfy the polynomial composition property.

\begin{lemma}\label{lem:two_sided_pcp_implies_quasi_conserved}
    Suppose that $c$ satisfies the polynomial composition property. Then $c$ is quasi-conserved.
\end{lemma}

This proves our main theorem for two-sided landscapes.

\begin{theorem}
    Let $c$ be a two-sided landscape. Then \cref{thm:qc_quasi_conserved_iff_polynomial_composition_property} is true for $c$.
\end{theorem}

\begin{proof}
    Combine \cref{lem:two_sided_quasi_conserved_implies_pcp} and \cref{lem:two_sided_pcp_implies_quasi_conserved}.
\end{proof}

\begin{remark}
    The polynomial composition property says that
    \begin{equation*}
        \gamma_j\circ \left(\gamma_0+\sum_{i=1}^s \alpha_i \gamma_i\right) = \sum_{i=0}^s \alpha_i \gamma_{i+j}
    \end{equation*}
    for all $j\geq 0$ and all $\alpha_1,\ldots,\alpha_s\in\F_2$. As a particular instance,
    \begin{equation}\label{eq:two_sided_weak_pcp}
        \gamma_j\circ(\id+c) = \gamma_j+\gamma_{j+1}
    \end{equation}
    for all $j\geq 0$ and
    \begin{equation}\label{eq:two_sided_very_weak_pcp}
        c\circ(\id+c) = c+\gamma_2.
    \end{equation}
    Our proof shows that in fact \cref{eq:two_sided_weak_pcp} already implies that $c$ is quasi-conserved. Furthermore, the property \cref{eq:two_sided_very_weak_pcp} is already enough to show that $S^\ell c\odot c=\zero$ for all $\ell\in\supp(c)\setminus\{q\}$ with $\ell<0$ or $0<\ell<q$. This in particular applies to all landscapes with $q=m$. 

    Numerical evidence suggests that \cref{eq:two_sided_very_weak_pcp} is enough to show that $c$ is quasi-conserved for all landscapes $c$. We leave this question open. 
\end{remark}

\begin{problem}
    Show that the following two statements are equivalent:
    \begin{enumerate}
        \item $c$ is quasi-conserved with special index $q$.
        \item $c\circ(\id+c)=c+\cq\odot S^q c$.
    \end{enumerate}
\end{problem}

\begin{remark}[added after seeing \cite{lyu2026classshiftinvariant}]
    \cref{thm:qc_quasi_conserved_iff_polynomial_composition_property} proves that quasi-conservedness and the polynomial composition property are equivalent over $\F_2^\Z$. By restricting the functions $\gamma_j$ to $\mP_n\cong \F_2^n$, this also shows that if $c$ is quasi-conserved over $\F_2^\Z$, then $(\{\gamma_j: j\geq 0\})$ satisfy the polynomial composition property over $\F_2^n$. 

    More generally, one can adapt \cref{def:qc_conserved_and_quasi_conserved_landscapes} to $\F_2^n$: A landscape $c$ is $n$-quasi-conserved if $c$ is not conserved and there exists one index $q\in\supp(c)$ such that whenever $x\in\F_2^n$ matches $c$ at position $i\in\{0,\ldots,n-1\}$, then $x$ does not match $c$ at all positions $i+\ell$ for $\ell\in\supp(c)\setminus\{q\}$. It is possible that a landscape $c$ is $n$-quasi-conserved for some $n$ while not being quasi-conserved over $\F_2^\Z$. For example, $0*\mt01$ is $6$-quasi-conserved with special index $2$, but not quasi-conserved. Furthermore, for small enough $n$ most landscapes are trivially $n$-quasi-conserved, since there are no vectors $x\in\F_2^n$ that match the landscape at any position. 

    In \cite{lyu2026classshiftinvariant} it is shown directly that $n$-quasi-conservedness and the polynomial composition property over $\F_2^n$ are equivalent for sufficiently large $n$. More precisely, they give the sharp lower bound $n>\max\{2m+k, 2k+m\}$. 

    Our proof of the fact that quasi-conservedness implies the polynomial composition property only uses that $c$ is quasi-conserved in the form of \cref{lem:qc_equivalent_characterization_conserved_quasi_conserved_with_S_ell_c_odot_c}. It is easy to adapt the statement of \cref{lem:qc_equivalent_characterization_conserved_quasi_conserved_with_S_ell_c_odot_c} to the case of $n$-quasi-conservedness, and hence our proof also shows this implication directly over $\F_2^n$ for all $n$. 
    
    For the converse implication, that the polynomial composition property implies quasi-conservedness, we use the concept of $j$-nice vectors. Since $j$-nice vectors also exist over $\F_2^n$ for sufficiently large $n$, the same argument proves this implication for sufficiently large $n$. A closer examination of the argument yields the bound $n>\max\{2m+k, 2k+m+q\}$, where $q>0$ since $c$ is normed. More generally, the argument applies for every $n$ for which some $j$-nice vector exists. 
\end{remark}

\section{Functions that behave like polynomials}\label{sec:pcp_studying_pcp}

Let $V\in\{\F_2^n, \F_2^\Z\}$ and let $F(V)$ be the vector space of all functions $V\to V$. Suppose that $T:F(V)\to F(V)$ is a linear operator. Take $\gamma_0:V\to V$ arbitrary and set $\gamma_{j+1} = T \gamma_j$ for $j\geq 0$. The subspace
\begin{equation*}
    \Gamma = \spn\{\gamma_0, \gamma_1, \gamma_2, \ldots\} \leq F(V)
\end{equation*}
is then called the cyclic subspace of $\gamma_0$ under $T$. The elements in $\Gamma$ have the form $p(T)(\gamma_0)$ for some polynomial $p\in \F_2[X]$. In other words, the evaluation map $\ev_T:\F_2[X]\to\Gamma$ given by $p(X)\mapsto p(T)\gamma_0$ is surjective. By the First Isomorphism Theorem,
\begin{equation*}
    \F_2[X]/\ker(\ev_T) \cong \Gamma.
\end{equation*}
One can quickly check that $\ker(\ev_T)$ is an ideal in $\F_2[X]$. As $\F_2[X]$ is a principal ideal ring, there exists a polynomial $m_T(X)\in\F_2[X]$ such that $\ker(\ev_T) = (m_T(X))$ is the ideal generated by $m_T(X)$. We call $m_T(X)$ the minimal polynomial of $T$. 

Denote by $\mR$ the quotient ring $\F_2[X]/(m_T(X))$. We write $[f]$ for the coset $f+(m_T(X))$ in $\mR$. The isomorphism $\psi:\mR\to\Gamma$ has the form
\begin{equation*}
    \psi\left(\left[\sum_{i=0}^s \alpha_i X^i\right]\right) = \sum_{i=0}^s \alpha_i \gamma_i
\end{equation*}
and the inverse $\varphi:\Gamma\to \mR$ is given by
\begin{equation*}
    \varphi\left(\sum_{i=0}^s \alpha_i \gamma_i\right) = \left[\sum_{i=0}^s \alpha_i X^i\right].
\end{equation*}
Define the coset
\begin{equation*}
    \mG = \gamma_0 + \spn\{\gamma_1, \gamma_2, \ldots\} \subseteq \Gamma.
\end{equation*}
Recall the polynomial composition property (\cref{def:qc_polynomial_composition_property}): the family $\{\gamma_j: j\geq 0\}$ satisfies the polynomial composition property if for every $j\geq 0$ and $\alpha_0,\ldots,\alpha_s\in\F_2$ with $\alpha_0=1$ it holds that
\begin{equation*}
    \gamma_j\circ\left(\sum_{i=0}^s \alpha_i \gamma_i\right) = \sum_{i=0}^s \alpha_i \gamma_{i+j}.
\end{equation*}
We then have the following lemma.

\begin{lemma}\label{lem:pcp_G_is_monoid}
    Suppose that $\gamma_0=\id$ is the identity function and the family $\{\gamma_j: j\geq 0\}$ satisfies the polynomial composition property. Then $\mG= \gamma_0 + \spn\{\gamma_j: j\geq 1\} $ is a monoid with respect to composition.
\end{lemma}

\begin{proof}
    Let $f,g\in \mG$ with $f=\sum_{j=0}^t \beta_j\gamma_j$ and $g=\sum_{i=0}^s \alpha_i\gamma_i$, where $\alpha_0=\beta_0=1$. Using the polynomial composition property, 
    \begin{equation*}
        \gamma_j \circ g = \gamma_j \circ \left(\sum_{i=0}^s \alpha_i\gamma_i\right) = \sum_{i=0}^s \alpha_i\gamma_{i+j} \in \spn\{\gamma_1, \gamma_2, \ldots\}.
    \end{equation*} 
    Hence, 
    \begin{equation*}
        f\circ g = \gamma_0\circ g + \sum_{j=1}^t \beta_j \gamma_j\circ g \in g+\spn\{\gamma_1,\gamma_2,\ldots\} = \mG.
    \end{equation*}
    Finally, note that $\gamma_0=\id\in \mG$, which is the neutral element with respect to composition.
\end{proof}

\begin{lemma}\label{lem:pcp_restriction_of_phi_is_homomorphism}
    Let $f\in \Gamma, g\in \mG$. Then $\varphi(f\circ g) = \varphi(f)\cdot\varphi(g)$. In particular, the restriction of $\varphi$ to $\mG$ is a monoid homomorphism. 
\end{lemma}

\begin{proof}
    Let first $f=\gamma_j$ and $g=\gamma_0+\sum_{i=1}^s \alpha_i \gamma_i\in \mG$. By the polynomial composition property
    \begin{equation*}
        \varphi(\gamma_j\circ g) = \varphi\left(\sum_{i=0}^s \alpha_i \gamma_{i+j}\right) = \left[\sum_{i=0}^s \alpha_i X^{i+j}\right] = \left[X^j\right]\cdot\left[\sum_{i=0}^s \alpha_i X^i\right]
            = \varphi(\gamma_j)\cdot \varphi(g).
    \end{equation*}
    The general case follows by linearity of $\varphi$. Set $f=\sum_{j=0}^t \beta_j\gamma_j\in\Gamma$. Then
    \begin{align*}
        \varphi(f\circ g) &= \varphi\left(\sum_{j=0}^t \beta_j\gamma_j\circ g\right) = \sum_{j=0}^t \beta_j \varphi(\gamma_j\circ g) \\
            &= \sum_{j=0}^t \beta_j\varphi(\gamma_j)\cdot\varphi(g) = \varphi\left(\sum_{j=0}^t \beta_j \gamma_j\right)\cdot\varphi(g) 
             = \varphi(f)\cdot\varphi(g).
    \end{align*}
    Note that $\varphi(\gamma_0)=[1]\in \mR$, so that $\varphi$ restricted to $\mG$ is a monoid homomorphism. 
\end{proof}

It is easy to see that the set $1+\spn\{X,X^2,\ldots\}\subseteq \F_2[X]$ is a monoid with respect to polynomial multiplication. Denote by $\pi$ the canonical ring homomorphism $\F_2[X]\to \mR, f\mapsto [f]$. As the image of a monoid under a ring homomorphism is again a monoid, it follows that the image 
\begin{equation*}
    \mM:=\pi(1+\spn\{X,X^2,\ldots\})\subseteq \mR
\end{equation*}
is again a monoid.

\begin{theorem}\label{thm:pcp_G_isomorphic_M_monoid}
    $\mG$ is isomorphic to $\mM$ as a monoid. In particular, $\mG$ is Abelian.
\end{theorem}

\begin{proof}
    By \cref{lem:pcp_restriction_of_phi_is_homomorphism} the map $\varphi:\mG\to \mM$ is a monoid homomorphism. It is easy to see that $\varphi:\mG\to \mM$ has the inverse $\psi:\mM\to \mG$. The claim follows.
\end{proof}

\cref{thm:pcp_G_isomorphic_M_monoid} implies that we can study the monoid $\mG$ by studying the monoid
\begin{equation*}
    \mM = \left\{\left[1+\sum_{i=1}^s \alpha_i X^i\right]:\alpha_i\in\F_2, s\in\N\right\}\subseteq \mR
\end{equation*}
consisting of (cosets of) polynomials. As polynomials are very well-studied, this allows us to gain insight into the structure of $\mG$. Of particular interest are the units in the monoid $\mG$, i.e. the set
\begin{equation*}
    \mG^\times = \{g\in \mG: \exists h\in \mG: hg=gh=\id\}.
\end{equation*}
Note that $\mG^\times$ is a group. 

\begin{corollary}\label{cor:pcp_G*_M*_iso_as_groups}
    $\mG^\times$ is isomorphic to $\mM^\times$ as a group.
\end{corollary}

\begin{remark}\label{rem:pcp_unit_iff_permutation_if_G_finite}
    Let $X$ be a set and $(H,\circ)$ be a monoid consisting of functions $X\to X$. Let $h\in H$ be a unit, i.e. there exists $h^{-1}\in H$ with $h\circ h^{-1}=h^{-1}\circ h=\id$. It then follows that $h$ is a permutation on $X$. However, if $h:X\to X$ is a permutation with inverse $h^{-1}$ then it does not follow that $h$ is a unit in $H$ because $h^{-1}$ is not necessarily contained in $H$.

    As an example, let $X=\Z$ and $H=\N_0$, where $h\in\N_0$ acts on $\Z$ by $z\mapsto z+h$. Then it is clear that $h=3:z\mapsto z+3$ is a permutation on $\Z$, but its inverse map $z\mapsto z-3$ is not in $H$ and $h$ is not a unit in $H$.

    On the other hand, if a permutation $h\in H$ has finite order, i.e. $h^s=\id$ for some $s\in\N$, then $h^{s-1}\in H$ is the inverse of $h$ and $h\in H$ is a unit. As a special case it follows that if $H$ is finite then $h\in H$ is a unit if and only if it is a permutation on $X$. In particular, if $m_T\neq 0$, then $\mG$ is finite, so that in the following we can use that $f\in\mG$ is a unit if and only if $f$ is a permutation on $V$.
\end{remark}

\subsection{Minimal polynomials of quasi-conserved landscapes}\label{subsec:pcp_computing_min_polys}

Let now $c$ be a normed nondegenerate quasi-conserved landscape. Let $V\in \{\F_2^n, \F_2^\Z\}$ and let $T:F(V)\to F(V)$ be the linear operator given by $Tf=\cq\odot S^q f$. Set $\gamma_0=\id$, then
\begin{equation*}
    T\gamma_j = \cq \odot S^q \gamma_j=\gamma_{j+1}
\end{equation*}
for $j\geq 0$. By the discussion at the beginning of \cref{sec:pcp_studying_pcp} the operator $T$ has a minimal polynomial $m_T(X)$. We then define the minimal polynomial of the landscape $c$.

\begin{definition}\label{def:pcp_min_poly_of_landscape}
    With the notation from above, we define the minimal polynomial $m_{c,V}(X)$ of $c$ as the minimal polynomial of $T$ restricted to $\Gamma$, i.e. 
    \begin{equation*}
        m_{c,V}(X)=m_T(X).
    \end{equation*} 
     If $V=\F_2^n$ we also denote it by $m_{c,n}(X)$ and if $V=\F_2^\Z$ we also write $m_{c,\Z}(X)$.
\end{definition}

Since $c$ is quasi-conserved, \cref{thm:qc_quasi_conserved_iff_polynomial_composition_property} implies that the family $\{\gamma_j: j\geq 0\}$ satisfies the polynomial composition property. Using the results of the previous section we then know that
\begin{equation*}
    \mG = \gamma_0 + \spn\{\gamma_j: j\geq 1\}
\end{equation*}
is an abelian monoid isomorphic to
\begin{equation*}
    \mM = \left\{\left[1+\sum_{i=1}^t \alpha_i X^i\right]\right\}\subseteq \F_2[X]/(m_{c,V}(X))
\end{equation*}
and its structure is completely determined by $m_{c,V}(X)$. Our goal from here is to determine all possible minimal polynomials $m_{c,V}(X)$. We start with an example. 

\begin{example}\label{ex:pcp_min_polys_of_chi}
    Let $c=*01$ and $V=\F_2^n$. In \cite{kriepke2024algebraic} it was shown that if $n$ is odd then $\gamma_0,\gamma_1,\ldots,\gamma_{\frac{n-1}{2}}$ are linearly independent and $\gamma_{\frac{n+1}{2}}=\zero$. For $n$ even in \cite{kriepke2026siblings} it was shown that $\gamma_0,\ldots,\gamma_{n-1}$ are linearly independent and $\gamma_n=\gamma_{n/2}$. 
    Now consider $V=\F_2^\Z$. By \cref{lem:right_sided_gammas_never_zero} we know that $\gamma_j\neq \zero$ for all $j\geq 0$. By \cref{lem:two_sided_support_gamma_i}, $\max\supp(\gamma_j)=2j$. In particular, the supports $\supp(\gamma_j)$ for $j\geq 0$ are distinct and hence the functions $\{\gamma_j: j\geq 0\}$ are linearly independent by \cref{lem:pcp_linear_dependence_landscapes}, which we prove below. In conclusion,
    \begin{equation*}
        m_{*01,V}(X) = 
            \begin{cases}
                X^{(n+1)/2},  & V=\F_2^n, n \text{ is odd,} \\
                X^n+ X^{n/2}, & V=\F_2^n, n \text{ is even,} \\
                0,            & V=\F_2^\Z.
            \end{cases}
    \end{equation*}
\end{example}

\begin{lemma}\label{lem:pcp_linear_dependence_landscapes}
    Let $c^{(1)},\ldots, c^{(s)}$ be landscapes such that the functions $c^{(1)},\ldots, c^{(s)}$ are linearly dependent and nonzero. Then there exist two integers $1\leq i,j\leq s$ such that $\supp(c^{(i)})=\supp(c^{(j)})$.
\end{lemma}

\begin{proof}
    Write $c^{(i)} = a_{-k_i}^{(i)}\ldots a_{m_i}^{(i)}$ for all $i=1,\ldots, s$ and let $c_0^{(i)}$ be the $0$-th coordinate function of $c^{(i)}$. Observe that if $c^{(1)},\ldots, c^{(s)}$ are linearly dependent, then also $c_0^{(1)},\ldots, c_0^{(s)}$ are linearly dependent. The function $c_0^{(i)}(x)$ is given in \cref{eq:prelim_coordinate_function_from_landscape}, i.e. 
    \begin{equation*}
        c_0^{(i)}(x) = \varepsilon_{-k_i}\cdot \varepsilon_{-k_i+1}\cdots \varepsilon_{m_i}
    \end{equation*}
    where
    \begin{equation*}
        \varepsilon_\ell = \begin{cases}
            x_\ell,         & \text{if } a_\ell^{(i)} = 1, \\
            1+x_\ell,       & \text{if } a_\ell^{(i)} = 0, \\
            1,              & \text{if } a_\ell^{(i)} = \mt.
        \end{cases}
    \end{equation*}
    The leading term of the polynomial $c_0^{(i)}$ is $\prod_{\ell\in\supp(c^{(i)})} x_\ell$. As nonzero polynomials with distinct leading terms are linearly independent, there exist at least two indices $i,j$ such that $c_0^{(i)}, c_0^{(j)}$ have the same leading term. This implies $\supp(c^{(i)})=\supp(c^{(j)})$. 
\end{proof}

In fact, every minimal polynomial has one of the three forms shown in \cref{ex:pcp_min_polys_of_chi}. We introduce terminology for these three forms. 

\begin{definition}
    \begin{enumerate}
        \item[(i)] If $m_{c,V}(X)=X^j$ for some integer $j>0$, then we say that $(c,V)$ is $\gamma$-nilpotent. The number $j$ is called the index of $\gamma$-nilpotency.
        \item[(ii)] If $m_{c,V}(X)=X^j+X^k=X^j(X^{k-j}+1)$ for some integers $j,k>0$ with $j<k$, then we say that $(c,V)$ is $\gamma$-periodic. The number $j$ is called the preperiod and $k-j$ is called the period.
        \item[(iii)] If $m_{c,V}(X)=0$, then we say that $(c,V)$ is $\gamma$-independent. 
    \end{enumerate}
    Instead of $(c,V)$ we also say $(c,n)$ or $(c,\Z)$, depending on $V=\F_2^n$ or $V=\F_2^\Z$, respectively.
\end{definition}

This terminology is motivated by the following observations. First, if $\gamma_j=\zero$ for some $j\geq 0$, then $T^j$ is the zero operator on $\Gamma$, i.e. $T$ is nilpotent. Next, if $\gamma_j=\gamma_k$ for some $j,k$, then the sequence of functions $\gamma_0,\gamma_1,\gamma_2,\ldots$ is eventually periodic. Finally, if $m_{c,V}(X)=0$, then the sequence $\gamma_0,\gamma_1,\gamma_2,\ldots$ is linearly independent.

The case $V=\F_2^\Z$ admits a simple characterization. In fact, the proof is almost identical to the argument from \cref{ex:pcp_min_polys_of_chi}. 

\begin{lemma}\label{lem:pcp_min_poly_F_2_Z}
    Let $V=\F_2^\Z$. Then $(c,V)$ is $\gamma$-nilpotent or $\gamma$-independent. 
\end{lemma}

\begin{proof}
    If $\gamma_j\neq \zero$ for all $j\geq 0$, then by \cref{lem:two_sided_support_gamma_i} we have $\max\supp(\gamma_j)=(j-1)q+m$. In particular the supports $\supp(\gamma_j)$ for $j\geq 0$ are distinct and hence the functions $\{\gamma_j: j\geq 0\}$ are linearly independent by \cref{lem:pcp_linear_dependence_landscapes}, giving $m_{c,\Z}(X)=0$.

    Otherwise, $\gamma_j=\zero$ for some $j>0$, i.e. $m_{c,\Z}(X)\mid X^j$. 
\end{proof}

In fact, we can compute the index of $\gamma$-nilpotency. For this we need the following auxiliary lemma.

\begin{lemma}\label{lem:pcp_product_c1_ldots_cs_equals_zero}
    Let $c^{(1)},\ldots, c^{(s)}$ be landscapes such that $c^{(1)}\odot\cdots \odot c^{(s)}=\zero$. Then there exist two indices $i,j$ such that $c^{(i)}\odot c^{(j)}=\zero$. 
\end{lemma}

\begin{proof}
    Write $c^{(i)}=a_{-k_i}^{(i)}\ldots a_{m_i}^{(i)}$ for all $1\leq i\leq s$ and suppose without loss of generality (by setting $a_\ell^{(i)}=\mt$ if necessary) that there exist some $k,m\geq 0$ with $k_i =k$ and $m_i=m$ for all $i$. We know
    \begin{equation}\label{eq:pcp_product_c1_ldots_cs_equals_zero}
        c^{(1)}\odot\cdots\odot c^{(s)}
            = \term[2pt]{
                a_{-k}^{(1)} & \ldots & a_{m}^{(1)} \\
                \vdots       & \ddots & \vdots        \\
                a_{-k}^{(s)} & \ldots & a_{m}^{(s)}
            }
            =\zero,
    \end{equation}
    which means that there is a conflict in \cref{eq:pcp_product_c1_ldots_cs_equals_zero}, i.e. there exist rows $i,j$ and a column $\ell$ with $a_\ell^{(i)}=1$ and $a_\ell^{(j)}=0$. Then $c^{(i)}\odot c^{(j)}$ has the same conflict. 
\end{proof}

We note that over $V=\F_2^n$ it may already be the case that $c^{(i)}=\zero$ for some index $i$, this is not a problem.

\begin{lemma}\label{lem:pcp_index_of_nilpotency}
    Let $h\geq 0$. The following statements are equivalent:
    \begin{enumerate}
        \item The number $h$ is minimal with the property that $S^{hq} c\odot \cq=\zero$.
        \item The minimal polynomial of $c$ is $m_{c,V}(X)=X^{h+1}$. 
    \end{enumerate}
\end{lemma}

\begin{proof}
    We first show the following two statements. Afterward, the claim follows from the minimality assumption on $h$ and the minimal polynomial.
    \begin{enumerate}
        \item[(A)] Let $j\geq 0$ with $S^{jq}c\odot \cq=\zero$. Then $m_{c,V}(X)\mid X^{j+1}$.
        \item[(B)] Let $j\geq 0$ with $m_{c,V}(X) \mid X^{j+1}$. Then there exists $\ell\leq j$ with $S^{\ell q}c\odot \cq =\zero$. 
    \end{enumerate}
    $(A)$: The product
    \begin{equation*}
        S^{jq}c\odot \cq = \zero
    \end{equation*}
    appears in the formula of $\gamma_{j+1}$, see \cref{eq:two_sided_formula_gamma_j}. This implies $\gamma_{j+1}=\zero$, i.e. $m_{c,V}(X)\mid X^{j+1}$.

    $(B)$: Assume that $\gamma_{j+1}=\zero$ for some $j\geq 0$, i.e.
    \begin{equation*}
        \zero = \gamma_{j+1}=\cq\odot\cdots\odot S^{\mu q}\cq \odot\cdots\odot S^{\rho q}\cq \odot\cdots\odot S^{jq}c.
    \end{equation*}
    By \cref{lem:pcp_product_c1_ldots_cs_equals_zero} there exist $\mu,\rho$ with $0\leq \mu < \rho\leq j$ and $\tilde c\in\{c,\cq\}$ such that
    \begin{equation*}
        S^{\mu q}\cq \odot S^{\rho q}\tilde c=\zero.
    \end{equation*}
    Observe that if $S^{\mu q}\cq \odot S^{\rho q}\cq=\zero$, then also $S^{\mu q}\cq \odot S^{\rho q}c=\zero$. Therefore,
    \begin{equation*}
        \zero = S^{\mu q}\cq\odot S^{\rho q}c = S^{\mu q}(\cq \odot S^{(\rho-\mu)q}c)
    \end{equation*}
    and hence $S^{\ell q}c\odot \cq=\zero$ with $\ell=\rho-\mu\leq j$.

    We can now show the original claims:

    $(a)\implies (b)$: Suppose $h\geq 0$ is minimal such that $S^{hq}c\odot \cq=\zero$. Assume $m_{c,V}(X)\neq X^{h+1}$. By claim $(A)$ we know $m_{c,V}(X)\mid X^{h+1}$, so there exists $0\leq j<h$ with $m_{c,V}(X)=X^{j+1}$. But then by $(B)$ there exists $\ell\leq j<h$ with $S^{\ell q}c\odot \cq=\zero$, contradicting the minimality of $h$.

    $(b)\implies (a)$: Suppose $m_{c,V}(X) = X^{h+1}$. By claim $(B)$ there exists $\ell\leq h$ with $S^{\ell q}c\odot \cq=\zero$. Assume $\ell < h$. Then by $(A)$ we have $m_{c,V}(X)\mid X^{\ell+1}$, a contradiction.
\end{proof}

\cref{lem:pcp_index_of_nilpotency} gives an effective procedure to compute the index of $\gamma$-nilpotency. Indeed, first compute the set
\begin{equation*}
    \mL_{\cq} = \{\ell\in\Z: S^\ell c\odot \cq=\zero\},
\end{equation*}
which is a finite computation, since the supports of $S^\ell c$ and $\cq$ are disjoint for $\ell$ large enough. If $V=\F_2^n$, then we may view $\mL_{\cq}$ as a subset of $\Z_n=\Z/n\Z$, the integers modulo $n$. Since $S^n=\id$ on $\F_2^n$, this is well-defined. Afterward, determine the minimal $h>0$ such that $hq\in\mL_{\cq}$, if it exists.

From now on we focus on the case $V=\F_2^n$ and start with the following simple observation. We define
\begin{equation*}
    Q = \{hq: h\in\N_0\} \leq \Z_n,
\end{equation*}
which is the subgroup generated by $q$ in $\Z_n$. 

\begin{lemma}\label{lem:pcp_S_ell_c_odot_cq_zero_iff_S_ell_c_odot_c_zero}
    The following are equivalent:
    \begin{enumerate}
        \item There exists $\ell\in Q$ with $S^\ell c\odot \cq=\zero$.
        \item There exists $\ell'\in Q$ with $S^{\ell'} c\odot c=\zero$.
    \end{enumerate}
\end{lemma}

\begin{proof}
    It is clear that $(a)$ implies $(b)$, so it suffices to show that $(b)$ implies $(a)$. 

    Take $\ell\in Q$ with $S^\ell c\odot c=\zero$. By \cref{lem:two_sided_S^ell_c_or_c*_equals_zero} combined with \cref{rem:two_sided_S^ell_c_or_c*_equals_zero_also_for_ell_notin_supp} at least one of $S^\ell c\odot \cq=\zero$ or $S^\ell \cq \odot c = S^\ell (\cq \odot S^{-\ell}c)=\zero$ holds. As $Q\leq \Z_n$ is a subgroup, both $\ell$ and $-\ell$ are in $Q$. The claim follows. 
\end{proof}

Recall that in \cref{sec:quasi_conserved_landscapes} we introduced the set
\begin{equation*}
    \mL_c = \{\ell\in\Z: S^\ell c\odot c=\zero\}.
\end{equation*}
Note that $S^\ell c\odot c = S^\ell (c\odot S^{-\ell} c)$, showing that $\ell\in\mL_c$ if and only if $-\ell\in\mL_c$. This symmetry allows a nicer treatment for $\mL_c$ compared to the set $\mL_{\cq}$. Again, we may view $\mL_c$ as a subset of $\Z_n$, since $S^n=\id$ on $\F_2^n$. By combining \cref{lem:pcp_index_of_nilpotency,lem:pcp_S_ell_c_odot_cq_zero_iff_S_ell_c_odot_c_zero} we obtain the following corollary.

\begin{corollary}\label{cor:pcp_min_poly_X^j_iff_exists_ell_in_Q_cap_L}
    The following statements are equivalent:
    \begin{enumerate}
        \item The set $Q\cap \mL_c$ is non-empty. 
        \item $(c,n)$ is $\gamma$-nilpotent.
    \end{enumerate}
\end{corollary}

As $Q\leq \Z_n$ is the subgroup generated by $\gcd(q,n)$, this leads to a characterization of all dimensions $n$ for which $(c,n)$ is $\gamma$-nilpotent.

\begin{lemma}\label{lem:pcp_char_of_n_which_are_nilpotent}
    Set
    \begin{equation*}
        \mD_c = \{d: d\mid q, \exists \ell\in\mL_c \text{ s.t. } d\mid \ell\}.
    \end{equation*}
    Then $(c,n)$ is $\gamma$-nilpotent if and only if $\gcd(q,n)\in\mD_c$.
\end{lemma}

\begin{proof}
    By \cref{cor:pcp_min_poly_X^j_iff_exists_ell_in_Q_cap_L} $(c,n)$ is $\gamma$-nilpotent if and only if there exists some $\ell\in \mL_c\cap Q$. Because $Q\leq\Z_n$ is generated by $\gcd(q,n)$, it holds that $\ell\in Q$ if and only if $\gcd(q,n)\mid \ell$ as integers, i.e. $\gcd(q,n)\in\mD_c$.
\end{proof}

\begin{example}\label{ex:pcp_0_01_nilpotent}
    Let $c=*0\mt01$. Then $q=4$ and $\mL_c=\{-3,-1,1,3\}$, so that $\mD_c=\{\pm 1\}$. Therefore, $(c,n)$ is $\gamma$-nilpotent if and only if $\gcd(4,n)=1$, i.e. $n$ is odd.
\end{example}

Observe that $1\in\mD_c$ as long as $\abs{\supp(c)}>1$ (which we have assumed in \cref{ass:qc_assume_not_affine} as otherwise $c$ is affine), since $\supp(c)\setminus\{q\}\subseteq \mL_c$. Then \cref{lem:pcp_char_of_n_which_are_nilpotent} implies that there are always infinitely many $n\in\N$ such that $(c,n)$ is $\gamma$-nilpotent. For the specific case of $q=1$ we get an even stronger statement.

\begin{corollary}\label{cor:pcp_landscape_with_q_equal_1_is_gamma_nilpotent_for_every_n}
    Suppose $q=1$. Then $(c,n)$ is $\gamma$-nilpotent for every $n\in\N$. 
\end{corollary}

\begin{proof}
    For every $n\in\N$ we have $\gcd(q,n)=\gcd(1,n)=1\in\mD_c$. 
\end{proof}

The following lemma determines the relationship between the infinite-dimensional and the finite-dimensional setting.

\begin{lemma}\label{lem:pcp_nilpotent_on_Z_iff_nilpotent_for_every_n}
    $(c,\Z)$ is $\gamma$-nilpotent if and only if $(c,n)$ is $\gamma$-nilpotent for every $n\in\N$. If $(c,\Z)$ is $\gamma$-nilpotent of index $s$, then $(c,n)$ is $\gamma$-nilpotent of index $\leq s$ for every $n\in\N$.
\end{lemma}

\begin{proof}
    Suppose $(c,\Z)$ is $\gamma$-nilpotent of index $s$, i.e. $\gamma_s=\zero$ on $\F_2^\Z$. But then also $\gamma_s=\zero$ as a function on $\F_2^n$, implying that $m_{c,n}(X)\mid X^s$. 

    Conversely, assume that $(c,n)$ is $\gamma$-nilpotent for every $n\in\N$. Our goal is to show that there exists some positive $h$ such that $S^{hq}c\odot\cq=\zero$, then we are done by \cref{lem:pcp_index_of_nilpotency}. By \cref{lem:pcp_char_of_n_which_are_nilpotent} we know that $\gcd(q,n)\in\mD_c$ for every $n\in\N$. Choose $n=q$, then $q\in\mD_c$. By the definition of the set $\mD_c$ this means that there exists some $h\in\Z$ such that $hq\in\mL_c$, i.e. $S^{hq}c\odot c=\zero$. Observe that $h\neq 0$, as $S^0c\odot c=c\odot c=c\neq \zero$. Choose $h$ with the property that $\abs{h}$ is minimal among all $h$ with $hq\in\mL_c$. If there are two such $h$ then choose the positive one.

    Suppose first $h>0$. By combining \cref{lem:two_sided_S^ell_c_or_c*_equals_zero} and \cref{rem:two_sided_S^ell_c_or_c*_equals_zero_also_for_ell_notin_supp} we know that $S^{hq}c\odot\cq=\zero$ or $q-hq=(1-h)q\in\supp(c)\setminus\{q\}\subseteq \mL_c$. Since $h$ is positive we have $\abs{1-h}<\abs{h}$ and by the minimality assumption on $h$ it follows that $q-hq\notin\mL_c$. Hence, $S^{hq}c\odot \cq=\zero$ and we are done. 
    
    Suppose now that $h<0$. Again, by \cref{lem:two_sided_S^ell_c_or_c*_equals_zero} and \cref{rem:two_sided_S^ell_c_or_c*_equals_zero_also_for_ell_notin_supp} imply $S^{hq}\cq \odot c=\zero$ or $q+hq\in\supp(c)\setminus\{q\}\subseteq \mL_c$. Since $h<0$ we have $\abs{1+h}<\abs{h}$, so that $q+hq\notin\mL_c$. Therefore, $S^{-hq}c\odot \cq = S^{-hq}(c\odot S^{hq}\cq)=\zero$ and we are done by choosing $-h$, which is positive.
\end{proof}

We obtain the following corollary.

\begin{corollary}\label{cor:pcp_equivalent_conditions_to_gamma_nilpotent_on_Z_and_every_n}
    The following are equivalent:
    \begin{enumerate}
        \item $(c,\Z)$ is $\gamma$-nilpotent.
        \item $(c,n)$ is $\gamma$-nilpotent for every $n\in\N$. 
        \item The set $\mD_c$ contains every divisor of $q$.
        \item There exists $h\in\Z$ such that $hq\in\mL_c$. 
    \end{enumerate}
\end{corollary}

\begin{proof}
    \cref{lem:pcp_nilpotent_on_Z_iff_nilpotent_for_every_n} shows that $(a)$ and $(b)$ are equivalent. \cref{lem:pcp_char_of_n_which_are_nilpotent} shows that $(b)$ and $(c)$ are equivalent, by choosing $n$ as a divisor of $q$. The equivalence of $(c)$ and $(d)$ follows from the observation that $q\in\mD_c$ if and only if $(d)$ holds. 
\end{proof}

From now on we assume that there exists some $n$ such that $(c,n)$ is not $\gamma$-nilpotent. As $\F_2^n$ is finite dimensional, the minimal polynomial of $c$ has to be nonzero. First, we present a simple example of two linearly dependent functions $\gamma_r$ and $\gamma_{2r}$. 

\begin{lemma}\label{lem:pcp_min_poly_divides_X^2r_plus_X^r}
    Let $r$ be the additive order of $q$ modulo $n$. Then 
    \begin{equation*}
        m_{c,n}(X) \mid X^{2r}+X^r=X^r(X^r+1).
    \end{equation*}
\end{lemma}

\begin{proof}
    Let $f=\cq\odot S^q \cq \odot\cdots\odot S^{(r-1)q}\cq$. Then
    \begin{align*}
        \gamma_{2r} &= \cq \odot S^q\cq \odot\cdots\odot S^{(r-1)q}\cq \odot S^{rq}\cq \odot S^{(r+1)q}\cq\odot\cdots\odot S^{(2r-1)q}(\cq \odot A_q) \\
            & = f\odot S^{rq}f\odot S^{(2r-1)q}A_q = f\odot S^{(r-1)q}A_q = \gamma_r,
    \end{align*}
    i.e. $\gamma_{2r}+\gamma_r=\zero$. The claim follows.
\end{proof}

Observe that the previous lemma also applies if $(c,n)$ is $\gamma$-nilpotent, giving a bound on the index of $\gamma$-nilpotency. 

\begin{corollary}\label{cor:index_of_nilpotency_bounded_by_additive_order}
    Assume that $(c,n)$ is $\gamma$-nilpotent of index $j$ and let $r$ be the additive order of $q$ modulo $n$. Then $j\leq r$.
\end{corollary}

\begin{proof}
    Follows from the fact that $X^j = m_{c,n}(X) \mid X^r(X^r+1)$. 
\end{proof}

The subsequent lemma determines $m_{c,n}(X)$ in the case that $(c,n)$ is not $\gamma$-nilpotent.

\begin{lemma}\label{lem:pcp_periodic_if_not_nilpotent}
    Suppose that $(c,n)$ is not $\gamma$-nilpotent. Let $r$ be the additive order of $q$ modulo $n$. Then $(c,n)$ is $\gamma$-periodic with period $r$ and preperiod $j\leq r$.
\end{lemma}

\begin{proof}
    Suppose that $\gamma_1,\ldots, \gamma_s$ are linearly dependent. Such a number $s$ exists since $\F_2^n$ is finite. By \cref{lem:pcp_linear_dependence_landscapes} there exist at least two indices $\mu,\rho$ with $\supp(\gamma_\mu)=\supp(\gamma_\rho)$. Take $k$ minimal such that $\supp(\gamma_j)=\supp(\gamma_k)$ for some $j<k$. Then $\gamma_0,\gamma_1,\ldots,\gamma_{k-1}$ are linearly independent, so $\deg m_{c,n}\geq k$. We now show that $\gamma_j=\gamma_k$. This gives $m_{c,n}(X)\mid X^j+X^k$ and hence $m_{c,n}(X)=X^j+X^k$. 
    
    By the assumption we know $\gamma_j,\gamma_k\neq \zero$. Furthermore, $A_q=S^q$ since $c$ is normed. We can write 
    \begin{equation*}
        \gamma_j = \cq\odot S^q\cq \odot\cdots\odot S^{(j-1)q}c = f\odot S^{(j-1)q}A_q = f\odot S^{jq}
    \end{equation*}
    and 
    \begin{equation*}
        \gamma_k = \cq\odot S^q\cq \odot\cdots\odot S^{(j-1)q}\cq\odot S^{jq}\cq\odot\cdots\odot S^{(k-1)q}c = f\odot g\odot S^{kq}
    \end{equation*}
    for some shift-invariant functions $f,g$. We observe that $jq\in\supp(\gamma_j)=\supp(\gamma_k)$. 

    Suppose first that $jq\in\supp(f\odot g)$. Then there exist $0\leq \xi\leq k-1$ and $\ell\in\supp(c)\setminus\{q\}$ such that $\xi q+\ell\equiv jq \bmod n$, implying that $\ell\in Q\leq \Z_n$. Because $\ell\in\supp(c)\setminus\{q\}\subseteq \mL_c$, we then have $\ell\in Q\cap \mL_c$. But by \cref{cor:pcp_min_poly_X^j_iff_exists_ell_in_Q_cap_L} $(c,n)$ is $\gamma$-nilpotent, a contradiction.

    Therefore, $jq\in\supp(S^{kq})=\{kq\}$, i.e. $S^{jq}=S^{kq}$. Then 
    \begin{equation*}
        \gamma_k = f\odot g\odot S^{kq} = f\odot S^{jq} \odot g = \gamma_j \odot g.
    \end{equation*}
    Since $\gamma_j\odot g\neq \zero$, 
    \begin{equation*}
        \supp(\gamma_j)\subseteq \supp(\gamma_j)\cup\supp(g) = \supp(\gamma_j\odot g) = \supp(\gamma_k)=\supp(\gamma_j).
    \end{equation*}
    It follows that $\supp(g)\subseteq \supp(\gamma_j)$. But then $\gamma_k=\gamma_j\odot g=\gamma_j$.

    We now know $m_{c,n}(X)=X^j+X^k = X^j(X^{k-j}+1)$. The only thing that is left to show is that $j\leq r$ and $k-j=r$. From $S^{jq}=S^{kq}$ we see that $jq\equiv kq \bmod n$, i.e. $(k-j)q\equiv 0 \bmod n$. As $r$ is the additive order of $q$ modulo $n$, we have $r\mid k-j$. On the other hand, by \cref{lem:pcp_min_poly_divides_X^2r_plus_X^r} we know that $m_{c,n}(X) \mid X^r(X^r+1)$. This implies $X^j\mid X^r$ and $X^{k-j}+1\mid X^r+1$, hence $j\leq r$ and $k-j\leq r$. This finishes the proof.
\end{proof}

The proof of \cref{lem:pcp_periodic_if_not_nilpotent} furthermore shows that the preperiod $j$ is the minimal positive integer $j$ such that $\supp(\gamma_j)=\supp(\gamma_{j+r})$. By \cref{lem:two_sided_support_gamma_i},
\begin{align*}
    \supp(\gamma_j)     &= \bigcup_{\zeta=0}^{j-1} (\zeta q+\supp(c)\setminus\{q\}) \cup \{jq\}, \\
    \supp(\gamma_{j+r}) &= \bigcup_{\zeta=0}^{j+r-1} (\zeta q+\supp(c)\setminus\{q\}) \cup \{(j+r)q\}.
\end{align*}
Since $r$ is the additive order of $q$ modulo n, $jq=(j+r)q \bmod n$. Moreover, the first union in $\supp(\gamma_j)$ is contained in the corresponding union for $\supp(\gamma_{j+r})$. Thus, it suffices to verify that
\begin{equation*}
    \bigcup_{\zeta=0}^{r-1} ((j+\zeta) q+\supp(c)\setminus\{q\}) \subseteq \bigcup_{\zeta=0}^{j-1} (\zeta q+\supp(c)\setminus\{q\}).
\end{equation*}
As $\zeta$ ranges from $0$ to $r-1$, both $\zeta q$ and $(j+\zeta)q$ range over the entire subgroup $Q$ generated by $q$, so that the above inclusion is equivalent to
\begin{equation*}
    \bigcup_{\zeta=0}^{r-1} (\zeta q+\supp(c)\setminus\{q\}) \subseteq \bigcup_{\zeta=0}^{j-1} (\zeta q+\supp(c)\setminus\{q\})
\end{equation*}
which can be efficiently used to determine $j$.

\begin{lemma}\label{lem:pcp_preperiod}
    Suppose $(c,n)$ is $\gamma$-periodic. The preperiod $j$ is the smallest positive integer $j$ with
    \begin{equation}\label{eq:pcp_defining_inclusion_preperiod}
        \bigcup_{\zeta=0}^{r-1} (\zeta q+\supp(c)\setminus\{q\}) \subseteq \bigcup_{\zeta=0}^{j-1} (\zeta q+\supp(c)\setminus\{q\}).
    \end{equation}
\end{lemma}

In conclusion, we obtain the following theorem. 

\begin{theorem}\label{thm:pcp_minimal_polynomial_conclusion}
    Let $c$ be a quasi-conserved landscape. 
    \begin{itemize}
        \item Let $V=\F_2^\Z$. Then $(c,V)$ is $\gamma$-nilpotent or $\gamma$-independent.
        \item Let $V=\F_2^n$. Then $(c,V)$ is $\gamma$-nilpotent or $\gamma$-periodic.
    \end{itemize}
\end{theorem}

\begin{example}\label{ex:pcp_example_0_01_min_polys}
    Consider the quasi-conserved landscape $c=*0\mt01$ with $q=4$. 
    \begin{enumerate}
        \item Let $V=\F_2^\Z$. By \cref{lem:right_sided_gammas_never_zero} we have $m_{c,\Z}(X)=0$.
        \item Let $V=\F_2^5$. By \cref{ex:pcp_0_01_nilpotent} $(*0\mt01,5)$ is nilpotent. We can check that $S^{4h}c\odot c^{(4)}=\zero$ with $h=1$, so that $m_{c,5}(X)=X^2$ by \cref{lem:pcp_index_of_nilpotency}.
        \item Let $V=\F_2^6$. Then $r=3$ and 
        \begin{equation*}
            \bigcup_{\zeta=0}^{r-1} (\zeta q+\supp(c)\setminus\{q\}) = \bigcup_{\zeta=0}^{2} (4\zeta+\{1,3\})=\{1,3,5\}\subseteq \Z_6.
        \end{equation*}
        Observe that 
        \begin{equation*}
            \bigcup_{\zeta=0}^{1} (4\zeta+\{1,3\})=\{1,3,5\}.
        \end{equation*}
        Then $m_{c,6}(X)=X^2(X^3+1)$ by \cref{lem:pcp_preperiod}.
    \end{enumerate}
\end{example}

\subsubsection{Right-sided landscapes}

In this section we compute more precise results about the minimal polynomials of right-sided quasi-conserved landscapes. We start with the infinite dimensional case.

\begin{lemma}\label{lem:pcp_min_poly_of_right_sided_over_Z}
    Let $c=*a_1\ldots a_m$ be a right-sided quasi-conserved landscape. Then $(c,\Z)$ is $\gamma$-independent.
\end{lemma}

\begin{proof}
    Combine \cref{lem:pcp_min_poly_F_2_Z} with \cref{lem:right_sided_gammas_never_zero}. 
\end{proof}

We now move on to the finite dimensional case. We make a case distinction depending on whether $m$ is a divisor of $n$ or not.

\begin{lemma}\label{lem:pcp_min_poly_right_sided_m_divides_n}
    Let $c=*a_1\ldots a_m$ be a right-sided quasi-conserved landscape and suppose $m\mid n$. Then $(c,n)$ is $\gamma$-periodic with period and preperiod $n/m$.
\end{lemma}

\begin{proof}
    Recall that by \cref{lem:right_sided_special_index_is_m} we know $q=m$. Therefore, $\gcd(q,n)=\gcd(m,n)=m$. Since $c$ is quasi-conserved we know $m\notin\mL_c$, so that $m\notin\mD_c$. By \cref{lem:pcp_char_of_n_which_are_nilpotent} we know that $(c,n)$ is not $\gamma$-nilpotent. Then by \cref{lem:pcp_periodic_if_not_nilpotent} $(c,n)$ is $\gamma$-periodic with period $r=n/m$ and preperiod $j\leq r$. We only need to show that $j=r$. 
    
    Suppose that $j<r$. Then 
    \begin{equation*}
        \bigcup_{\zeta=0}^{j-1} (\zeta q+\supp(c)\setminus\{q\}) \subseteq \{1,\ldots jm\}.
    \end{equation*}
    Take $\ell\in\supp(c)\setminus\{m\}$, then $jm < (r-1)m+\ell < rm=n$ because $c$ is right-sided and 
    \begin{equation*}
        (r-1)m+\ell\in \bigcup_{\zeta=0}^{r-1} (\zeta q+\supp(c)\setminus\{q\}),
    \end{equation*}
    a contradiction to \cref{lem:pcp_preperiod}. Therefore, $j=r$. 
\end{proof}

This gives a complete characterization of the minimal polynomial of right-sided quasi-conserved landscapes with full support. As a special case this determines the minimal polynomial of $*01, *0\ldots 01$ and more generally all landscapes given by a bifix-free string (see \cref{rem:right_sided_bifix_free}), i.e. the results of \cite{kriepke2024algebraic,kriepke2026siblings,lyu2026generalized,zhao2026bifix}.

\begin{lemma}\label{lem:pcp_min_poly_if_full_support}
    Let $c=*a_1\ldots a_m$ be a right-sided quasi-conserved landscape with full support, i.e. $\supp(c)=\{1,\ldots,m\}$. Then 
    \begin{equation*}
        m_{c,V}(X) = 
            \begin{cases}
                0,                   & V=\F_2^\Z, \\
                X^{n/m}(X^{n/m}+1),  & V=\F_2^n, m\mid n, \\
                X^{\ceil{n/m}},      & \text{otherwise.}
            \end{cases}
    \end{equation*}
\end{lemma}

\begin{proof}
    If $V=\F_2^\Z$ or $m\mid n$ we are done. Hence, suppose that $m\nmid n$. We need to show that $\gamma_j=\zero$ for some $j$ and the minimal such $j$ is $\ceil{n/m}$. 
    
    We first show that $j\leq \ceil{n/m}$. Take $0<\ell<m$ which is determined by $n\equiv \ell \mod m$. It holds that
    \begin{equation*}
        (\ceil{n/m}-1)m+\ell=n. 
    \end{equation*} 
    As $\supp(c)=\{1,\ldots,m\}$ we have $\ell\in\supp(c)\setminus\{m\}$ and hence $S^\ell \cm\odot c=\zero$ by \cref{rem:right_sided_no_diff_between_S_ell_cm_odot_c_and_S_ell_c_odot_c}. Therefore,
    \begin{equation*}
        \gamma_{\ceil{n/m}} = \cm \odot\cdots\odot S^{(\ceil{n/m}-1)m}c = \cm \odot S^{-\ell}c \odot(\ldots) = S^{-\ell}(S^\ell \cm \odot c)\odot(\ldots) = \zero.
    \end{equation*}
    
    We now need to show that $j\geq \ceil{n/m}$. By \cref{lem:right_sided_gammas_never_zero} we know that $\gamma_j\neq \zero$ for all $j\geq 0$ if we consider $\gamma_j$ as a function on $\F_2^\Z$. Hence, for $\gamma_j$ to vanish over $\F_2^n$, we need to have at least one element $\ell\in\supp(\gamma_j)$ with $\ell>n$, so that after reducing the exponents of $S$ modulo $n$, we obtain a product $S^\ell\odot(S^\ell+\one)=\zero$. As $\supp(\gamma_j)\subseteq \{1,\ldots,jm\}$, we need at least $j>n/m$ for $\gamma_j=\zero$ to hold. As $j$ is a natural number this implies $j\geq \ceil{n/m}$ and we are done.
\end{proof}

\begin{remark}
    The condition that $c$ has full support in \cref{lem:pcp_min_poly_if_full_support} is needed in general. Indeed, we saw in \cref{ex:pcp_example_0_01_min_polys} that $m_{*0\mt01,6}(X)=X^2(X^3+1)$ instead of $X^2$. The proof shows more generally that the lemma holds for all $n\in\N$ for which there exists $\ell$ such that $n\equiv \ell \mod m$ and $S^\ell \cm\odot c=\zero$.
\end{remark}

\subsection{Shift-invariant permutations in finite dimensions}

Using the results of the previous section we can now study the functions in $\mG$. By \cref{thm:pcp_G_isomorphic_M_monoid} we know that $\mG\cong \mM \subseteq \F_2[X]/(m_{c,V}(X))$ via the isomorphism $\mG\to\mM$ given by 
\begin{equation*}
    \varphi\left(\sum_{i=0}^s \alpha_i \gamma_i\right) = \left[\sum_{i=0}^s \alpha_i X^i\right].
\end{equation*}
By \cref{thm:pcp_minimal_polynomial_conclusion} the only possible minimal polynomials have the form
\begin{equation*}
    m_{c,V}(X) \in \{0, X^j, X^j(X^r+1)\}
\end{equation*}
for some positive integers $j,r$, implying that the study of $\mG$ is essentially equivalent to the study of the rings
\begin{equation*}
    \F_2[X]/(X^j), \qquad \F_2[X]/(X^j(X^r+1))\qquad  \text{ and }\qquad \F_2[X]/(0)\cong\F_2[X]. 
\end{equation*}
If $m_{c,V}(X)\neq 0$, then by \cref{rem:pcp_unit_iff_permutation_if_G_finite} a function $f\in\mG$ is a permutation on $V$ if and only if $\varphi(f)\in\mR$ is a unit. Recall that the unit group of $\mR=\F_2[X]/(m(X))$ is given by 
\begin{equation*}
    \mR^\times = \{[F(X)]: \gcd(F(X),m(X))=1\}
\end{equation*}
We obtain a full classification of the permutations from $\mG$ in the case that the minimal polynomial $m_{c,V}(X)$ is nonzero.

\begin{theorem}\label{thm:pcp_full_classification_permutations_in_G_nonzero_minpoly}
    Let $c$ be a quasi-conserved landscape. 
    \begin{enumerate}
        \item Suppose $m_{c,V}(X)=X^j$. Then any $f\in \mG$ is a permutation.
        \item Suppose $m_{c,V}(X)=X^j(X^r+1)$. Then $f\in \mG$ is a permutation if and only if $F(X)\in\F_2[X]$ given by $\varphi(f)=[F(X)]$ is co-prime with $X^r+1$.
    \end{enumerate}
\end{theorem}

\begin{proof}
    As $\mG^\times\cong \mM^\times=\mM\cap \mR^\times$, we only need to study $\mR^\times$ for each of the two cases. 
    \begin{enumerate}
        \item Set $\varphi(f)=[F(X)]$. Then $[F(X)]\in\mR^\times$ if and only if $\gcd(F(X),X^j)=1$. This is equivalent to $F(0)\neq 0$, which is satisfied because $F(X)$ has constant coefficient $1$. Therefore, $[F(X)]$ is a unit and $f\in\mG$ is a permutation. 
        \item $[F(X)]\in\mR^\times$ if and only if $\gcd(F(X),X^j(X^r+1))=1$. As $F(X)$ has constant coefficient $1$ this is equivalent to $\gcd(F(X),X^r+1)=1$.
    \end{enumerate}
\end{proof}

In this section we focus mostly on the case $V=\F_2^n$, which is completely covered by \cref{thm:pcp_full_classification_permutations_in_G_nonzero_minpoly} by \cref{thm:pcp_minimal_polynomial_conclusion}. Nonetheless, some of the results also apply to the case $V=\F_2^\Z$.

If $(c,V)$ is $\gamma$-nilpotent, then \cref{thm:pcp_full_classification_permutations_in_G_nonzero_minpoly} states that every element in $\mG$ is a permutation. Since $\mG$ is a monoid it follows that $\mG$ is a group. In fact, $\mG\cong \mM=\mM^\times=\mR^\times$ is the unit group of $\mR$. By a simple counting argument we know $\abs{\mM}=2^{j-1}$, which means that there are $2^{j-1}$ permutations in $\mG$. We record this in the following lemma.

\begin{lemma}\label{lem:pcp_unit_group_if_min_poly_X^j}
    Let $\mR=\F_2[X]/(X^j)$. Then $\mR^\times = \mM$ and $\abs{\mR^\times} = 2^{j-1}$. 
\end{lemma}

The case $m_{c,V}(X)=X^j(X^r+1)$ is more involved. The factorization of $X^r+1$ is a well-known problem, see for example \cite[Chapter 2, Section 4]{Lidl:1997}. We also refer to \cite{graner2024irreducible} for the more general problem of factorizing $X^r-a$ over any finite field. Write $r=2^s t$ with $t$ odd. Then $X^r+1=(X^t+1)^{2^s}$ and the factorization of $X^r+1$ follows from the factorization of $X^t+1$. The roots of $X^t+1$ are called $t$-th roots of unity and form a cyclic subgroup of the splitting field of $X^t+1$ over $\F_2$. A generator $\alpha$ of that subgroup is called a primitive $t$-th root of unity. The $t$-th cyclotomic polynomial $\Phi_t(X)$ is given by 
\begin{equation*}
    \Phi_t(X) = \prod_{\substack{i=1 \\ \gcd(i,t)=1}}^t (X-\alpha^i)
\end{equation*}
and is independent of the specific choice of $\alpha$. It splits into $\phi(t)/d$ distinct irreducible polynomials $\rho_{t,1},\ldots,\rho_{t,\phi(t)/d}$ of degree $d$, where $\phi(t)$ is the Euler totient function and $d=\ord_t(2)$ is the multiplicative order of $2$ modulo $t$, i.e. the minimal positive integer with $2^d \equiv 1 \mod t$. The factorization of $X^t+1$ is then
\begin{equation*}
    X^t+1 = \prod_{\ell\mid t} \Phi_{\ell}(X) = \prod_{\ell\mid t} \prod_{i=1}^{\phi(\ell)/\ord_\ell(2)} \rho_{\ell,i}.
\end{equation*}
Finally, 
\begin{equation*}
    X^r+1 = \prod_{\ell\mid t} \prod_{i=1}^{\phi(\ell)/\ord_\ell(2)} \rho_{\ell,i}^{2^s}.
\end{equation*}
By the Chinese Remainder Theorem
\begin{equation}\label{eq:pcp_crt_on_F_2_X_mod_X^jX^r_plus_1}
    \F_2[X]/(X^j(X^r+1)) \cong \F_2[X]/(X^j) \times \prod_{\ell\mid t} \prod_{i=1}^{\phi(\ell)/\ord_\ell(2)} \F_2[X]/( \rho_{\ell,i}^{2^s}),
\end{equation}
i.e. the unit group of $\F_2[X]/(X^j(X^r+1))$ is isomorphic to the product of the unit groups on the right-hand side, which are easier to study. The following lemma describes their size.

\begin{lemma}\label{lem:pcp_unit_group_size_modulo_irreducible_to_power_n}
    Let $\rho(X)\in\F_2[X]$ be an irreducible polynomial of degree $d$. Then the unit group of $\mR=\F_2[X]/(\rho^w(X))$ has $2^{d(w-1)}(2^d-1)$ many elements.
\end{lemma}

\begin{proof}
    An element $[F(X)]\in\mR$ is a unit if and only if $\gcd(F(X),\rho^w(X))=1$. As $\rho(X)$ is irreducible this is equivalent to $\gcd(F(X),\rho(X))=1$, i.e. $\rho(X)$ is not a divisor of $F(X)$. As $\rho^w(X)$ has degree $dw$, each element in $\mR$ corresponds uniquely to a polynomial of degree at most $dw-1$ and the multiples of $\rho(X)$ are exactly those polynomials of the form $\rho (X)\theta(X)$ with $\deg(\theta) < dw-d=d(w-1)$. Since there are $2^{d(w-1)}$ such polynomials $\theta(X)$, there are $2^{wd}-2^{d(w-1)}=2^{d(w-1)}(2^d-1)$ polynomials which are not multiples of $\rho(X)$.
\end{proof}

\begin{lemma}
    Set $\mR = \F_2[X]/(X^j(X^r+1))$ and write $r=2^st$ with $t$ odd. Then
    \begin{equation*}
        \abs{\mR^\times} = 2^{j-1+(2^s-1)t} \prod_{\ell\mid t} \left(2^{\ord_\ell(2)}-1\right)^{\phi(\ell)/\ord_\ell(2)}.
    \end{equation*}
\end{lemma}

\begin{proof}
    By \cref{eq:pcp_crt_on_F_2_X_mod_X^jX^r_plus_1} combined with \cref{lem:pcp_unit_group_if_min_poly_X^j} and \cref{lem:pcp_unit_group_size_modulo_irreducible_to_power_n} we have
    \begin{align*}
        \abs{\mR^\times} 
            &= 2^{j-1} \prod_{\ell\mid t} \prod_{i=1}^{\phi(\ell)/\ord_\ell(2)} 2^{\ord_\ell(2)(2^s-1)}(2^{\ord_\ell(2)}-1) \\
            &= 2^{j-1} \prod_{\ell\mid t} 2^{\phi(\ell)(2^s-1)} (2^{\ord_\ell(2)}-1)^{\phi(\ell)/\ord_\ell(2)} \\
            &= 2^{j-1+(2^s-1)\sum_{\ell\mid t} \phi(\ell)}  \prod_{\ell\mid t} (2^{\ord_\ell(2)}-1)^{\phi(\ell)/\ord_\ell(2)}.
    \end{align*}
    Using the well-known identity $\sum_{\ell\mid t}\phi(\ell)=t$ finishes the proof. 
\end{proof}

As a rough estimate we note that
\begin{equation*}
    \prod_{\ell\mid t} (2^{\ord_\ell(2)}-1)^{\phi(\ell)/\ord_\ell(2)} 
        \approx \prod_{\ell\mid t} 2^{\phi(\ell)} = 2^t,
\end{equation*}
so that
\begin{equation*}
    \abs{\mR^\times} \approx 2^{j-1+(2^s-1)t+t} = 2^{j-1+2^st} = 2^{j+r-1}.
\end{equation*}

If $[F(X)]\in\F_2[X]/(m(X))$ is invertible, then the inverse $[F(X)]^{-1}$ can be efficiently calculated using the Extended Euclidean algorithm. This gives an easy tool to compute the inverse of any permutation $f\in\mG$ for any given landscape $c$. A particular example are binomials $[1+X^i]$. Note that these correspond to the functions $\id+\gamma_i\in\mG$, in particular for $i=1$ they correspond to $\id+c$. Observe that $1+X$ is a common divisor of $1+X^i$ and $X^r+1$, implying that $\gcd(1+X^i, X^r+1)\neq 1$. By \cref{thm:pcp_full_classification_permutations_in_G_nonzero_minpoly} it follows that $\gamma_0+\gamma_i$ is a permutation on $\F_2^n$ if and only if $(c,n)$ is $\gamma$-nilpotent. We can then compute their inverse. 

\begin{lemma}\label{lem:pcp_formula_for_inverse_1_plus_X^i}
    Let $(c,V)$ be $\gamma$-nilpotent of index $j$. Let $i\geq 1$ and $s=\max\{t\in\N: i\mid t\leq j\}$ be the largest multiple of $i$ which does not exceed $j$. Take $f=\gamma_0+\gamma_i\in\mG$. Then the inverse of $f$ is 
    \begin{equation*}
        f^{-1} = \gamma_0+\gamma_i+\gamma_{2i}+\cdots+\gamma_s.
    \end{equation*}
\end{lemma}

\begin{proof}
    We have $\varphi(f)=[1+X^i]\in\F_2[X]/(X^j)$. By \cref{lem:pcp_unit_group_if_min_poly_X^j} combined with Lagrange's Theorem we know that $[1+X^i]^{2^{j-1}}=[1]$ which implies that $[1+X^i]^{-1}=[1+X^i]^{2^{j-1}-1}$. We obtain
    \begin{align*}
        [1+X^i]^{2^{j-1}-1} &= \frac{[1+X^i]^{2^{j-1}}}{[1+X^i]} = \frac{\left[1+X^{2^{j-1}i}\right]}{[1+X^i]} \\
            &= [1+X^i+X^{2i}+\cdots + X^{(2^{j-1}-1)i}] \\
            &= [1+X^i+X^{2i}+\cdots + X^s].
    \end{align*}
    By applying $\psi$ we get
    \begin{equation*}
        f^{-1} = \psi\left([1+X^i]^{-1}\right) = \gamma_0+\gamma_i+\gamma_{2i}+\cdots+\gamma_s.
    \end{equation*}
\end{proof}

As a particular instance of \cref{lem:pcp_formula_for_inverse_1_plus_X^i} for the specific choice of the landscapes $c=*01$ and $c=*0\ldots01$ we obtain the inverses of $\chi$ and $\chi^{(m)}$ which were found in \cite{kriepke2024algebraic,lyu2026generalized}. The case $i=1$ for right-sided landscapes with full support, i.e. landscapes given by bifix-free strings, is covered in \cite{zhao2026bifix}. The inverse of $\chi$ was first found in \cite{liu2022inverse} using different tools.

\begin{remark}\label{rem:pcp_elementary_observation_to_compute_inverse}
    The following elementary observation gives a different way to compute the inverse 
    \begin{equation*}
        [G(X)]=[F(X)]^{-1}\in\F_2[X]/(m(X)), 
    \end{equation*}
    assuming it exists. By the definition of the inverse we know that
    \begin{equation*}
        F(X) G(X) = 1 + t(X) m(X)
    \end{equation*}
    for some polynomial $t(X)\in\F_2[X]$. But then $[t(X)] [m(X)] =[1]\in\F_2[X]/(F(X))$. Since we assumed that $\gcd(F(X),m(X))=1$, we then have $[t(X)]=[m(X)]^{-1}\in\F_2[X]/(F(X))$. After finding $t(X)$ we may compute $G(X)$ by 
    \begin{equation*}
        G(X) = \frac{1+t(X)m(X)}{F(X)},
    \end{equation*}
    which is an exact division in $\F_2[X]$. The advantage here is that computing the inverse $[m(X)]^{-1}\in\F_2[X]/(F(X))$ is sometimes simpler than computing $[F(X)]^{-1}\in \F_2[X]/(m(X))$, especially if $F(X)$ has a small degree. Note that in particular $\deg t(x) < \deg F(X)$.
\end{remark}

Recall from \cref{subsec:prelim_finite_dim} the set
\begin{equation*}
    \inv(f) = \{n\in\N: f\text{ is a permutation on }\F_2^n\},
\end{equation*}
which was introduced in \cite{omland2022permutation}. \cref{lem:prelim_inv_f_is_lower_set} shows that the set $\inv(f)$ is closed under taking divisors. Furthermore, the set $\N\setminus\inv(f)$ has a generating set $\xi(f)$: it holds that $n\in\N\setminus\inv(f)$ if and only if $n$ is a multiple of some element from $\xi(f)$. Determining $\inv(f)$ for a general shift-invariant function is hard. For functions $f\in\mG$ we are able to fully characterize $\xi(f)$ and $\inv(f)$. Recall the set 
\begin{equation*}
    \mD_c = \{d: d\mid q, \exists \ell\in\mL_c \text{ s.t. } d\mid \ell\}
\end{equation*}
from \cref{lem:pcp_char_of_n_which_are_nilpotent}, where $q$ is the special index of $c$. Next, we recall the concept of order of a polynomial. Let $F(X)\in\F_2[X]$ be a polynomial with $F(0)=1$. We say that the order of $F$ is the smallest natural number $r$ such that $F(X) \mid X^r+1$. We write $r=\ord(F)$. If $F(X)$ is irreducible, then $F(X)\mid X^\ell+1$ if and only if $\ord(F) \mid \ell$. For more information on the order of a polynomial, see for example \cite{Lidl:1997}.

\begin{theorem}\label{thm:pcp_xi_of_f}
    Let $f\in\mG$ with $\varphi(f)=[F(X)]$ and $F=g_1^{e_1}\cdots g_t^{e_t}\in\F_2[X], e_i\geq 1$ be its factorization into irreducible factors. Then
    \begin{equation*}
        \xi(f) = \{d\cdot \ord(g_i): d\mid q\text{ s.t. } d\notin\mD_c, i\in\{1,\ldots,t\}\}.
    \end{equation*}
    In particular, $\xi(f)$ is finite and $\inv(f)$ is infinite. 
\end{theorem}

\begin{proof}
    Suppose first that $d=\gcd(q,n)\in\mD_c$. We show that $n$ is not a multiple of an element in $\xi(f)$ and that $f$ is a permutation on $\F_2^n$. Assume otherwise, i.e. $d'\ord(g_i)\mid n$ for some $d'\mid q$ and $i\in\{1,\ldots,t\}$. Then $d'\mid \gcd(q,n)=d$. Because $d\in\mD_c$ there exists $\ell\in \mL_c$ such that $d\mid \ell$. But then $d' \mid d\mid \ell$, a contradiction to the assumption that $d'\notin\mD_c$. To finish this case, observe that by \cref{lem:pcp_char_of_n_which_are_nilpotent} $(c,n)$ is $\gamma$-nilpotent and thus by \cref{thm:pcp_full_classification_permutations_in_G_nonzero_minpoly} $f$ is a permutation on $\F_2^n$.

    Suppose now that $\gcd(q,n)=d\notin \mD_c$. Then $(c,n)$ is $\gamma$-periodic with period $n/d$ by \cref{lem:pcp_char_of_n_which_are_nilpotent} combined with \cref{lem:pcp_periodic_if_not_nilpotent}. By \cref{thm:pcp_full_classification_permutations_in_G_nonzero_minpoly} $f$ is a permutation on $\F_2^n$ if and only if $\gcd(F(X),X^{n/d}+1)=1$. Using the factorization $F=g_1^{e_1}\cdots g_t^{e_t}$ this holds if and only if $\gcd(g_i(X),X^{n/d}+1)=1$ for all $i=1,\ldots,t$. As $g_i$ is irreducible this means that $g_i$ is not a divisor of $X^{n/d}+1$ which is equivalent to $\ord(g_i)$ not being a divisor of $n/d$. Hence, $f$ being a permutation is equivalent to $n$ not being a multiple of $d\cdot\ord(g_i)$ for any $i=1,\ldots,t$ and the claim is proven.
\end{proof}

We give several corollaries of \cref{thm:pcp_xi_of_f}.

\begin{corollary}\label{cor:pcp_xi_of_id_plus_c}
    It holds that
    \begin{equation*}
        \xi(\id+c) = \{d: d\mid q\text{ s.t. } d\notin\mD_c\}
    \end{equation*}
    and 
    \begin{equation*}
        \inv(\id+c) = \{n\in\N: \gcd(q,n)\in \mD_c\}.
    \end{equation*}
\end{corollary}

\begin{proof}
    Follows immediately from $\varphi(\id+c)=[1+X]$ and $\ord(1+X)=1$. 
\end{proof}

\begin{example}
    Let $v>0$ and define $\chi_{n,v}=\id+(\one+S^v)\odot S^{2v}$ on $\F_2^n$. This map is a generalization of $\chi$ which is considered in \cite{feng2026generalization}, where it is shown that it is invertible on $\F_2^n$ if and only if $2^k\mid v$, where $n=2^k n_0$ with $n_0$ odd. We show how to obtain this result with our tools. 

    Set $c=(\one+S^v)\odot S^{2v}=*\mt\ldots\mt0\mt\ldots\mt1$, then $\supp(c)=\{v,2v\}$. It is easy to show that $(c,2v)$ is quasi-conserved with $\mL_c=\{\pm v\}$. Therefore,
    \begin{equation*}
        \mD_c = \{d: d\mid 2v \text{ s.t. } d\mid v\} = \{d: d\mid v\}
    \end{equation*}
    is the set of divisors of $v$. By \cref{cor:pcp_xi_of_id_plus_c}, $f$ is invertible on $\F_2^n$ if and only if $\gcd(2v,n)\mid v$. Write $v=2^\ell v_0$ with $v_0$ odd. We show that $\gcd(2v,n)\mid v$ if and only if $k\leq \ell$, which is equivalent to $2^k\mid v$. 

    It is clear that 
    \begin{equation*}
        \gcd(2v,n) = \gcd(2^{\ell+1}v_0,2^k n_0) = 2^{\min(\ell+1,k)}\gcd(2^{\ell+1-\min(\ell+1,k)}v_0,2^{k-\min(\ell+1,k)}n_0).
    \end{equation*}
    Suppose first that $\gcd(2v,n)\mid v=2^\ell v_0$. It follows that $\min(\ell+1,k)\leq \ell$ and hence, $k\leq \ell$. Suppose now that $k\leq \ell$, then $\min(\ell+1,k)=k$ and 
    \begin{equation*}
        \gcd(2v,n)=2^k \gcd(2^{\ell+1-k}v_0, n_0) = 2^k \gcd(v_0, n_0) \mid v.
    \end{equation*}
    The same calculation also applies to the map $\tilde{f}=\id+(\one+S^{-2v})\odot S^{-v}$ considered in \cite{feng2026generalization}. Observe that $(\one+S^{-2v})\odot S^{-v}$ and $(\one+S^v)\odot S^{2v}$ are conjugate equivalent. 
\end{example}

\begin{corollary}\label{cor:pcp_xi_of_f_for_right_sided_full_supp}
    Let $c=*a_1\ldots a_m$ be a right-sided quasi-conserved landscape with full support. Let $f\in\mG$ with $\varphi(f)=[F(X)]$ and $F=g_1^{e_1}\cdots g_t^{e_t}\in\F_2[X], e_i\geq 1$ be its factorization into irreducible factors. Then
    \begin{equation*}
        \xi(f) = \{m\cdot \ord(g_i): i\in\{1,\ldots,t\}\}.
    \end{equation*}
\end{corollary}

\begin{proof}
    Since $c$ has full support we have $\{1,\ldots,m-1\}\subseteq \mL_c$. Recall that $q=m$ for right-sided landscapes. It holds that $\mD_c = \{d\mid q: d<m\}$ is the set of all proper divisors of $m$. The result follows by applying \cref{thm:pcp_xi_of_f}.
\end{proof}

\begin{corollary}\label{cor:pcp_xi_of_id_plus_c_full_support}
    Let $c=*a_1\ldots a_m$ be a right-sided quasi-conserved landscape with full support. Then
    \begin{equation*}
        \xi(\id+c) = \{m\}.
    \end{equation*}
\end{corollary}

\begin{proof}
    Combine \cref{cor:pcp_xi_of_id_plus_c} and \cref{cor:pcp_xi_of_f_for_right_sided_full_supp}.
\end{proof}

\begin{example}
    Let $F(X)=1+X+X^2$. The associated map $\kappa=\gamma_0+\gamma_1+\gamma_2$ was first studied in \cite{kriepke2026siblings} for the landscape $*01$. The polynomial $F$ is irreducible and $\ord(F)=3$, since $(1+X)(1+X+X^2)=1+X^3$. By \cref{thm:pcp_xi_of_f} it follows that
    \begin{equation*}
        \xi(\gamma_0+\gamma_1+\gamma_2) =  \{3d: d\mid q\text{ s.t. } d\notin\mD_c\}.
    \end{equation*}
    If $c$ is right-sided with full support,
    \begin{equation*}
        \xi(\gamma_0+\gamma_1+\gamma_2) = \{3m\}.
    \end{equation*}
    By \cref{rem:pcp_elementary_observation_to_compute_inverse}, the inverse of $[1+X+X^2]\in\F_2[X]/(m_{c,V}(X))$ is given by 
    \begin{equation*}
        [1+X+X^2]^{-1} = \left[\frac{1+t(X)m_{c,V}(X)}{1+X+X^2}\right],
    \end{equation*}
    with $[t(X)] = [m_{c,V}(X)]^{-1}\in\F_2[X]/(1+X+X^2)$. In particular, $t(X)\in\{1,X,1+X\}$. One can quickly show that if $m_{c,V}(X)=X^j$, then
    \begin{equation*}
        t(X) = 
            \begin{cases}
                1,   & j\equiv 0 \bmod 3, \\
                1+X, & j\equiv 1 \bmod 3, \\
                X,   & j\equiv 2 \bmod 3,
            \end{cases}
    \end{equation*}
    and if $m_{c,V}(X)=X^j+X^k$, then
    \begin{equation*}
        t(X) = 
            \begin{cases}
                1,   & (j,k) \equiv (1,2), (2,1) \bmod 3, \\
                X,   & (j,k) \equiv (0,1), (1,0) \bmod 3, \\
                1+X, & (j,k) \equiv (0,2), (2,0) \bmod 3.
            \end{cases}
    \end{equation*}
    An explicit form of the inverse is given in \cite{kriepke2026siblings} for the specific case $k=2j$. 
    Note that if $j\equiv k\bmod 3$ then $1+X+X^2\mid 1+X^3 \mid 1+X^{k-j}$, so in that case $[1+X+X^2]\in\F_2[X]/(X^j+X^k)$ is not invertible. 
\end{example}

Another important property of a permutation is its order, i.e. the minimal positive integer $s$ such that $f^s = \id$. It is clear that the order of $f\in\mG$ is equal to the multiplicative order of $\varphi(f)=[F(X)]$ in $\mR$. The multiplicative order of an element $[F(X)]\in\F_2[X]/(X^j(X^r+1))$ is in general difficult to determine, as it depends on the factorization of $X^r+1$. We therefore focus on the case $[F(X)]\in\F_2[X]/(X^j)$.

\begin{lemma}\label{lem:pcp_order_of_F_in_R}
    Let $[F(X)]=[1+X^\ell+\sum_{i=\ell+1}^s \alpha_i X^i]\in \F_2[X]/(X^j)$. Then the multiplicative order of $[F(X)]$ is given by $2^\rho$ with $2^{\rho-1}<\frac{j}{\ell}\leq 2^\rho$.
\end{lemma}

\begin{proof}
    By \cref{lem:pcp_unit_group_if_min_poly_X^j} the unit group of $\F_2[X]/(X^j)$ has $2^{j-1}$ elements. By Lagrange's theorem the order of every unit is a power of $2$. Let $\rho$ be arbitrary, then
    \begin{equation*}
        [F(X)]^{2^\rho} = \left[1+X^\ell+\sum_{i=\ell+1}^s \alpha_i X^i\right]^{2^\rho} = \left[1+X^{2^\rho \ell}+\sum_{i=\ell+1}^s \alpha_i X^{2^\rho i}\right]
    \end{equation*}
    and therefore $[F(X)]^{2^\rho}=[1]$ if and only if $2^\rho \ell \geq j$. As we are interested in the minimal such $\rho$, we obtain the claim.
\end{proof}

Related to the order of a permutation is a general formula for all its iterates $f^t$ for $t\geq 0$. However, one cannot generally expect such a formula. For polynomials $F(X)$ it is possible to give a formula for its powers $F(X)^t$. This allows to give expressions for all iterates $f^t$ for all functions $f\in\mG$. We restrict ourselves here to the simplest case of $F(X)=1+X^i$ being a binomial. For non-negative integers $t=t_0+2t_1+2^2t_2+\cdots+2^st_s$ and $h=h_0+2h_1+2^2h_2+\cdots+2^sh_s$ written in base $2$, we write $h\preceq t$ if $h_i\leq t_i$ for all $i=0,\ldots,s$ and say that $h$ is covered by $t$. We then obtain the following result.

\begin{lemma}\label{lem:pcp_formula_for_iterates}
    Let $f=\gamma_0+\gamma_i\in\mG$. Then
    \begin{equation*}
        f^t = \sum_{h=0}^t \alpha_h \gamma_{ih}
    \end{equation*}
    with $\alpha_h=1$ if and only if $h\preceq t$. 
\end{lemma}

\begin{proof}
    It holds that $\varphi(f)=[1+X^i]$. By the Binomial theorem,
    \begin{equation*}
        [1+X^i]^t = \left[\sum_{h=0}^t \binom{t}{h} (X^i)^h\right]  = \left[\sum_{h=0}^t \binom{t}{h} X^{ih}\right].
    \end{equation*}
    Lucas' theorem states that $\binom{t}{h}$ is odd if and only if $h\preceq t$. Hence, 
    \begin{equation*}
        f^t = \psi\left([1+X^i]^t\right) = \psi\left(\left[\sum_{h=0}^t \binom{t}{h} X^{ih}\right]\right) = \sum_{h=0}^t \alpha_h \gamma_{ih}.
    \end{equation*}
\end{proof}

Until now, we have mostly used that 
\begin{equation*}
    \mG \cong \mM \subseteq \F_2[X]/(m_{c,V}(X))
\end{equation*}
as a monoid, so that the multiplicative structure of $\mR=\F_2[X]/(m_{c,V}(X))$ allows us to study the unit group $\mG^\times$. However, we also know that
\begin{equation*}
    \Gamma \cong \mR
\end{equation*}
as a vector space, giving us control about the additive behavior of $\mG = \gamma_0+\spn\{\gamma_j: j\geq 1\}\subseteq \Gamma$. For an interesting application of this structure we refer to \cite{krompholz2026cmtriples}.

\subsection{Cellular automata}

We now consider the case $V=\F_2^\Z$ and start by stating the result corresponding to \cref{thm:pcp_full_classification_permutations_in_G_nonzero_minpoly} in the infinite-dimensional case. We repeat the $\gamma$-nilpotent case here for completeness.

\begin{theorem}\label{thm:pcp_full_classification_permutations_in_G_infinite_dims}
    Let $c$ be a quasi-conserved landscape. Then:
    \begin{enumerate}
        \item Suppose $m_{c,V}(X)=X^j$. Then any $f\in \mG$ is a permutation.
        \item Suppose $m_{c,V}(X)=0$. Then the only permutation in $\mG$ is the identity.
    \end{enumerate}
\end{theorem}

\begin{proof}
    Statement $(a)$ is the same as in \cref{thm:pcp_full_classification_permutations_in_G_nonzero_minpoly}. 

    Now suppose $m_{c,V}(X)=0$ and take $f\in\mG$ with $f\neq \id$. By \cref{lem:pcp_nilpotent_on_Z_iff_nilpotent_for_every_n} we know that there exists some $n\in\N$ for which $(c,n)$ is $\gamma$-periodic. \cref{lem:pcp_char_of_n_which_are_nilpotent} then gives a divisor $d\mid q$ with $d\notin\mD_c$. By \cref{thm:pcp_xi_of_f} the set $\xi(f)$ is nonempty, implying that there exists some $\hat{n}\in\N$ for which $f$ is not a permutation on $\F_2^{\hat{n}}$. But then $f:\F_2^{\hat{n}}\to\F_2^{\hat{n}}$ is not injective and hence $f:\F_2^\Z\to\F_2^\Z$ is also not injective.
\end{proof}

\begin{example}
    Let $c=0\mt_r*1\mt_r0 = (\one+S^{-r})\odot S\odot (\one+S^{r+1})$ for some $r\geq 1$ from \cref{tb:landscapes_from_daemen}, i.e. Table A.1 from Daemen \cite{daemen1995cipher}. Then $(c,1)$ is quasi-conserved and it is easy to see that $S^{(r+1)\cdot 1}c\odot c^{(1)}=\zero$ and $h=r+1$ is the minimal such number. Therefore, $(c,\Z)$ is $\gamma$-nilpotent of index $r+2$. By \cref{lem:pcp_unit_group_if_min_poly_X^j} we obtain $2^{r+2-1}=2^{r+1}$ permutations on $\F_2^\Z$.  
\end{example}

There exist landscapes $c$ with $m_{c,\Z}(X)=0$, for example all right-sided landscapes by \cref{lem:pcp_min_poly_of_right_sided_over_Z}. By \cref{thm:pcp_full_classification_permutations_in_G_infinite_dims} these do not give any (nontrivial) permutations on $\F_2^\Z$. However, we can show that for every quasi-conserved landscapes all functions $f\in\mG$ are always surjective on $\F_2^\Z$. This was previously proved for the specific case of $f=\chi$ in \cite{schoone2024state}. 

Before showing this claim we discuss some general results about shift-invariant functions $f:\F_2^\Z\to\F_2^\Z$. Recall that if $f_0$ is a local map, i.e. it has the form $f_0(x)=f_0(x_{-k},\ldots,x_m)$ in \cref{eq:prelim_f_0_is_local_map}, then $f$ is also called a cellular automaton. For a thorough introduction to cellular automata we refer to the survey \cite{kari2005CA} by Kari. We only give a proof sketch of the following result and for more details refer to \cite{kari2005CA} and the references therein.

\begin{lemma}\label{lem:pcp_cellular_automata_surjective_iff_permutation_for_infinite_n}
    Let $f:\F_2^\Z\to\F_2^\Z$ be a cellular automaton and $\mN\subseteq \N$ be an infinite subset. If $f:\F_2^n\to\F_2^n$ is surjective for all $n\in\mN$, then $f:\F_2^\Z\to\F_2^\Z$ is surjective.
\end{lemma}

\begin{proof}
    The space $\F_2^\Z$ equipped with the product topology is compact. Define
    \begin{equation*}
        \mP_\mN = \bigcup_{n\in \mN} \mP_n = \bigcup_{n\in\mN} \{x\in\F_2^\Z: x_i = x_{i+n} \text{ for all } i\in\Z\}.
    \end{equation*}
    Then $\mP_\mN$ is dense in $\F_2^\Z$ since $\mN$ is infinite and $f:\mP_\mN\to\mP_\mN$ is surjective. By the Curtis-Hedlund-Lyndon theorem $f$ is continuous. A standard compactness argument shows that $f:\F_2^\Z\to\F_2^\Z$ is surjective. 
\end{proof}

\begin{theorem}\label{thm:pcp_all_f_are_surjective}
    Let $f\in\mG$ be arbitrary. Then $f:\F_2^\Z\to\F_2^\Z$ is surjective.
\end{theorem}

\begin{proof}
    By definition of the set $\inv(f)$ the function $f:\F_2^n\to\F_2^n$ is surjective for every $n\in\inv(f)$ and by \cref{thm:pcp_xi_of_f} the set $\inv(f)$ is infinite. \cref{lem:pcp_cellular_automata_surjective_iff_permutation_for_infinite_n} finishes the proof.
\end{proof}

\section{On the cycle structure of \texorpdfstring{$\chi^{(m)}$}{χ(m)}}\label{sec:cycle_structure}

Write $0_r$ for a string consisting of $r$ consecutive $0$s. In this chapter we aim to study the cycle structure of $\chi^{(m)}=\id+0_{m-1}1$. For $m=2$ we have $\chi^{(m)}=\chi^{(2)}=\chi$. Previous descriptions of the cycle structure of $\chi$ appear in \cite{kriepke2024algebraic,schoone2024state} and for $\chi^{(m)}$ in \cite{lyu2026generalized}. Our contribution here is that we give the exact count of the number of cycles of every length in the form of a generating function. 

In general, determining the exact cycle structure of a permutation on a finite field or finite vector space is hard and it has been done only for relatively small number of families. Some important examples include monomials \cite{ahmad1969monomial}, Dickson polynomials \cite{rudolf1991dickson}, linear polynomials \cite{mullen1988linear}, Rédei functions \cite{qureshi2015redei}, functions of the form $x\mapsto x+\gamma\Tr(h(x))$ \cite{gerike2020determining} and generalized cyclotomic polynomials \cite{bors2025cyclotomic}. We also refer to the survey \cite{martins2019survey} for more information.

We start by determining the order of $\chi^{(m)}$. This has been done before in \cite{kriepke2024algebraic,schoone2024state} for $m=2$ and \cite{lyu2026generalized} for general $m$.

\begin{corollary}\label{cor:cycle_structure_order_of_chi_n_m}
    Let $m\nmid n$. Then $\chi^{(m)}$ has order $2^\rho$, where $2^{\rho-1}<\ceil{\frac{n}{m}}\leq 2^{\rho}$.
\end{corollary}

\begin{proof}
    By \cref{lem:pcp_min_poly_if_full_support} the minimal polynomial of $*0_{m-1}1$ is $m_{*0_{m-1}1,n}(X)=X^{\ceil{n/m}}$. The claim follows by applying \cref{lem:pcp_order_of_F_in_R}.
\end{proof}

It follows from \cref{cor:cycle_structure_order_of_chi_n_m} that every cycle in the cycle decomposition of $\chi^{(m)}$ has length equal to a power of $2$. Observe that $x\in\F_2^n$ is contained in a cycle of length dividing $t$ of a function $f$ if and only if $f^t(x)=x$, i.e. $x$ is a fixed point of $f^t$. In that sense the study of the cycle structure of $f$ is equivalent to finding the fixed points of the iterates $f^t$. \cref{lem:pcp_formula_for_iterates} enables us to compute $f^t$ for binomials $f\in\mG$. This gives us a complete characterization of the elements $x\in \F_2^n$ occurring in cycles of a given length.

\begin{lemma}\label{lem:cycle_structure_in_cycle_iff_gamma_j_x_equals_0}
    A vector $x\in \F_2^n$ belongs to a cycle of length dividing $2^\mu$ of $\chi^{(m)}$ if and only if $\gamma_{2^\mu}(x)=0$.
\end{lemma}

\begin{proof}
    By \cref{lem:pcp_formula_for_iterates},
    \begin{equation*}
        (\id+c)^{2^\mu} = \id+\gamma_{2^\mu}.
    \end{equation*}
    Therefore, $(\id+c)^{2^\mu}(x)=x$ if and only if $\gamma_{2^\mu}(x)=0$. 
\end{proof}

Motivated by \cref{lem:cycle_structure_in_cycle_iff_gamma_j_x_equals_0} we now want to compute more generally the number of $x\in\F_2^n$ with $\gamma_j(x)=0$ for every $j\geq 0$. For $j=0$ the answer is trivial, since $\gamma_0=\id$ is the identity.

\begin{lemma}[\cite{kriepke2024algebraic,lyu2026generalized}]\label{lem:cycle_structure_gamma_j_of_x_zero_iff_no_substring}
    Let $j\geq 1$ and $x\in\F_2^n$. Then $\gamma_j(x)=0$ if and only if $x$ does not contain a substring of the form $0_{m-1}\mt0_{m-1}\mt\ldots0_{m-1}\mt0_{m-1}1$ of length $jm$. More precisely, there exists no index $i$ such that
    \begin{equation*}
        x_i x_{i+1} x_{i+2} \ldots x_{i+jm-1} = 0_{m-1} b_1 0_{m-1} b_2 \ldots 0_{m-1} b_{j-1} 0_{m-1} 1
    \end{equation*}
    for some $b_1,b_2,\ldots,b_{j-1}\in\F_2$. As usual, the indices are computed modulo $n$. 
\end{lemma}

\begin{proof}
    Follows immediately from the fact that 
    \begin{equation*}
        \gamma_j = \cq \odot S^q\cq \odot\cdots\odot S^{(j-2)q}\cq \odot S^{(j-1)q}c
    \end{equation*}
    is given by the landscape $*0_{m-1}\mt0_{m-1}\mt\ldots0_{m-1}\mt0_{m-1}1$.
\end{proof}

By \cref{lem:cycle_structure_gamma_j_of_x_zero_iff_no_substring}, the study of the cycle structure of $\chi^{(m)}$ reduces to a question in combinatorics on words. We introduce the required terminology and results in the next section.

\subsection{Generating functions and combinatorics on words}

In this section we discuss some results about generating functions and combinatorics on words. The discussion of generating functions will be brief, for more details we refer to \cite{flajoletAnalyticCombinatorics2009,Stanley2011enumerativeVol1,wilf2006generatingfunctionology}. 

Given a sequence $(z_n)=(z_n)_{n\in\N_0}$ of non-negative integers, we define the generating function
\begin{equation*}
    \mZ(\lambda) = \sum_{n=0}^\infty z_n \lambda^n.
\end{equation*}
The expression $\mZ(\lambda)$ can either be viewed as a formal power series in $\Z[[\lambda]]$ or as a complex-valued function for those values of $\lambda\in\C$ for which it converges. The theory of generating functions allows deriving many results about the properties of the sequence $(z_n)$. A special class of generating functions are rational functions, i.e. generating functions of the form
\begin{equation*}
    \mZ(\lambda) = \frac{P(\lambda)}{Q(\lambda)},
\end{equation*}
where $P(\lambda),Q(\lambda)$ are polynomials. Rational generating functions occur frequently in combinatorics and are well understood. The following theorem is well-known, see \cite[Theorem 4.1.1]{Stanley2011enumerativeVol1}. 

\begin{theorem}\label{thm:cycle_structure_rational_gf_properties}
    Let $(z_n)$ be a sequence of numbers and suppose
    \begin{equation*}
        Q(\lambda)= \sum_{i=0}^d \alpha_i \lambda^i= \prod_{i=1}^k (\lambda-\xi_i)^{d_i}
    \end{equation*}
    is a polynomial, where the $\xi_i$'s are distinct and nonzero. Then the following are equivalent:
    \begin{enumerate}
        \item The generating function of $(z_n)$ is rational of the form
        \begin{equation*}
            \sum_{n\geq 0} z_n \lambda^n = \frac{P(\lambda)}{Q(\lambda)},
        \end{equation*}
        for some polynomial $P(\lambda)$ of degree less than $d$.
        \item The sequence $(z_n)$ satisfies the recurrence
        \begin{equation*}
            \alpha_0 z_{n+d} + \alpha_1 z_{n+d-1} + \cdots + \alpha_d z_n = 0
        \end{equation*}
        for all $n\geq 0$.
        \item For all $n\geq 0$ it holds that
        \begin{equation*}
            z_n = \sum_{i=1}^k P_i(n)\xi_i^{-n}
        \end{equation*}
        for some polynomials $P_i(n)$ of degree less than $d_i$. 
    \end{enumerate}
\end{theorem}

\cref{thm:cycle_structure_rational_gf_properties} $(c)$ gives a way to determine the asymptotic behavior of the sequence $(z_n)$ for $n\to\infty$. For simplicity suppose that $\abs{\xi_1} < \abs{\xi_i}$ for $i\geq 2$ and that $\xi_1$ is a simple root of $Q(\lambda)$. Then $z_n\sim c\,\xi_1^{-n}$ for some constant $c$, which can be determined explicitly. This is a standard consequence of partial fractions, see for example \cite{flajoletAnalyticCombinatorics2009}. 

\begin{lemma}\label{lem:cycle_structure_explicit_coefficient_of_asymptotic_rational_function}
    With the notation from above,
    \begin{equation*}
        z_n \sim -\frac{P(\xi_1)}{\xi_1 Q'(\xi_1)}\xi_1^{-n}.
    \end{equation*}
\end{lemma}

\begin{proof}
    The partial fraction decomposition of $P(\lambda)/Q(\lambda)$ is given by 
    \begin{equation*}
        \frac{P(\lambda)}{Q(\lambda)} = \frac{P(\xi_1)}{Q'(\xi_1)}\frac{1}{\lambda-\xi_1} + \sum_{i=2}^k \sum_{j=1}^{d_i} \frac{c_{i,j}}{(\lambda-\xi_i)^j}
    \end{equation*}
    for some complex numbers $c_{i,j}$. Focusing just on the first term, 
    \begin{equation*}
        \frac{P(\xi_1)}{Q'(\xi_1)} \frac{1}{\lambda-\xi_1}
            = \frac{-P(\xi_1)}{\xi_1 Q'(\xi_1)} \frac{1}{1-\xi_1^{-1}\lambda}
            = \sum_{n=0}^\infty \frac{-P(\xi_1)}{\xi_1 Q'(\xi_1)} \xi_1^{-n} \lambda^n.
    \end{equation*}
    The claim follows from comparing the coefficients and using that $\xi_1$ is the root of smallest absolute value, which dominates the asymptotic growth. 
\end{proof}

Next we discuss some results about combinatorics on words. We closely follow \cite[Section 4.7.4]{Stanley2011enumerativeVol1}. For the reader's convenience we also include proofs of the results. 

Let $\mA$ be a finite set, called the alphabet. We write a tuple $a=(a_1,a_2,\ldots,a_n)\in\mA^n$ as a string $a_1\ldots a_n$ and call it a word. The length $\ell(a)$ of $a$ is $n$. We denote the empty word by $\varepsilon$, it has length $0$. The set $\mA^*$ is the set of all words in the alphabet $\mA$, i.e.
\begin{equation*}
    \mA^* = \bigcup_{n\geq 0} \mA^n,
\end{equation*}
where $\mA^0\coloneqq \{\varepsilon\}$. On $\mA^*$ we consider the operation $\concat:\mA^*\times \mA^*\to \mA^*$ obtained from concatenating two words $a_1\ldots a_n, b_1\ldots b_m\in \mA^*$, i.e. 
\begin{equation*}
    a_1\ldots a_n \concat b_1\ldots b_m = a_1\ldots a_nb_1\ldots b_m.
\end{equation*}
It is easy to see that this operation is associative and that $\varepsilon\concat u=u\concat \varepsilon=u$ for all words $u\in \mA^*$, so that $(\mA^*, \concat)$ is a monoid with identity $\varepsilon$. 

\begin{example}
    \begin{enumerate}
        \item Let $\mA=\{1\}$. Then $\mA^* = \{\varepsilon, 1, 11, 111, 1111, \ldots\}$ is a monoid isomorphic to $(\N_0,+)$.
        \item Let $\mA=\F_2$. Then $\mA^*=\F_2^*$ contains for example the words $0, 111, 01010$ and $000011$. 
    \end{enumerate}
\end{example}

For a subset $\mB\subseteq \mA^*$ we consider the submonoid $\mB^*$ generated by $\mB$, i.e. the set of all finite products of elements from $\mB$. Observe that by choosing $\mB=\mA$, we obtain $\mB^*=\mA^*$ again. The submonoid $\mB^*$ is called freely generated by $\mB$ if every word $u\in\mB^*$ has a unique factorization into elements of $\mB$. More precisely, if $u=u_1\concat\cdots \concat u_n = v_1\concat\cdots \concat v_m$ with $u_i,v_i\in\mB$, then $n=m$ and $u_i=v_i$ for all $i=1,\ldots,n$. 

\begin{example}
    \begin{enumerate}
        \item Let $\mA=\{1\}$. If $\mB$ contains at least $2$ elements, then $\mB^*$ is not freely generated. Suppose for example that $11,111\in\mB$, then 
        \begin{equation*}
            111111=11\concat11\concat11=111\concat111
        \end{equation*}
        cannot be uniquely factored. If $\mB$ has exactly one element, then $\mB^*$ is freely generated.
        \item Let $\mA=\F_2$ and $\mB=\{1, 101, 001\}\subseteq \mA^*$. Then $\mB^*$ is freely generated by $\mB$. Indeed, take for example the word $10011001101100110111\in\mB^*$. Then from the positions of the zeros we immediately read off that the unique factorization is given by 
        \begin{equation*}
            10011001101100110111 = 1\concat 001\concat 1\concat 001\concat 101\concat 1\concat 001\concat 101\concat 1\concat 1.
        \end{equation*}
    \end{enumerate}
\end{example}

For any subset $\mC\subseteq \mA^*$, define the generating function
\begin{equation*}
    \mC(\lambda) = \sum_{x\in\mC} \lambda^{\ell(x)} = \sum_{n\geq 0} \abs{\mC_n}\lambda^n,
\end{equation*}
where $\mC_n = \mC\cap \mA^n$ is the set of all words of length $n$ from $\mC$. If $\mB^*$ is freely generated by $\mB$, the generating function $B^*(\lambda)$ can be computed using the generating function $\mB(\lambda)$.

\begin{lemma}\label{lem:cycle_structure_comb_on_words_freely_generated_monoid_gf}
    Let $\mB^*\subseteq \mA^*$ be freely generated by $\mB\subseteq \mA^*$. Then
    \begin{equation*}
        \mB^*(\lambda) = \sum_{n\geq 0} \mB(\lambda)^n = \frac{1}{1-\mB(\lambda)}.
    \end{equation*}
\end{lemma}

\begin{proof}
    The generating function $\mB(\lambda)^n$ counts the number of words in $\mB^*$ whose unique factorization into elements of $\mB$ contains exactly $n$ factors. As $\mB^*$ is freely generated by $\mB$, every word in $\mB^*$ is counted exactly once by exactly one of the generating functions $\mB(\lambda)^n$ for some $n\in\N_0$, thus their sum is the generating function of $\mB^*$. The last equality is the geometric series.
\end{proof}

As mentioned earlier, we are interested in counting nonzero elements $x\in\F_2^n$ such that $\gamma_j(x)=0$, for some fixed landscape $c=*0_{m-1}1$ and fixed $j, n$. Rephrasing this using the terminology we introduced, this means we are interested in counting all nonzero words $x$ of length $n$ from $\F_2^*$ with $x\in\mZ^{(m)}_j$, where
\begin{equation*}
    \mZ^{(m)}_j = \{x\in\F_2^*\setminus\{0\}: \gamma_j(x)=0\}
\end{equation*}
and $\gamma_j$ is interpreted as a function on $\F_2^{\ell(x)}$. We cannot, however, apply \cref{lem:cycle_structure_comb_on_words_freely_generated_monoid_gf}, as $\mZ^{(m)}_j$ is not a monoid. Indeed, take for example $c=*01$ and $j=2$, so that $\gamma_2 = *0\mt01$ and $\gamma_2(x)=0$ if and only if $x$ does not contain the pattern $0\mt01$. We easily see that $x=1101$ and $y=011011$ belong to $\mZ_2^{(2)}$, but 
\begin{equation*}
    x\concat y=11101\concat011011=111\underline{0101}1011\notin\mZ_2^{(2)}.
\end{equation*}
On the other hand, we have not used that $\gamma_j$ is shift-invariant, implying that $\mZ^{(m)}_j$ is closed under $S$. If we replace $y=011011$ by $Sy=110110$ then
\begin{equation*}
    x\concat Sy=11101110110\in\mZ_2^{(2)}. 
\end{equation*}
In fact, $\mZ^{(m)}_j$ is the set of shifts of a certain monoid. We use this observation to compute $\mZ_j^{(m)}(\lambda)$. We first need to introduce a definition.

\begin{definition}\label{def:unique_circular_factorization}
    Suppose that $\mB^*$ is freely generated by $\mB$, i.e. every word in $\mB^*$ can be uniquely factored into words from $\mB$. The monoid $\mB^*$ is called very pure if every word in $\mB^*$ has a unique circular factorization into words from $\mB$. More precisely, let $x=a_1a_2\ldots a_n\in\mB^*$ with letters $a_i\in\mA$. Then there are unique integers $0=n_0<n_1<n_2<\ldots<n_k<n$ such that $b_j=a_{n_j+1}\ldots a_{n_{j+1}}\in\mB$, i.e. $x=b_0\concat \cdots \concat b_k$ is the unique factorization of $x$ into words from $\mB$. Suppose that there is some shift $S^ix=a_i a_{i+1}\ldots a_n a_1 \ldots a_{i-1}\in \mB^*$. Then it has to hold that $i=n_j+1$ for some $0\leq j\leq k$.

    In other words, suppose we write $x$ on a circle, without specifying the beginning position $a_1$. Then there is a unique way to insert concatenation symbols $\concat$ between letters of $x$, so that the words formed by reading these letters is a word in $\mB$. 

\end{definition}

\begin{figure}[!ht]
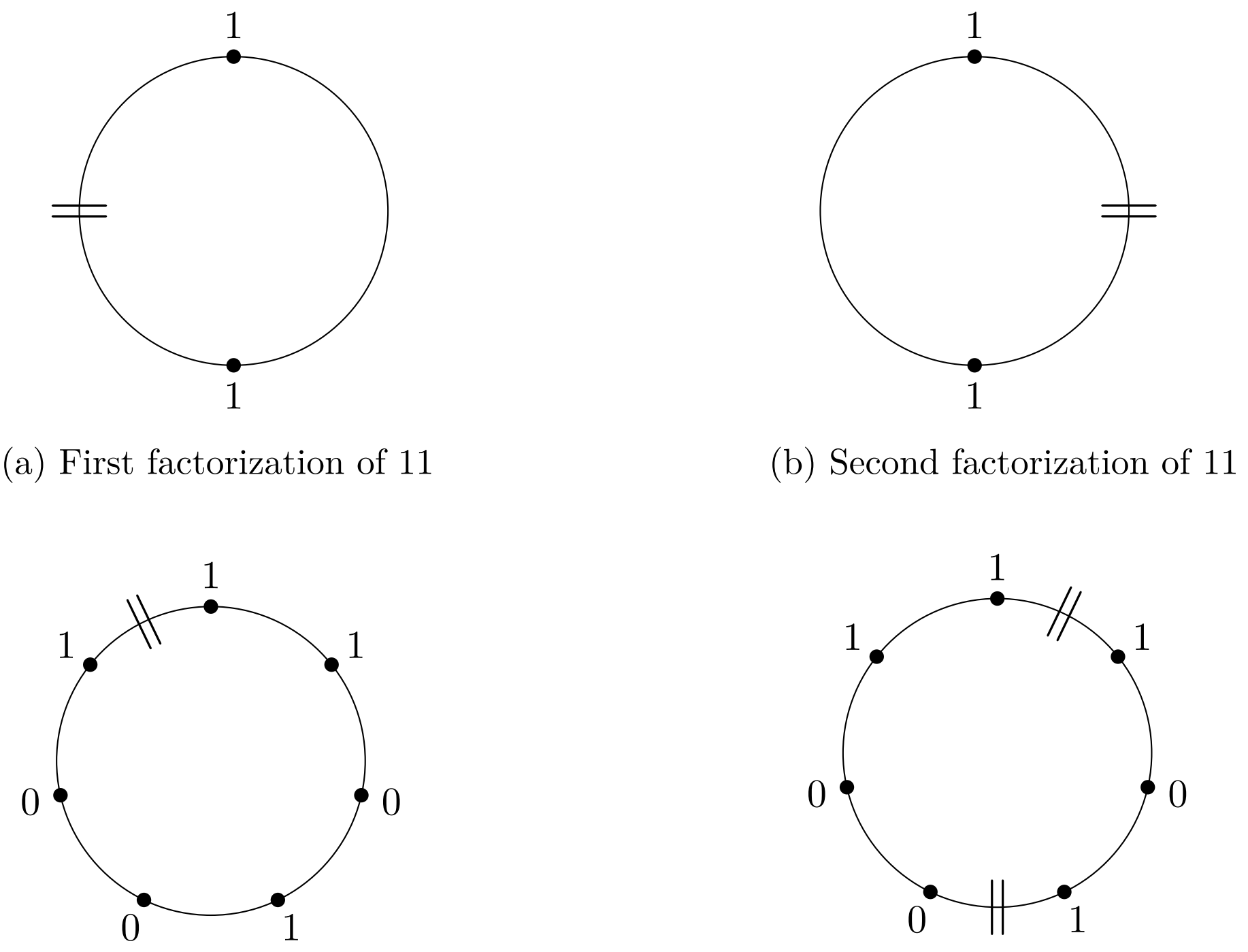

\centering

\begin{subfigure}{0.48\textwidth}
\centering
\wordcircle{1.5}{2}{1,1}{2}
\caption{First factorization of $11$}
\end{subfigure}
\hfill
\begin{subfigure}{0.48\textwidth}
\centering
\wordcircle{1.5}{2}{1,1}{1}
\caption{Second factorization of $11$}
\end{subfigure}

\vspace{0.5cm}
\begin{subfigure}{0.48\textwidth}
\centering
\wordcircle{1.5}{7}{1,1,0,1,0,0,1}{7}
\caption{First factorization of $1101001$}
\end{subfigure}
\hfill
\begin{subfigure}{0.48\textwidth}
\centering
\wordcircle{1.5}{7}{1,1,0,1,0,0,1}{1,4}
\caption{Second factorization of $1101001$}
\end{subfigure}
\caption{The factorizations of the words from \cref{ex:cycle_structure_failure_unique_circular_factorization}\label{fig:cycle_structure_failure_unique_circular_factorization}}
\end{figure}

\begin{example}\label{ex:cycle_structure_failure_unique_circular_factorization}
    \begin{enumerate}
        \item Let $\mA=\{1\}$ and $\mB=\{11\}$. Then $\mB^*$ is freely generated by $\mB$, but it is not very pure. Indeed, the word $11=a_1a_2\in \mB^*$ has two circular factorizations, namely $11\concat$ and $1\concat1$.
        \item Let $\mA=\F_2$ and $\mB=\{101,0011,1101001\}$. Then $1101001\in\mB^*$ has the two circular factorization $1101001\concat $ and $1\concat101\concat001$.
    \end{enumerate}
    The different circular factorizations of these words are depicted in \cref{fig:cycle_structure_failure_unique_circular_factorization}
\end{example}

Suppose now that $\mB^*$ is very pure. Define the set
\begin{equation*}
    \widetilde{\mB} = \{x\in\mA^*\setminus\{\varepsilon\}: x \text{ is a shift of a word in } \mB^*\}.
\end{equation*}
Our goal is to count the elements of $\widetilde{\mB}$. We have the following result.

\begin{lemma}\label{lem:cycle_structure_comb_on_words_very_pure_gf}
    Let $\mB^*$ be very pure. Then
    \begin{equation*}
        \widetilde{\mB}(\lambda) = \frac{\lambda \frac{\d}{\d\lambda}\mB(\lambda)}{1-\mB(\lambda)}.
    \end{equation*}
\end{lemma}

\begin{proof}
    Let $x\in\widetilde{\mB}$. Then there is a word $b=b_1\ldots b_m\in\mB$, an index $i\in\{1,\ldots,m\}$ and a word $y\in\mB^*$ such that $x=b_i\ldots b_m y b_1\ldots b_{i-1}$. Moreover, $b,i,y$ are unique as $\mB^*$ is very pure. We say that $x$ is a shift that splits up $b$.
    
    Now flip the perspective, fix $b\in\mB$ and count the number of shifts of length $n$ that split up $b$. Suppose the length of $b$ is $\ell(b)=m$, then there are $m$ positions at which $b$ can be split up and there are $\abs{\mB_{n-m}^*}$ options to choose the word $y\in\mB^*$ of length $n-m$. Thus,
    \begin{equation*}
        \abs{\widetilde{\mB}_n} = \sum_{b\in \mB} \ell(b) \abs{\mB_{n-\ell(b)}^*} = \sum_{m=0}^n m \abs{\mB_m} \abs{\mB_{n-m}^*}
    \end{equation*}
    and therefore
    \begin{align*}
        \widetilde{\mB}(\lambda) 
            &= \sum_{n\geq 0} \abs{\widetilde{\mB}_n} \lambda^n 
             = \sum_{n\geq 0} \sum_{m=0}^n m \abs{\mB_m} \abs{\mB_{n-m}^*} \lambda^n \\
            &= \sum_{m\geq 0} \sum_{n\geq m} m \abs{\mB_m} \abs{\mB_{n-m}^*} \lambda^m \lambda^{n-m} \\
            &= \left(\sum_{m\geq 0} m\abs{\mB_m} \lambda^m \right)\left(\sum_{n\geq m} \abs{\mB_{n-m}^*} \lambda^{n-m}\right) \\
            &= \lambda \left(\frac{\d}{\d\lambda} \mB(\lambda)\right) \mB^*(\lambda) 
             = \frac{\lambda \frac{\d}{\d\lambda}\mB(\lambda)}{1-\mB(\lambda)},
    \end{align*}
    where the last equality follows from \cref{lem:cycle_structure_comb_on_words_freely_generated_monoid_gf}.
\end{proof}

\subsection{Cycle structure of \texorpdfstring{$\chi^{(m)}$}{χ(m)}}

We now want to apply the results of the previous section to $\chi^{(m)}$, i.e. we want to count the number of nonzero $x$ of length $n$ in
\begin{equation*}
    \mZ_j^{(m)} = \{x\in\F_2^*: \gamma_j(x)=0\}.
\end{equation*}
We start by setting 
\begin{align*}
    \mZ_{j,n}^{(m)} &= \{x\in\F_2^n\setminus\{0\}: \gamma_j(x)=0\} \\
        &= \{x\in \F_2^n\setminus\{0\}: x \text{ does not contain } 0_{m-1}\mt0_{m-1}\mt\ldots0_{m-1}\mt0_{m-1}1 \text{ cyclically} \}
\end{align*}
and define $z^{(m)}(j,n) = \abs{\mZ_{j,n}^{(m)}}$, so that
\begin{equation*}
    \mZ^{(m)}_j(\lambda) = \sum_{n\geq 0} z^{(m)}(j,n)\lambda^n.
\end{equation*}

The following theorem gives a complete description of the cycle structure of $\chi^{(m)}$.

\begin{theorem}\label{thm:cycle_structure_of_chi_m}
    Let $m\nmid n$. The cycles in the cycle decomposition of $\chi^{(m)}$ have lengths $2^\mu$ for $0\leq \mu\leq \rho$, where $\rho$ is given by $2^{\rho-1} < \ceil{\frac{n}{m}}\leq 2^\rho$. There are exactly $z^{(m)}(2^\mu,n)$ nonzero vectors $x\in\F_2^n$ which belong to cycles whose lengths divide $2^\mu$. The number of cycles of length $1$ is $z^{(m)}(1,n)+1$. For $\mu\geq 1$, the number of cycles of length $2^\mu$ is given by 
    \begin{equation*}
        \frac{z^{(m)}(2^\mu,n)-z^{(m)}(2^{\mu-1},n)}{2^\mu}.
    \end{equation*}
\end{theorem}

The rest of the section is dedicated to studying the numbers $z^{(m)}(j,n)$. Our main goal is to obtain a closed form for their generating function
\begin{equation*}
    \mZ^{(m)}_j(\lambda) = \sum_{n\geq 0} z^{(m)}(j,n)\lambda^n.
\end{equation*}
\cref{tab:values_z_2_j_n,tab:values_z_3_j_n} show the values of $z^{(m)}(j,n)$ for $m\in{2,3}$ and small values of $j,n$. Before turning to the generating function, we record several elementary properties of the numbers $z^{(m)}(j,n)$. 

\begin{table}
    \begin{subtable}{\textwidth}
        \centering
        \begin{tabular}{r|rrrrrrrr}
            $n\backslash j$ & $1$ & $2$ & $3$ & $4$ & $5$ & $6$ & $7$ & $8$ \\ \hline
            $1$ & $1$ & $1$ & $1$ & $1$ & $1$ & $1$ & $1$ & $1$ \\
            $2$ & $1$ & $1$ & $1$ & $1$ & $1$ & $1$ & $1$ & $1$ \\
            $3$ & $1$ & $7$ & $7$ & $7$ & $7$ & $7$ & $7$ & $7$ \\
            $4$ & $1$ & $9$ & $9$ & $9$ & $9$ & $9$ & $9$ & $9$ \\
            $5$ & $1$ & $11$ & $31$ & $31$ & $31$ & $31$ & $31$ & $31$ \\
            $6$ & $1$ & $25$ & $49$ & $49$ & $49$ & $49$ & $49$ & $49$ \\
            $7$ & $1$ & $43$ & $71$ & $127$ & $127$ & $127$ & $127$ & $127$ \\
            $8$ & $1$ & $65$ & $161$ & $225$ & $225$ & $225$ & $225$ & $225$ \\
            $9$ & $1$ & $115$ & $295$ & $367$ & $511$ & $511$ & $511$ & $511$ \\
            $10$ & $1$ & $201$ & $561$ & $801$ & $961$ & $961$ & $961$ & $961$ \\
            $11$ & $1$ & $331$ & $1079$ & $1519$ & $1695$ & $2047$ & $2047$ & $2047$ \\
            $12$ & $1$ & $561$ & $1953$ & $3009$ & $3585$ & $3969$ & $3969$ & $3969$ \\
            $13$ & $1$ & $963$ & $3719$ & $5903$ & $6943$ & $7359$ & $8191$ & $8191$ \\
            $14$ & $1$ & $1625$ & $7057$ & $11425$ & $13889$ & $15233$ & $16129$ & $16129$ \\
            $15$ & $1$ & $2747$ & $13207$ & $22447$ & $27487$ & $29887$ & $30847$ & $32767$ \\
        \end{tabular}
        \caption{Values $z^{(2)}(j,n)$ for $1\leq j\leq 8$ and $1\leq n\leq 15$.}\label{tab:values_z_2_j_n}
    \end{subtable}

    \vspace{0.5cm}
    \begin{subtable}{\textwidth}
        \centering
        \begin{tabular}{r|rrrrrrr}
            $n\backslash j$ & $1$ & $2$ & $3$ & $4$ & $5$ & $6$ & $7$ \\ \hline
            $1$ & $1$ & $1$ & $1$ & $1$ & $1$ & $1$ & $1$ \\
            $2$ & $3$ & $3$ & $3$ & $3$ & $3$ & $3$ & $3$ \\
            $3$ & $4$ & $4$ & $4$ & $4$ & $4$ & $4$ & $4$ \\
            $4$ & $7$ & $15$ & $15$ & $15$ & $15$ & $15$ & $15$ \\
            $5$ & $11$ & $31$ & $31$ & $31$ & $31$ & $31$ & $31$ \\
            $6$ & $18$ & $54$ & $54$ & $54$ & $54$ & $54$ & $54$ \\
            $7$ & $29$ & $99$ & $127$ & $127$ & $127$ & $127$ & $127$ \\
            $8$ & $47$ & $191$ & $255$ & $255$ & $255$ & $255$ & $255$ \\
            $9$ & $76$ & $382$ & $490$ & $490$ & $490$ & $490$ & $490$ \\
            $10$ & $123$ & $743$ & $943$ & $1023$ & $1023$ & $1023$ & $1023$ \\
            $11$ & $199$ & $1431$ & $1871$ & $2047$ & $2047$ & $2047$ & $2047$ \\
            $12$ & $322$ & $2754$ & $3762$ & $4050$ & $4050$ & $4050$ & $4050$ \\
            $13$ & $521$ & $5331$ & $7463$ & $7983$ & $8191$ & $8191$ & $8191$ \\
            $14$ & $843$ & $10335$ & $14815$ & $15935$ & $16383$ & $16383$ & $16383$ \\
            $15$ & $1364$ & $20014$ & $29434$ & $31954$ & $32674$ & $32674$ & $32674$ \\
            $16$ & $2207$ & $38719$ & $58495$ & $63743$ & $65023$ & $65535$ & $65535$ \\
            $17$ & $3571$ & $74903$ & $116111$ & $127263$ & $129983$ & $131071$ & $131071$ \\
            $18$ & $5778$ & $144954$ & $230418$ & $254178$ & $260226$ & $261954$ & $261954$ \\
            $19$ & $9349$ & $280555$ & $457559$ & $507567$ & $520031$ & $523071$ & $524287$ \\
            $20$ & $15127$ & $542975$ & $908735$ & $1013375$ & $1039615$ & $1046015$ & $1048575$ \\
            $21$ & $24476$ & $1050774$ & $1804618$ & $2023186$ & $2078626$ & $2092738$ & $2096770$ \\
        \end{tabular}
        \caption{Values $z^{(3)}(j,n)$ for $1\leq j\leq 7$ and $1\leq n\leq 21$.}\label{tab:values_z_3_j_n}
    \end{subtable}
\end{table}

\begin{lemma}\label{lem:cycle_structure_properties_numbers_z_m_j_n}
    The following statements hold:
    \begin{enumerate}
        \item If $j=1$ and $m=2$ then $z^{(2)}(1,n)=1$ for all $n$.
        \item\label{itm:cycle_structure_z_j_n_exponential_in_n} If $j\geq2$ or $m\geq 3$ then there exists a constant $\zeta_{j,m}>1$ such that $\zeta_{j,m}^n < z^{(m)}(j,n)<2^n$ for all $n\geq 3$.
        \item\label{itm:z_j_n_monotone_in_j} The map $j\mapsto z^{(m)}(j,n)$ is monotone increasing for fixed $n,m$, in particular it is eventually constant.
        \item Let $m\nmid n$. Then $z^{(m)}(j,n)=2^n-1$ for $j\geq\ceil{\frac{n}{m}}$. 
        \item Let $m\mid n$. Then $z^{(m)}(j,n)$ is constant for $j\geq \frac n m$. 
        \item The map $n\mapsto z^{(m)}(j,n)$ is monotone increasing for fixed $j,m$.
    \end{enumerate}
\end{lemma}

\begin{proof}
    \begin{enumerate}
        \item We have $\mZ^{(2)}_{1,n} = \{x\in\F_2^n\setminus\{0\}: x\text{ does not contain } 01 \text{ cyclically} \} = \{\one\}$.
        \item Suppose that $n\equiv 0\bmod 3$. Set 
        \begin{equation*}
            \overline{S_{n,3}}=\{11b_1 11b_2\ldots 11b_{n/3}: b_i\in\F_2\}.
        \end{equation*}
        Then $\overline{S_{n,3}}\subseteq \mZ_j^{(m)}$ and the claim follows as $\abs{\overline{S_{n,3}}}=\left(2^{1/3}\right)^n$. The cases $n\equiv 1,2\bmod 3$ are similar. We note that the constant $\zeta_{j,m}=2^{1/3}$ is not optimal. 
        \item Suppose that $x\in\mZ_{j,n}^{(m)}$, i.e. $x$ does not contain $0_{m-1}\mt0_{m-1}\mt\ldots0_{m-1}\mt0_{m-1}1$ of length $jm$. Then it also does not contain $0_{m-1}\mt0_{m-1}\mt\ldots0_{m-1}\mt0_{m-1}1$ of length $(j+1)m$ and hence $x\in\mZ_{j+1,n}^{(m)}$. 
        \item By \cref{lem:pcp_min_poly_if_full_support} we have $\gamma_{\ceil{n/m}}=\zero$, i.e. $\gamma_{\ceil{n/m}}(x)=0$ for all $x\in\F_2^n$. The claim follows (recall that by definition $0\notin \mZ_{j,n}^{(m)}$). 
        \item By \cref{lem:pcp_min_poly_if_full_support} we have $\gamma_{n/m} = \gamma_{2n/m}$ which implies $z^{(m)}(\frac n m, n) = z^{(m)}(\frac{2n}{m},n)$. The claim follows together with \cref{itm:z_j_n_monotone_in_j}.
        \item Take $x\in\mZ_{j,n}^{(m)}$. As $x\neq 0$ there exists a minimal index $i\in\{1,\ldots,n\}$ with $x_i=1$. Define $\hat x\in\F_2^{n+1}$ by inserting a $1$ at position $i+1$, i.e. $\hat x_j=x_j$ for $j\leq i$, $\hat x_{i+1}=1$ and $\hat x_j=x_{j-1}$ for $i+1\leq j\leq n+1$. As $x$ does not contain a pattern $0_{m-1}\mt0_{m-1}\mt\ldots0_{m-1}\mt0_{m-1}1$ it is clear that $\hat x$ also does not contain this pattern. The map $x\mapsto \hat x$ gives an injective map $\mZ_{j,n}^{(m)}\to \mZ_{j,n+1}^{(m)}$, implying the claim.
    \end{enumerate}
\end{proof}

We remark that several of these properties could also be proved after finding the generating function $\mZ_j^{(m)}(\lambda)$. In particular, later we obtain a precise asymptotic of $z^{(m)}(j,n)$ for fixed $j$ and $n\to\infty$ using \cref{thm:cycle_structure_rational_gf_properties} and \cref{lem:cycle_structure_explicit_coefficient_of_asymptotic_rational_function}. We think it is better to give elementary proofs here as they might generalize more easily to other landscapes.

We now claim that $\mZ^{(m)}_j = \widetilde{\mB}$ for some subset $\mB\subseteq \F_2^*$ with the property that $\mB^*$ is very pure. To show this we introduce the following definition, which is a generalization of the concept of anchors from \cite[Definition 10]{schoone2024state}.

\begin{definition}\label{def:anchor}
    Let $m\geq 2$ and fix $x\in\F_2^n$. A position $i\in\{1,\ldots,n\}$ is called an anchor of $x$ if $x_i=1$ and $x_{i+d}=1$ for minimal $d>0$ implies that $m\nmid d$. As usual, the indices are computed modulo $n$. 

    The anchor factorization of $x$ is obtained by inserting the concatenation symbol $\concat$ after every anchor.
\end{definition}

\begin{example}\label{ex:cycle_structure_anchor_factorization}
    \begin{enumerate}
        \item Let $m=2$. Then $i$ is an anchor if $x_i=1$ and the distance to the next $1$ is odd. The following shows the anchors highlighted as $\hat1$ and the resulting anchor factorization for some $x\in\F_2^{47}$:
        \begin{align*}
            x &= 01001100100101010100101000110111001101001010111 \\
            &= 0\hat100\hat1\hat100\hat100101010\hat100101000\hat110\hat1\hat1\hat100\hat110\hat1001010\hat1\hat11 \\
            &= 0\hat1\concat00\hat1\concat\hat1\concat00\hat1\concat00101010\hat1\concat00101000\hat1\concat10\hat1\concat\hat1\concat\hat1\concat00\hat1\concat10\hat1\concat001010\hat1\concat\hat1\concat1.
        \end{align*}
        \item Let $m=3$. Then the anchor factorization of the same $x$ is given by 
        \begin{align*}
            x &= 01001100100101010100101000110111001101001010111 \\
            &= 0100\hat1100100\hat10\hat10\hat10100\hat10\hat1000\hat1\hat10\hat1\hat1100\hat1\hat10100\hat10\hat10\hat1\hat1\hat1 \\
            &= 0100\hat1\concat100100\hat1\concat0\hat1\concat0\hat1\concat0100\hat1\concat0\hat1\concat000\hat1\concat\hat1\concat0\hat1\concat\hat1\concat100\hat1\concat\hat1\concat0100\hat1\concat0\hat1\concat0\hat1\concat\hat1\concat\hat1\concat.
        \end{align*}
    \end{enumerate}
\end{example}

Of course, the first natural question is whether an anchor factorization always exists. This is false. Take $m=2$, then $x=010001$ has no anchors and therefore no anchor factorization. Nevertheless, we can easily determine the sets of elements $x$ which have an anchor factorization. This was done for $m=2$ in \cite{schoone2024state} and the statement can be generalized. For $m\mid n$ and $0\leq r\leq m-1$, define the subspaces $S_{n,r}\subseteq\F_2^n$ by 
\begin{equation*}
    S_{n,r} = \{x\in\F_2^n: \text{if } x_i=1 \text{ then } i\equiv r\bmod m\} = \{x\in\F_2^n: x_i = 0 \text{ for } i\not\equiv r\bmod n\}.
\end{equation*}

\begin{lemma}\label{lem:cycle_structure_existence_anchor_factorization}
    Let $m\geq 2$. 
    \begin{enumerate}
        \item Suppose $m\nmid n$. Then $x\in\F_2^n\setminus\{0\}$ always has an anchor factorization.
        \item Suppose $m\mid n$. Then $x\in\F_2^n\setminus\{0\}$ has an anchor factorization if and only if $x\notin \bigcup_{r=0}^{m-1} S_{n,r}$.
    \end{enumerate}
\end{lemma}

\begin{proof}
    \begin{enumerate}
        \item Suppose $m\nmid n$ and take $x\in\F_2^n\setminus\{0\}$. We only need to show that $x$ has at least one anchor. For each $i$ with $x_i=1$ let $d_i$ be the distance to the next $1$ in $x$, i.e. $x_{i+d_i}=1$ with minimal $d_i>0$. It is clear that the sum of all distances $d_i$ equals $n$. As $m\nmid n$ it follows that at least one of the distances is not divisible by $m$, meaning that $x$ has an anchor.
        \item Suppose $m\mid n$ and $x\in S_{n,r}$ for some $0\leq r\leq m-1$. Then it is clear that the distance between any $1$s in $x$ is always divisible by $m$ and $x$ does not have an anchor.
        
        If $x\notin S_{n,r}$ for any $0\leq r\leq m-1$, then there exist two indices $i<j$ with $x_i=x_j=1$ and $i\not\equiv j\bmod m$. Consider the distances from each $1$ between $x_i$ and $x_j$ to the next $1$. These distances add up to a number not divisible by $m$, so at least one of them is also not divisible by $m$ and $x$ has an anchor. 
    \end{enumerate}
\end{proof}

\cref{lem:cycle_structure_existence_anchor_factorization} shows that for any $x\in\mZ_{j}^{(m)}$ there exists an anchor factorization. Indeed, if $x\in\bigcup_{r=0}^{m-1} S_{n,r}$ is nonzero, then $x$ contains at least one substring of the form $0_{m-1}\mt0_{m-1}\mt\ldots0_{m-1}\mt0_{m-1}1$ of length $jm$ for some $j$ and $\gamma_j(x)\neq 0$, so that $x\notin\mZ^{(m)}_j$. 

\begin{remark}
    The previous argument also shows that 
    \begin{equation*}
        z^{(m)}(j,n) = 2^n - \abs{\bigcup_{r=0}^{m-1} S_{n,r}} = 2^n - m2^{n/m}+m-1
    \end{equation*}
    for $m\mid n$ and $j\geq n/m$, as every vector $x\in\F_2^n$ is either counted by $z^{(m)}(j,n)$ for some large enough $j$ or has no anchor factorization. The second equality follows because the subspaces $S_{n,r}$ are in direct sum. Combined with \cref{lem:cycle_structure_properties_numbers_z_m_j_n} we now know the exact value of $z^{(m)}(j,n)$ for $n$ fixed and every $j\geq n/m$. 
\end{remark}

From here it suffices to determine the set $\mB_j^{(m)}$ such that $\mZ^{(m)}_j=\widetilde{\mB_j}^{(m)}$. By construction the anchor factorization gives a unique circular factorization of elements in $\mZ^{(m)}_j$. We can now read off the possible factors to obtain the set $\mB_j^{(m)}$. As the anchor factorization is obtained by placing the concatenation symbol after every anchor, every element in $\mB_j^{(m)}$ has to end in a $1$. Next, if an element of $\mB_j^{(m)}$ contains other $1$s, then the distances between those $1$s has to be divisible by $m$, as otherwise one of them would be an anchor. Furthermore, an element of $\mB_j^{(m)}$ cannot start with the substring $0_{m-1}1$, as otherwise in the anchor factorization of $x$ we would have $x_i\concat0_{m-1}1$ and $x_i=1$ because $i$ is an anchor, but then the distance from $i$ to the next $1$ is $m$ which is divisible by $m$. A similar argument excludes elements that start with $0_m0_m\ldots 0_m0_{m-1}1$. 

We obtain
\begin{equation*}
    \mB_j^{(m)} = \{0_r (b_10_{m-1})(b_2 0_{m-1}) \ldots (b_i 0_{m-1})1: 0\leq r\leq m-2, b_1,\ldots,b_i \in\F_2, 0\leq i<j\}, 
\end{equation*}
where we added some parenthesis for clarity. In other words, every element in $\mB_j^{(m)}$ contains a block $0_r$, then $i$ blocks of the form $(b_k0_{m-1})$ and then a $1$. 

\begin{example}
    \begin{enumerate}
        \item If $m=2$ then 
        \begin{equation*}
            \mB_j^{(2)} = \{(b_1 0)(b_2 0) \ldots (b_i 0)1: b_1,\ldots,b_i\in\F_2, 0\leq i<j\},
        \end{equation*}
        so that 
        \begin{align*}
            \mB_1^{(2)} &= \{1\}, \\ 
            \mB_2^{(2)} &= \{1, 001, 101\}, \\
            \mB_3^{(2)} &= \{1, 001, 101, 00001, 00101, 10001, 10101\}, \\
            \mB_4^{(2)} &= \{1, 001, 101, 00001, 00101, 10001, 10101, 0000001, 0000101, 0010001, \ldots \}.
        \end{align*}

        \item If $m=3$, then
        \begin{equation*}
            \mB_j^{(3)} = \{0_r(b_1 00)(b_2 00) \ldots (b_i 00)1: 0\leq r\leq 1, b_1,\ldots,b_i\in\F_2, 0\leq i < j\},
        \end{equation*}
        and hence,
        \begin{align*}
            \mB_1^{(3)} &= \{1, 01\}, \\
            \mB_2^{(3)} &= \{1, 01, 0001, 1001, 00001, 01001\}, \\
            \mB_3^{(3)} &= \{1, 01, 0001, 1001, 00001, 01001, 0000001, 0001001, 1000001, 1001001, \ldots\}.
        \end{align*}
    \end{enumerate}
\end{example}

By applying \cref{lem:cycle_structure_comb_on_words_very_pure_gf} we obtain the generating function $\mZ^{(m)}_j(\lambda)$.

\begin{theorem}\label{thm:cycle_structure_gf_of_z_m_j_n}
    It holds that
    \begin{equation*}
        \mZ^{(m)}_j(\lambda) = \frac{\lambda \frac{\d}{\d\lambda}\mB_j^{(m)}(\lambda)}{1-\mB_j^{(m)}(\lambda)}
    \end{equation*}
    with
    \begin{equation*}
        \mB_j^{(m)}(\lambda) = \sum_{i=0}^{j-1} \sum_{r=0}^{m-2} 2^i \lambda^{im+r+1} = \lambda\cdot \frac{1-\lambda^{m-1}}{1-\lambda}\cdot\frac{1-(2\lambda^m)^j}{1-2\lambda^m}.
    \end{equation*}
\end{theorem}

\begin{proof}
    We are first going to show that $\mZ^{(m)}_j=\widetilde{\mB_j}^{(m)}$, then the first claim follows from \cref{lem:cycle_structure_comb_on_words_very_pure_gf}. Afterward we only need to determine $\mB_j^{(m)}(\lambda)$. 

    Let $z\in\mZ^{(m)}_j$ arbitrary, then by the construction of $\mB_j^{(m)}$, the anchor factorization of $z$ shows that $z\in\widetilde{\mB_j}^{(m)}$. Let now $z\in\widetilde{\mB_j}^{(m)}$, we need to show that $z\in\mZ^{(m)}_j$, i.e. we need to show that $z$ does not contain a substring $0_{m-1}\mt0_{m-1}\mt\ldots0_{m-1}\mt0_{m-1}1$ of length $jm$. Note that the occurrence of such a substring implies that there are at least two anchors of $z$ with distance at least $jm$. However, by the construction of $\mB_j^{(m)}$, every last entry of every element $b\in\mB_j^{(m)}$ is a $1$ and is an anchor in any word $x\in\widetilde{\mB_j}^{(m)}$ that contains $b$ as a factor. As every element $b\in\mB_j^{(m)}$ satisfies $\ell(b)\leq jm-1$, then all anchors of $z$ have distance at most $jm-1$, a contradiction.

    We now need to compute $\mB_j^{(m)}(\lambda)$. We know that 
    \begin{equation*}
    \mB_j^{(m)} = \{0_r b_10_{m-1}b_2 0_{m-1} b_3\ldots 0_{m-1}b_i 0_{m-1}1: 0\leq r\leq m-2, b_1,\ldots,b_i \in\F_2, i<j\}.
    \end{equation*}
    There are $2^i$ strings of the form $b_10_{m-1}b_20_{m-1}b_3\ldots0_{m-1} b_i0_{m-1}1$ of length $im+1$ with $0\leq i\leq j-1$ and therefore $2^i$ strings of the form $0_rb_10_{m-1}b_20_{m-1}b_3\ldots0_{m-1}b_i0_{m-1}1$ of length $im+1+r$ for each $0\leq r\leq m-2$. Thus,
    \begin{equation*}
        \mB_j^{(m)}(\lambda) = \sum_{i=0}^{j-1} \sum_{r=0}^{m-2} 2^i \lambda^{im+1+r} 
            = \lambda \sum_{r=0}^{m-2} \lambda^r  \sum_{i=0}^{j-1} (2\lambda^m)^i 
            = \lambda\cdot \frac{1-\lambda^{m-1}}{1-\lambda}\cdot\frac{1-(2\lambda^m)^j}{1-2\lambda^m}.
    \end{equation*}
\end{proof}

Observe that $\mZ^{(m)}(\lambda)$ is a rational generating function since $\mB_j^{(m)}(\lambda)$ is a polynomial. \cref{thm:cycle_structure_rational_gf_properties} $(b)$ implies that the sequence $(z^{(m)}(j,n))_n$ satisfies a linear recurrence. We pose the following problem.

\begin{problem}
    Is there a combinatorial interpretation of this recurrence? Does it generalize to other landscapes?
\end{problem}

\subsection{Asymptotic Results}

By \cref{thm:cycle_structure_rational_gf_properties} $(c)$ combined with \cref{lem:cycle_structure_explicit_coefficient_of_asymptotic_rational_function} we can improve the asymptotic for $z^{(m)}(j,n)$ for $n\to\infty$, first considered in \cref{lem:cycle_structure_properties_numbers_z_m_j_n}, \cref{itm:cycle_structure_z_j_n_exponential_in_n}. For this purpose we need to study the roots of $1-\mB_j^{(m)}(\lambda)$.

\begin{lemma}\label{lem:cycle_structure_solution_B_j_lambda_equal_1}
    Let $m\geq 2$. The polynomial $1-\mB_j^{(m)}(\lambda)$ has a unique positive root $\xi_{j,m}$, which is a simple root. The sequence $(\xi_{j,m})_j$ is monotone decreasing and $\xi_{j,m}\to 1/2$ for $j\to\infty$. If $\xi\in\C$ is any other root of $1-\mB_j^{(m)}(\lambda)$, then $\abs{\xi_{j,m}} < \abs{\xi}$. 
\end{lemma}

\begin{proof}
    We first show that a sequence of roots exists which converges to $1/2$. Recall 
    \begin{equation*}
        \mB_j^{(m)}(\lambda) = \sum_{i=0}^{j-1} \sum_{r=0}^{m-2} 2^i \lambda^{im+r+1} = \lambda\cdot \frac{1-\lambda^{m-1}}{1-\lambda}\cdot\frac{1-(2\lambda^m)^j}{1-2\lambda^m}.
    \end{equation*}
    First set $\lambda=1/2$, then 
    \begin{equation*}
        \mB_j^{(m)}\left(\frac{1}{2}\right) = \frac{1}{2} \cdot\frac{1-\frac{1}{2^{m-1}}}{\frac{1}{2}} \cdot \frac{1-\frac{1}{2^{(m-1)j}}}{1-\frac{1}{2^{m-1}}} = 1-\frac{1}{2^{(m-1)j}} < 1.
    \end{equation*}
    On the other hand, for large enough $j$, setting $\lambda=1/2+\varepsilon$ for $\varepsilon>0$ gives
    \begin{align*}
        \mB_j^{(m)}\left(\frac{1}{2}+\varepsilon\right) 
            &= \left(\frac{1}{2}+\varepsilon\right)\cdot \frac{1-\left(\frac{1}{2}+\varepsilon\right)^{m-1}}{\frac{1}{2}-\varepsilon}\cdot\frac{1-\left(2\left(\frac{1}{2}+\varepsilon\right)^m\right)^j}{1-2\left(\frac{1}{2}+\varepsilon\right)^m} \\
            &= \frac{\frac{1}{2}+\varepsilon}{\frac{1}{2}-\varepsilon}\cdot \frac{1-\left(\frac{1}{2}+\varepsilon\right)^{m-1}}{1-(1+2\varepsilon)\left(\frac{1}{2}+\varepsilon\right)^{m-1}}\cdot \left(1-\left(2\left(\frac{1}{2}+\varepsilon\right)^m\right)^j\right) > 1.
    \end{align*}
    By the intermediate value theorem there exists a root $\xi_{j,m}\in(1/2, 1/2+\varepsilon)$. Observe that 
    \begin{equation*}
        \frac{\d}{\d\lambda} \mB_j^{(m)}(\lambda) = \sum_{i=0}^{j-1} \sum_{r=0}^{m-2} 2^i (im+r+1)\lambda^{im+r} > 0
    \end{equation*}
    for all $\lambda>0$, implying that $\lambda\mapsto\mB_j^{(m)}(\lambda)$ is monotone increasing for $\lambda>0$ for every $j,m$. By monotonicity $\xi_{j,m}$ is the unique positive root of $1-\mB_j^{(m)}(\lambda)$ and it is a simple root. Now consider $\mB_{j+1}^{(m)}(\xi_{j,m})$. Then
    \begin{align*}
        \mB_{j+1}^{(m)}(\xi_{j,m}) 
            &= \sum_{i=0}^{j} \sum_{r=0}^{m-2} 2^i \xi_{j,m}^{im+r+1} \\
            &= \sum_{i=0}^{j-1} \sum_{r=0}^{m-2} 2^i \xi_{j,m}^{im+r+1} + \sum_{r=0}^{m-2} 2^j \xi_{j,m}^{jm+r+1} \\
            &= \mB_j^{(m)}(\xi_{j,m}) + \sum_{r=0}^{m-2} 2^j \xi_{j,m}^{jm+r+1} > 1 
    \end{align*}
    and because of the monotonicity of $\mB_j^{(m)}(\lambda)$ it follows that $\xi_{j+1,m}<\xi_{j,m}$. 

    Let $\xi\in\C$ be any other solution to $\mB_j^{(m)}(\lambda)=1$, then necessarily $\xi\notin\R$ or $\xi<0$. The triangle inequality $\abs{z_1+z_2}\leq \abs{z_1}+\abs{z_2}$ with $z_1,z_2\in\C$ holds with equality if and only if $z_1,z_2$ have the same argument. Assume that $\abs{\xi}\leq \abs{\xi_{j,m}}$, then
    \begin{align}
            \mB_j^{(m)}(\abs{\xi})
            &\leq \mB_j^{(m)}(\abs{\xi_{j,m}})
            = 1 
            = \abs{\mB_j^{(m)}(\xi)} 
            = \abs{\sum_{i=0}^{j-1} \sum_{r=0}^{m-2} 2^i \xi^{im+r+1} } \notag \\
            &\leq \sum_{i=0}^{j-1} \sum_{r=0}^{m-2} 2^i \abs{\xi}^{im+r+1} \label{eq:cycle_structure_showing_other_roots_more_abs_value} \\
            &= \mB_j^{(m)}(\abs{\xi}). \notag
    \end{align}
    Comparing left-hand side and right-hand side gives equality in \cref{eq:cycle_structure_showing_other_roots_more_abs_value}. 
    If $m>2$, then the sum 
    \begin{equation*}
        \sum_{i=0}^{j} \sum_{r=0}^{m-2} 2^i \xi^{im+r+1}=\xi+\xi^2+\ldots
    \end{equation*}
    contains two terms with different arguments, giving a contradiction. If $m=2$, then 
    \begin{equation*}
        \mB_j^{(2)} = \sum_{i=0}^{j} 2^i \xi^{2i+1}
    \end{equation*} 
    is an odd function, implying that $\xi\notin\R$ and then 
    \begin{equation*}
        \sum_{i=0}^{j} 2^i \xi^{2i+1}=\xi+2\xi^3+\ldots
    \end{equation*}
    also contains two terms with different arguments. We again obtain a contradiction to our assumption $\abs{\xi}\leq \abs{\xi_{j,m}}$, finishing the proof.
\end{proof}

\begin{corollary}
    With the notation from above,
    \begin{equation*}
        z^{(m)}(j,n) \sim \xi_{j,m}^{-n}
    \end{equation*}
    for $n\to\infty$.
\end{corollary}

\begin{proof}
    Set $P(\lambda) = \lambda \frac{\d}{\d\lambda}\mB_j^{(m)}(\lambda)$ and $Q(\lambda)=1-\mB_j^{(m)}(\lambda)$. Then by \cref{lem:cycle_structure_explicit_coefficient_of_asymptotic_rational_function} we only need to observe that
    \begin{equation*}
        -\frac{P(\xi_1)}{\xi_1Q'(\xi_1)} = 1,
    \end{equation*}
    which follows immediately from the definitions of $P$ and $Q$.
\end{proof}

\subsection{A possible application to cryptanalysis}

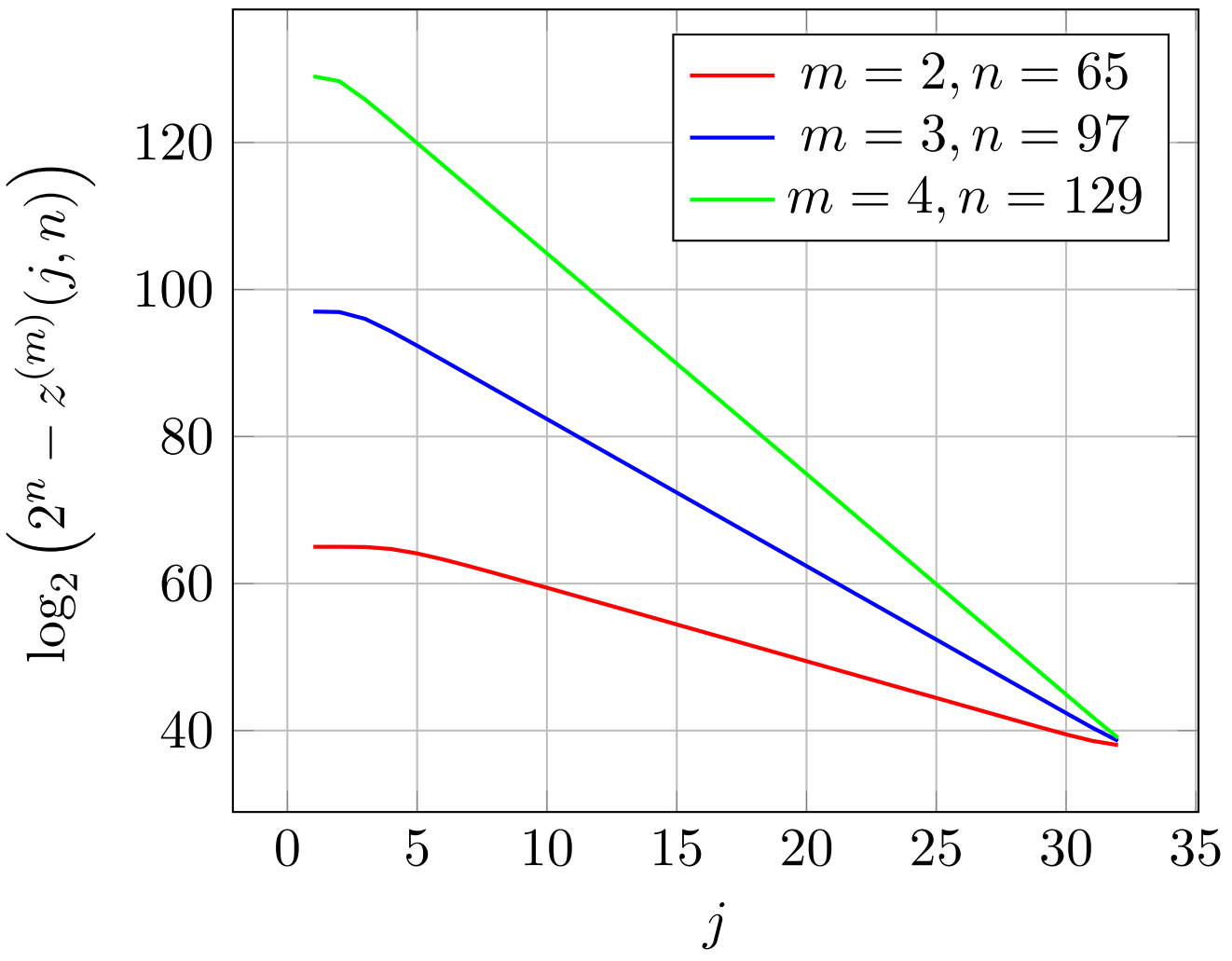
\begin{figure}
\centering
\begin{tikzpicture}
\begin{axis}[
    xlabel={$j$},
    ylabel={$\log_2\left(2^n-z^{(m)}(j,n)\right)$},
    legend pos=north east,
    grid=major,
    cycle list name=color,
]

\addplot[
    red,
    no marks, 
    thick,
    unbounded coords=discard,
    y filter/.code={
        \ifdim\pgfmathresult pt=0pt
            \def\pgfmathresult{inf}
        \fi
    }
] 
table[x=j,y={n=65},col sep=comma]{z_2_j_n_data.csv};
\addlegendentry{$m=2,n=65$}

\addplot[
    blue,
    no marks, 
    thick,
    unbounded coords=discard,
    y filter/.code={
        \ifdim\pgfmathresult pt=0pt
            \def\pgfmathresult{inf}
        \fi
    }
] 
table[x=j,y={n=97},col sep=comma]{z_3_j_n_data.csv};
\addlegendentry{$m=3,n=97$}

\addplot[
    green,
    no marks, 
    thick,
    unbounded coords=discard,
    y filter/.code={
        \ifdim\pgfmathresult pt=0pt
            \def\pgfmathresult{inf}
        \fi
    }
] 
table[x=j,y={n=129},col sep=comma]{z_4_j_n_data.csv};
\addlegendentry{$m=4,n=129$}
    
\end{axis}
\end{tikzpicture}
\caption{Logarithmic plot of $\abs{\F_2^n\setminus\mZ_{j,n}^{(m)}}=2^n -z^{(m)}(j,n)$ for $1\leq j\leq 32$ and suitable values of $m$ and $n$. }
\label{fig:cycle_structure_log_plot_z_m_j_n}
\end{figure}

We now know the behavior of $z^{(m)}(j,n)$ for $n$ fixed and $j\geq n/m$ and for $j$ fixed and $n\to\infty$. The next interesting question is how $z^{(m)}(j,n)$ behaves for $n$ fixed and $j<n/m$. We look instead at $2^n-z^{(m)}(j,n)$, i.e. the number of vectors in
\begin{equation*}
    \F_2^n\setminus\mZ_{j,n}^{(m)} = \{x\in\F_2^n: x \text{ contains } 0_{m-1}\mt0_{m-1}\mt\ldots0_{m-1}\mt0_{m-1}1 \text{ of length } jm \text{ cyclically} \}.
\end{equation*}
\cref{fig:cycle_structure_log_plot_z_m_j_n} indicates that $j\mapsto \log_2\left(2^n-z^{(m)}(j,n)\right)$ is approximately linear for $j$ neither too small nor too large. The slope appears to be $-(m-1)$, independently of $n$. Equivalently,
\begin{equation*}
    2^n - z^{(m)}(j,n) \approx 2^{m-1} \left(2^n - z^{(m)}(j+1,n)\right).
\end{equation*} 
Intuitively, this is to be expected. A vector $x\in\F_2^n\setminus \mZ_{j,n}^{(m)}$ contains the substring $0_{m-1}\mt0_{m-1}\mt\ldots0_{m-1}\mt0_{m-1}1$ of length $jm$ while $x\in\F_2^n\setminus \mZ_{j+1,n}^{(m)}$ must contain the corresponding substring of length $(j+1)m$. Thus, increasing $j$ by one roughly speaking fixes $m-1$ more positions, leading to a factor of $2^{m-1}$. 

This leads to a possible application of our results to cryptanalysis. We start by observing that the identity
\begin{equation*}
    (1+X+\cdots+X^{j-1})(1+X)=1+X^j
\end{equation*}
combined with the results of \cref{sec:pcp_studying_pcp} leads to
\begin{equation*}
    (\gamma_0+\gamma_1+\cdots+\gamma_{j-1})\circ (\gamma_0+\gamma_1) = \gamma_0+\gamma_j
\end{equation*}
for any landscape $c$. In particular,
\begin{equation}\label{eq:cycle_structure_near_inverse_identity}
    (\gamma_0+\gamma_1+\cdots+\gamma_{j-1})(\chi^{(m)}(x)) = x + \gamma_j(x)
\end{equation}
for all $x\in\F_2^n$ by choosing the landscape $c=*0_{m-1}1$. Now suppose that $\gamma_j(x)=0$, i.e. $x\in\mZ_j^{(m)}$, then 
\begin{equation*}
    (\gamma_0+\gamma_1+\cdots+\gamma_{j-1})(\chi^{(m)}(x)) = x.
\end{equation*}
This means that $f_j=\gamma_0+\gamma_1+\cdots+\gamma_{j-1}$ is a near-inverse of $\chi^{(m)}$, which may be used in a probabilistic attack.

In fact, the identity \cref{eq:cycle_structure_near_inverse_identity} for $m=2$, i.e. $\chi^{(2)}=\chi$ has appeared in an equivalent form in an earlier version of the preprint \cite{lambooij2025meet}. There it was used to approximate a single bit of $x$. Our contribution is that it is now possible to make precise statements about the probability that $f_j(\chi^{(m)}(x))=x$ for $x\in\F_2^n$. By our previous observation, increasing $j$ by $1$ roughly decreases the probability
\begin{equation*}
    \mathbb{P}\left[f_j(\chi^{(m)}(x))\neq x\right] = \mathbb{P}\left[x\in \F_2^n\setminus\mZ_{j,n}^{(m)}\right]
\end{equation*}
by a factor of $2^{m-1}$. Especially for larger $n$ this may lead to promising probabilistic attacks. We leave the explicit construction of an attack on a primitive based on $\chi^{(m)}$ as future work.

\subsection{A formula}

To aid in the construction of such an attack we present an explicit formula for $z^{(m)}(j,n)$ which is useful for running numerical experiments, especially for large $n$ and $m=2$. We start by noting
\begin{equation*}
    \mZ^{(m)}_j(\lambda) = \frac{\lambda \frac{\d}{\d\lambda}\mB_j^{(m)}(\lambda)}{1-\mB_j^{(m)}(\lambda)} 
        = \sum_{s\geq 0} \lambda \left(\mB_j^{(m)}\right)^s(\lambda) \frac{\d }{\d\lambda}\mB_j^{(m)}(\lambda)
        = \sum_{s\geq 0} \frac{\lambda}{s+1} \frac{\d }{\d\lambda}\left(\mB_j^{(m)}\right)^{s+1}(\lambda)
\end{equation*}
and now want to determine $\left(\mB_j^{(m)}\right)^{s+1}(\lambda)$. Recall for $u\in\C\setminus\{0\}$ the Taylor series
\begin{equation*}
    (1+x)^u = 1+u x + \frac{u(u-1)}{2}x^2+\frac{u(u-1)(u-2)}{3!}x^3+\ldots = \sum_{k\geq 0} \binom{u}{k} x^k,
\end{equation*}
where the binomial coefficient is $\binom{u}{k} = \frac{u(u-1)\cdots(u-k+1)}{k!}$. 
Suppose that $u\in\N$, then 
\begin{align*}
    \binom{-u}{k} &= \frac{(-u)(-u-1)\cdots(-u-k+1)}{k!} \\
        & = (-1)^k \frac{u(u+1)\cdots(u+k-1)}{k!} = (-1)^k \binom{u+k-1}{u-1}.
\end{align*}

\begin{lemma}\label{lem:cycle_structure_coefficient_delta_N_u_t}
    Let $u,t,N\in\N$. The coefficient of $\lambda^N$ in $(1+\lambda+\cdots+\lambda^t)^u$ is 
    \begin{equation*}
        \delta(N,u,t) = \sum_{k=0}^{\floor{N/(t+1)}} (-1)^k \binom{u}{k}\binom{u+N-k(t+1)-1}{u-1}.
    \end{equation*}
\end{lemma}

\begin{proof}
    We have
    \begin{align*}
        (1+\lambda+\cdots+\lambda^t)^u &= \left(\frac{1-\lambda^{t+1}}{1-\lambda}\right)^u = (1-\lambda^{t+1})^u (1-\lambda)^{-u} \\
            &= \sum_{k\geq 0} \binom{u}{k} (-1)^k (\lambda^{t+1})^k \sum_{\ell\geq 0} \binom{-u}{\ell}(-1)^\ell \lambda^\ell \\
            &= \sum_{k,\ell\geq 0} (-1)^k \binom{u}{k} \binom{u+\ell-1}{u-1} \lambda^{(t+1)k+\ell}.
    \end{align*}
    Set $N=(t+1)k+\ell$, so that $\ell=N-(t+1)k$. As $\ell\geq 0$ we get $0\leq k\leq \frac{N-\ell}{t+1}\leq \frac{N}{t+1}$. We obtain
    \begin{equation*}
        (1+\lambda+\cdots+\lambda^t)^u = \sum_{N\geq 0} \sum_{k=0}^{\floor{N/(t+1)}} (-1)^k \binom{u}{k}\binom{u+N-k(t+1)-1}{u-1} \lambda^N
    \end{equation*}
    and the claim follows. 
\end{proof}

With \cref{lem:cycle_structure_coefficient_delta_N_u_t} we now get
\begin{align*}
    \left(\mB_j^{(m)}\right)^{s+1}(\lambda) &= \left(\sum_{i=0}^{j-1} \sum_{r=0}^{m-2} 2^i \lambda^{im+r+1}\right)^{s+1} \\
        &= \lambda^{s+1} (1+\lambda+\cdots+\lambda^{m-2})^{s+1}(1+2\lambda^m+\cdots+(2\lambda^m)^{j-1})^{s+1} \\
        &= \lambda^{s+1} \sum_{k\geq 0} \delta(k,s+1,m-2)\lambda^k \sum_{\ell\geq 0} \delta(\ell,s+1,j-1)(2\lambda^m)^\ell.
\end{align*}
Thus,
\begin{align*}
    \mZ_j^{(m)}(\lambda) &= \sum_{s\geq 0} \frac{\lambda}{s+1} \frac{\d }{\d\lambda}\left(\mB_j^{(m)}\right)^{s+1}(\lambda) \\
        &= \sum_{s\geq 0} \frac{\lambda}{s+1} \frac{\d }{\d\lambda }\sum_{k,\ell\geq 0} 2^\ell \delta(k,s+1,m-2)\delta(\ell,s+1,j-1)\lambda^{\ell m+k+s+1} \\
        &= \sum_{s,k,\ell\geq 0} \frac{2^\ell(\ell m+k+s+1)}{s+1}  \delta(k,s+1,m-2)\delta(\ell,s+1,j-1)\lambda^{\ell m+k+s+1}.
\end{align*}
Setting $n=\ell m+k+s+1$ gives
\begin{equation*}
    \mZ_j^{(m)}(\lambda) = \sum_{n\geq 0} \sum_{k,\ell\geq 0} \frac{n2^\ell }{n-\ell m-k} \delta(k,n-\ell m-k,m-2)\delta(\ell,n-\ell m-k,j-1) \lambda^n
\end{equation*}
and therefore the following lemma.

\begin{lemma}\label{lem:cycle_structure_formula_z_m_j_n}
    We have 
    \begin{equation*}
        z^{(m)}(j,n) = \sum_{\substack{k,\ell\geq 0 \\ \ell m+k<n}} \frac{n2^\ell }{n-\ell m-k} \delta(k,n-\ell m-k,m-2)\delta(\ell,n-\ell m-k,j-1).
    \end{equation*}
\end{lemma}

\begin{corollary}
    For $m=2$ we get
    \begin{equation*}
        z^{(2)}(j,n) = \sum_{\ell=0}^{\floor{n\cdot\frac{j-1}{2j-1}}} \frac{n 2^\ell }{n-2\ell} \sum_{h=0}^{\floor{\ell/j}} (-1)^h \binom{n-2\ell}{h}\binom{n-\ell-hj-1}{n-2\ell-1}.
    \end{equation*}
\end{corollary}

\begin{proof}
    We start with \cref{lem:cycle_structure_formula_z_m_j_n} and use that $\delta(N,u,t)=0$ for all $N>ut$. In particular, $\delta(k,n-\ell m-k,m-2)=\delta(k,n-2\ell-k,0)=0$ for all $k>0$, so that we only have nonzero terms in the sum if $k=0$, in which case $\delta(0,n-2\ell,0)=1$. Plugging this in we obtain
    \begin{equation*}
        z^{(2)}(j,n) = \sum_{\ell\geq 0} \frac{n2^\ell }{n-2\ell} \delta(\ell,n-2\ell,j-1).
    \end{equation*}
    Similarly, $\delta(\ell,n-2\ell,j-1)=0$ for $\ell>(n-2\ell)(j-1)$ which is equivalent to $\ell > n\cdot\frac{j-1}{2j-1}$, giving
    \begin{equation*}
        z^{(2)}(j,n) = \sum_{\ell=0}^{\floor{n\cdot\frac{j-1}{2j-1}}} \frac{n2^\ell }{n-2\ell} \delta(\ell,n-2\ell,j-1).
    \end{equation*}
    Using \cref{lem:cycle_structure_coefficient_delta_N_u_t} we obtain the claim. 
\end{proof}

We conclude this section with two open problems. For the first, we have obtained the generating function
\begin{equation*}
    \sum_{n\geq 0} z^{(m)}(j,n)\lambda^n = \frac{\lambda \frac{\d}{\d\lambda}\mB_j^{(m)}(\lambda)}{1-\mB_j^{(m)}(\lambda)}
\end{equation*}
for a polynomial $\mB_j^{(m)}(\lambda)$. It is natural to ask whether this result can be extended to generating functions in the parameter $j$, or simultaneously in both $j$ and $n$. This leads to the following open problem.

\begin{problem}
    Determine whether the generating functions
    \begin{equation*}
        \sum_{j\geq 0} z^{(m)}(j,n)\mu^j
    \end{equation*}
    and
    \begin{equation*}
        \sum_{n,j\geq 0} z^{(m)}(j,n)\lambda^n\mu^j
    \end{equation*}
    admit closed-form expressions, and, if so, determine such expressions.
\end{problem}

The second open problem is to determine the cycle structure of other shift-invariant maps $f\in\mG$.

\begin{problem}
    Determine the cycle structure of other shift-invariant maps $f\in\mG$ for a given quasi-conserved landscape $c$. 
\end{problem}

We expect that the methods developed here can be adapted to study the cycle structure of other maps arising from quasi-conserved landscapes. A natural first class to consider is that of landscapes given by bifix-free strings; see \cref{rem:right_sided_bifix_free}. Numerical experiments suggest that the cycle structure may in fact be independent of the choice of the bifix-free string. More precisely, they suggest the following conjecture.

\begin{conjecture}
    The cycle structure of $\id+*a_1\ldots a_m$ is identical to the cycle structure of $\chi^{(m)}$ for any bifix-free string $a_1\ldots a_m$. 
\end{conjecture}

\bibliographystyle{splncs04_doi.bst}
\bibliography{arxiv_refs.bib}

@InProceedings{kriepke2024algebraic,
  author={Kriepke, Björn and Kyureghyan, Gohar},
  title={{Algebraic} {Structure} of the {Iterates} of $\chi$},
  booktitle={Advances in Cryptology -- CRYPTO 2024},
  year={2024},
  publisher={Springer},
  volume={14923},
  pages={412--424},
  doi={10.1007/978-3-031-68385-5_13},
  bibchecked = {21.09.2026}
}

@article{kriepke2026siblings,
  title     = {There are siblings of $\chi$ which are permutations for $n$ even},
  author    = {Kriepke, Björn and Kyureghyan, Gohar},
  journal   = {Designs, Codes and Cryptography},
  publisher = {Springer},
  volume    = {94},
  number    =  {6},
  eid = {125},
  year      =  {2026},
  doi = {10.1007/s10623-026-01830-0},
  bibchecked = {21.09.2026}
}

@article{lyu2026generalized,
  title={A {Generalized} $\chi_n$-{Function}},
  author={Lyu, Cheng and Yuan, Mu and Zheng, Dabin and Sun, Siwei and Li, Shun},
  journal={IEEE Transactions on Information Theory},
  year={2026},
  publisher={IEEE},
  volume={72},
  number={4},
  pages={2543--2556},
  doi = {10.1109/TIT.2026.3662335},
  bibchecked = {21.09.2026}
}

@phdthesis{daemen1995cipher,
  title={{Cipher} and {Hash} {Function} {Design} {Strategies} based on linear and differential cryptanalysis},
  author={Daemen, Joan},
  year={1995},
  school={KU Leuven},
  bibchecked = {21.09.2026}
}

@article{toffoli1990invertible,
  title={{Invertible} cellular automata: {A} review},
  author={Toffoli, Tommaso and Margolus, Norman H},
  journal={Physica D: Nonlinear Phenomena},
  volume={45},
  number={1--3},
  pages={229--253},
  year={1990},
  publisher={Elsevier},
  doi={10.1016/0167-2789(90)90185-R},
  bibchecked = {21.09.2026}
}

@article{liu2026finding,
  author    = {Liu, Fukang and Dixit, Vaibhav and Sarkar, Santanu and Meier, Willi and Isobe, Takanori},
  title     = {{Finding} the {Inverse} of some {Shift} {Invariant} {Transformations}},
  journal   = {Journal of Cryptology},
  volume    = {39},
  number    = {1},
  eid       = {1},
  year      = {2026},
  doi       = {10.1007/s00145-025-09559-0},
  publisher = {Springer},
  bibchecked = {21.09.2026}
}

@article{liu2022inverse,
  title={The {Inverse} of $\chi$ and {Its} {Applications} to {Rasta}-{Like} {Ciphers}},
  author={Liu, Fukang and Sarkar, Santanu and Meier, Willi and Isobe, Takanori},
  journal={Journal of Cryptology},
  volume={35},
  number={4},
  pages={28},
  year={2022},
  doi={10.1007/s00145-022-09439-x},
  publisher={Springer},
  bibchecked = {21.09.2026}
}

@article{schoone2024algebraic,
  title={Algebraic properties of the maps $\chi_n$},
  author={Schoone, Jan and Daemen, Joan},
  journal={Designs, Codes and Cryptography},
  volume={92},
  number={8},
  pages={2341--2365},
  year={2024},
  doi={10.1007/s10623-024-01395-w},
  publisher={Springer},
  bibchecked = {21.09.2026}
}

@article{graner2025bijectivity,
  title={On the bijectivity of the map $\chi$},
  author={Graner, Anna-Maurin and Kriepke, Björn and Krompholz, Lucas and Kyureghyan, Gohar},
  journal={Cryptography and Communications},
  volume={18},
  number={1},
  pages={97--105},
  year={2026},
  doi={10.1007/s12095-025-00786-7},
  publisher={Springer},
  bibchecked = {21.09.2026}
}

@article{schoone2024state,
  title={The state diagram of $\chi$},
  author={Schoone, Jan and Daemen, Joan},
  journal={Designs, Codes and Cryptography},
  volume={92},
  number={5},
  pages={1393--1421},
  year={2024},
  doi={10.1007/s10623-023-01349-8},
  publisher={Springer},
  bibchecked = {21.09.2026}
}

@techreport{nistSHA3StandardPermutationBased2015,
  title = {{SHA-3} {Standard}: {Permutation}-{Based} {Hash} and {Extendable}-{Output} {Functions}},
  author = {{NIST}},
  year = {2015},
  number = {FIPS 202},
  institution = {National Institute of Standards and Technology},
  doi = {10.6028/NIST.FIPS.202},
  bibchecked = {21.09.2026}
}

@article{dobraunig2021ascon,
  title={\textsc{Ascon} v1.2: {Lightweight} {Authenticated} {Encryption} and {Hashing}},
  author={Dobraunig, Christoph and Eichlseder, Maria and Mendel, Florian and Schläffer, Martin},
  journal={Journal of Cryptology},
  volume={34},
  number={33},
  pages={1--42},
  year={2021},
  doi={10.1007/s00145-021-09398-9},
  publisher={Springer},
  bibchecked = {21.09.2026}
}

@book{Stanley2011enumerativeVol1, 
  title={Enumerative {Combinatorics}}, 
  author={Stanley, Richard P.}, 
  volume = {1},
  edition = {2},
  publisher={Cambridge University Press}, 
  year={2011}, 
  doi={10.1017/CBO9781139058520},
  bibchecked = {21.09.2026}
}

@book{wilf2006generatingfunctionology,
  author    = {Wilf, Herbert S.},
  title     = {generatingfunctionology},
  edition   = {3},
  publisher = {A K Peters},
  year      = {2006},
  doi       = {10.1201/b10576},
  bibchecked = {21.09.2026}
}

@article{haugland2025shift,
  title={Shift-invariant transformations and almost liftings},
  author={Haugland, Jan Kristian and Omland, Tron},
  journal={Cryptography and Communications},
  volume={18},
  number={3},
  pages={705--726},
  year={2026},
  publisher={Springer},
  doi = {10.1007/s12095-025-00848-w},
  bibchecked = {21.09.2026}
}

@article{omland2022permutation,
  title={Permutation rotation-symmetric {S-boxes}, liftings and affine equivalence},
  author={Omland, Tron and Stănică, Pantelimon},
  journal={arXiv preprint arXiv:2203.00778},
  year={2022},
  doi={10.48550/arXiv.2203.00778},
  bibchecked = {21.09.2026}
}

@book{Lidl:1997,
  title = {{Finite} {Fields}},
  author = {Lidl, Rudolf and Niederreiter, Harald},
  publisher = {Cambridge University Press},
  year = {2009},
  edition = {2},
  doi={10.1017/CBO9780511525926}
}

@article{amoroso1972decision,
  title={Decision procedures for surjectivity and injectivity of parallel maps for tessellation structures},
  author={Amoroso, Serafino and Patt, Yale N.},
  journal={Journal of Computer and System Sciences},
  volume={6},
  number={5},
  pages={448--464},
  year={1972},
  publisher={Academic Press},
  doi={10.1016/S0022-0000(72)80013-8},
  bibchecked ={21.09.2026}
}

@article{haugland2026new,
  title={New classes of reversible cellular automata},
  author={Haugland, Jan Kristian and Omland, Tron},
  journal={Designs, Codes and Cryptography},
  volume={94},
  number={6},
  eid={131},
  year={2026},
  publisher={Springer},
  doi={10.1007/s10623-026-01869-z},
  bibchecked = {21.09.2026}
}

@misc{lambooij2025meet,
  title = {Meet-in-the-{Middle} {Attacks} on {Full} {ChiLow}-32},
  author = {Eran Lambooij and Patrick Neumann and Michiel Verbauwhede and Shichang Wang and Tianyu Zhang},
  howpublished = {Cryptology {ePrint} Archive, Paper 2025/1601},
  year = {2025},
  url = {https://eprint.iacr.org/archive/2025/1601/20250905:183726},
  bibchecked = {21.09.2026}
}

@article{nielsen1973note,
  title={A {Note} on {Bifix}-{Free} {Sequences}},
  author={Nielsen, Peter Tolstrup},
  journal={IEEE Transactions on Information Theory},
  volume={19},
  number={5},
  pages={704--706},
  year={1973},
  publisher={IEEE},
  doi = {10.1109/TIT.1973.1055065},
  bibchecked = {21.09.2026}
}

@misc{zhao2026bifix,
  author = {Xiao-Xin Zhao and Deng Tang and Zhong-Xiao Wang},
  title = {A family of invertible shift-invariant maps with strong arithmetic properties},
  howpublished = {Cryptology {ePrint} Archive, Paper 2026/1249},
  year = {2026},
  url = {https://eprint.iacr.org/2026/1249},
  bibchecked = {21.09.2026}
}

@book{flajoletAnalyticCombinatorics2009,
  title = {Analytic Combinatorics},
  author = {Flajolet, Philippe and Sedgewick, Robert},
  year = {2011},
  publisher = {Cambridge University Press},
  doi = {10.1017/CBO9780511801655},
}

@article{grassi2024monolith,
  author    = {Grassi, Lorenzo and Khovratovich, Dmitry and Lüftenegger, Reinhard and Rechberger, Christian and Schofnegger, Markus and Walch, Roman},
  title     = {{Monolith}: {Circuit}-{Friendly} {Hash} {Functions} with {New} {Nonlinear} {Layers} for {Fast} and {Constant}-{Time} {Implementations}},
  journal   = {IACR Transactions on Symmetric Cryptology},
  volume    = {2024},
  number    = {3},
  pages     = {44--83},
  year      = {2024},
  doi       = {10.46586/tosc.v2024.i3.44-83}
}

@article{kari2005CA,
  title = {{Theory} of {Cellular} {Automata}: {A} {Survey}},
  shorttitle = {Theory of Cellular Automata},
  author = {Kari, Jarkko},
  year = {2005},
  journal = {Theoretical Computer Science},
  volume = {334},
  number = {1--3},
  pages = {3--33},
  doi = {10.1016/j.tcs.2004.11.021},
  bibchecked = {21.09.2026}
}

@phdthesis{graner2024irreducible,
  title={Irreducible polynomials over finite fields for {coding} and cryptography},
  author={Graner, Anna-Maurin},
  year={2024},
  school={University of Rostock},
  doi = {10.18453/rosdok_id00004733},
  bibchecked = {21.09.2026}
}

@article{ahmad1969monomial,
  title={Cycle {Structure} of {Automorphisms} of {Finite} {Cyclic} {Groups}},
  author={Ahmad, Shair},
  journal={Journal of Combinatorial Theory},
  volume={6},
  number={4},
  pages={370--374},
  year={1969},
  publisher={Elsevier},
  doi = {10.1016/S0021-9800(69)80032-3}
}

@article{rudolf1991dickson,
  title={Cycle structure of {Dickson} permutation polynomials},
  author={Lidl, Rudolf and Mullen, Gary L.},
  journal={Mathematical Journal of Okayama University},
  volume={33},
  number={1},
  pages={1--11},
  year={1991},
  publisher={Okayama University},
  doi = {10.18926/mjou/33705},
  bibchecked = {21.09.2026}
}

@article{mullen1988linear,
  title = {{Cycles} of {Linear} {Permutations} over a {Finite} {Field}},
  author={Mullen, Gary L. and Vaughan, Theresa P.},
  journal={Linear Algebra and Its Applications},
  volume={108},
  number={C},
  pages={63--82},
  year={1988},
  publisher={Elsevier},
  doi = {10.1016/0024-3795(88)90179-6},
  bibchecked = {21.09.2026}
}

@article{qureshi2015redei,
  title = {{Rédei} {Actions} on {Finite} {Fields} and {Multiplication} {Map} in {Cyclic} {Group}},  author={Qureshi, Claudio and Panario, Daniel},
  journal={SIAM Journal on Discrete Mathematics},
  volume={29},
  number={3},
  pages={1486--1503},
  year={2015},
  publisher={SIAM},
  doi = {10.1137/140993338},
  bibchecked = {21.09.2026}
}

@phdthesis{gerike2020determining,
  title = {{Determining} the {Cycle} {Structure} of {Permutation} {Polynomials} of {Shape} $X^t+\gamma\Tr(X^k)$},
  author = {Gerike, Daniel},
  year = {2020},
  school = {Otto-von-Guericke-Universität Magdeburg},
  pages = {94},
  doi = {10.25673/33931},
  bibchecked = {21.09.2026}
}

@book{bors2025cyclotomic,
  title = {{Functional} {Graphs} of {Generalized} {Cyclotomic} {Mappings} of {Finite} {Fields}},
  author={Bors, Alexander and Panario, Daniel and Wang, Qiang},
  year={2025},
  publisher={European Mathematical Society-EMS-Publishing House GmbH},
  doi = {10.4171/MEMS/23},
  bibchecked = {21.09.2026}
}

@incollection{martins2019survey,
  title={A survey on iterations of mappings over finite fields},
  author={Martins, Rodrigo and Panario, Daniel and Qureshi, Claudio},
  booktitle={Combinatorics and Finite Fields},
  pages={135--172},
  year={2019},
  publisher={De Gruyter},
  doi = {10.1515/9783110642094-008},
  bibchecked = {21.09.2026}
}

@misc{otal2023conjecture,
  title = {A {Solution} to a {Conjecture} on the {Maps} $\chi_n^{(k)}$},
  author = {Otal, Kamil},
  howpublished = {Cryptology {ePrint} Archive, Paper 2023/1782},
  year = {2023},
  url = {https://eprint.iacr.org/2023/1782},
  bibchecked = {21.09.2026}
}

@misc{kriepke2026wcc,
  author = {Kriepke, Björn},
  title = {Shift-invariant functions that behave like polynomials},
  year = {2026},
  howpublished = {Extended abstract, 14th International Workshop on Coding and Cryptography (WCC 2026)},
  url = {https://wcc2026.inria.fr/assets/final_versions/WCC2026_paper_48.pdf},
  bibchecked = {21.09.2026}
}

@misc{kriepke2026pott65,
  title = {Shift-invariant functions that behave like polynomials},
  author = {Kriepke, Björn},
  year = {2026},
  howpublished = {Slides, Finite Geometry, Combinatorics, Boolean Functions and their Links},
  url = {https://cloud.ovgu.de/s/97KmQnKZTHDiZz3},
  bibchecked = {13.10.2026}
}

@misc{krompholz2026cmtriples,
  author = {Krompholz, Lucas},
  title = {Construction of {CM}-triples from a new class of permutations},
  year = {2026},
  howpublished = {Extended abstract, Finite Geometry, Combinatorics, Boolean Functions and their Links},
  url = {https://cloud.ovgu.de/s/PNS8SWq5FZwREWP},
  bibchecked = {13.10.2026}
}

@article{feng2026generalization,
  title={A generalization of the map $\chi$},
  author={Feng, Xiutao and Wang, Qiang and Yu, Jingyi and Zhang, Anpeng},
  journal={arXiv preprint arXiv:2609.19548},
  year={2026},
  doi = {10.48550/arXiv.2609.19548},
  bibchecked = {21.09.2026}
}

@article{lyu2026classshiftinvariant,
  title={A {Class} of {Shift}-{Invariant} {Permutations} and {Their} {Algebraic} {Structure}}, 
  author={Cheng Lyu and Dabin Zheng and Mu Yuan and Lei Hu},
  journal={arXiv preprint arXiv:2609.31174},
  year={2026},
  doi={10.48550/arXiv.2609.31174}, 
  bibchecked = {30.09.2026}
}

@inproceedings{belkheyar2025chilow,
  title={\textsc{ChiLow} and \textsc{ChiChi}: {New} {Constructions} for {Code} {Encryption}},
  author={Belkheyar, Yanis and Derbez, Patrick and Ghosh, Shibam and Leander, Gregor and Mella, Silvia and Perrin, Léo and Rasoolzadeh, Shahram and Stennes, Lukas and Sun, Siwei and Van Assche, Gilles and Vizár, Damian},
  booktitle = {Advances in Cryptology -- EUROCRYPT 2025},
  pages = {212--243},
  year = {2025},
  publisher = {Springer},
  doi = {10.1007/978-3-031-91107-1_8},
  bibchecked = {21.09.2026}
}

@article{andreoli2025generalizations,
  title={{Generalizations} of {ChiChi}: {Families} of {Low-Latency} {Permutations} in {Any} {Even} {Dimension}},
  author={Andreoli, Samuele and Leander, Gregor and Piccione, Enrico and Stennes, Lukas},
  journal={IACR Transactions on Symmetric Cryptology},
  volume={2025},
  number={3},
  pages={800--826},
  year={2025},
  doi = {10.46586/tosc.v2025.i3.800-826},
  bibchecked = {21.09.2026}
}

@inproceedings{dobraunig2018rasta,
  title={{Rasta}: {A} {Cipher} with {Low} {ANDdepth} and {Few} {ANDs} per {Bit}},
  author={Dobraunig, Christoph and Eichlseder, Maria and Grassi, Lorenzo and Lallemand, Virginie and Leander, Gregor and List, Eik and Mendel, Florian and Rechberger, Christian},
  booktitle = {Advances in Cryptology -- CRYPTO 2018},
  pages={662--692},
  year={2018},
  publisher={Springer},
  doi = {10.1007/978-3-319-96884-1_22},
  bibchecked = {21.09.2026}
}

@article{daemen2020subterranean,
  title   = {The {Subterranean} 2.0 {Cipher} {Suite}},
  author  = {Daemen, Joan and Massolino, Pedro Maat Costa and Mehrdad, Alireza and Rotella, Yann},
  journal = {IACR Transactions on Symmetric Cryptology},
  year    = {2020},
  volume  = {2020},
  number  = {S1},
  pages   = {262--294},
  doi     = {10.13154/tosc.v2020.iS1.262-294}
}

\end{document}